# GOAL-ORIENTED ADAPTIVE MESH REFINEMENT FOR NON-SYMMETRIC FUNCTIONAL SETTINGS

BRENDAN KEITH, ALI VAZIRI ASTANEH, AND LESZEK DEMKOWICZ

ABSTRACT. In this article, a new unified duality theory is developed for Petrov–Galerkin finite element methods. This novel theory is then used to motivate goal-oriented adaptive mesh refinement strategies for use with discontinuous Petrov–Galerkin (DPG) methods. The focus of this article is mainly on broken ultraweak variational formulations of stationary boundary value problems, however, many of the ideas presented within are general enough that they be extended to any such well-posed variational formulation. The proposed goal-oriented adaptive mesh refinement procedures require the construction of refinement indicators for both a primal problem and a dual problem. In the DPG context, the primal problem is simply the system of linear equations coming from a standard DPG method and the dual problem is a similar system of equations, coming from *a new method* which is dual to DPG. This new method has the same coefficient matrix as the associated DPG method but has a different load. We refer to this new finite element method as a DPG* method. A thorough analysis of DPG* methods, as stand-alone finite element methods, is not given here but will be provided in subsequent articles. For DPG methods, the current theory of *a posteriori* error estimation is reviewed and the reliability estimate in [13, Theorem 2.1] is improved on. For DPG* methods, three different classes of refinement indicators are derived and several contributions are made towards rigorous *a posteriori* error estimation. At the closure of the article, results of numerical experiments with Poisson's boundary value problem in a three-dimensional domain are provided. These results clearly demonstrate the utility of the goal-oriented adaptive mesh refinement strategies for quantities of interest with either interior or boundary terms.

## 1. INTRODUCTION

With a natural capacity to accommodate problems with non-symmetric functional settings, the discontinuous Petrov–Galerkin (DPG) methodology [24, 27] has proven its merit in many areas of engineering interest [38, 40, 43, 44, 51, 63, 68, 69]. Desirable features of DPG methods include intrinsic numerical stability [16, 26, 30], a positive definite stiffness matrix for all well-posed boundary value problems [52], flexibility of variational formulation for the treatment of said boundary value problems [14, 40, 50], and a "built-in" *a posteriori* error estimator which can be used for adaptive mesh refinement [13, 29]. The methodology incorporates a user-defined test norm, which, if chosen well, allows DPG methods to deliver a nearly optimal projection of the solution variable in a desired norm [9, 42, 75]. By exploiting properties of the underlying variational formulation, some DPG methods (ultraweak DPG methods) also permit irregular polytopal meshes [73].

For steady boundary value problems, all DPG methods begin with the declaration of a minimum residual principle. It is the first-order optimality condition for this minimum principle which delivers the variational equation determining the ideal DPG solution [52]. One advantage of this process is that it induces a very robust framework for analysis [12–14, 46]. From the implementation point of view, this universal methodology also has many significant advantages which has led to very general and very user-friendly DPG software [65–67].

Adaptive mesh refinement (AMR) and *a posteriori* error estimation [1, 3, 15, 33, 35, 58] has been fundamental since the inception of DPG methods and so AMR has been incorporated into most applied research involving DPG [9, 17, 29, 34, 50, 65]. In each of these studies, the same "built-in" *a posteriori* error estimator was used to

**Acknowledgments.** The work of the first author was supported in part by NSF grant DMS-1418822 and the 2017 Graduate School University Graduate Continuing Fellowship at The University of Texas at Austin. The work of the second author was supported in part by the 2016 Peter O'Donnell, Jr. Postdoctoral Fellowship from The Institute for Computational Engineering and Sciences at The University of Texas at Austin. This work of the third author was supported in part from grants by NSF (DMS-1418822), AFOSR (FA9550-12-1-0484), and ONR (N00014-15-1-2496). All authors express their gratitude to Socratis Petrides for his assistance with certain features of the $hp$3D software during the numerical experiments.

drive AMR. This error estimator solely relies upon an intrinsic so-called "energy norm" estimate of the solution error. This energy norm is induced by the chosen underlying minimum residual principle, and cannot otherwise be influenced by any specified output of the simulation. In the context of AMR, this is an example of what we will later on refer to as solution-adaptive mesh refinement (SMR).

Although SMR is sufficient in many circumstances, AMR geared towards a specified quantity of interest (QoI)—an output of the simulation which is defined through a function of the solution variables—instead of a global intrinsic energy, can yield far more computationally efficient results. In most practical circumstances, a simulation is invoked to yield an estimate of a particular (or various) QoI. As a motivating principle, the effect of solution error on any such estimate will depend upon the QoI under consideration. Therefore, an intelligent AMR strategy should take the influence of the given QoI into account; therein possibly sacrificing a minimal global solution error (energy error) for a far lesser QoI error. In many circumstances, if it is done well, an extrinsically-motivated AMR strategy can significantly reduce the ratio of QoI error to computational expenditure. We will broadly refer to this form of AMR as goal-oriented adaptive mesh refinement (GMR).

In the GMR framework developed in this manuscript, a primal and a dual problem is solved and refinement indicators for both problems are used to mark elements for refinement. In the DPG context, the associated dual problem can be identified with ***a new class of methods*** which we call DPG* methods [49] (see Section 4).

## 2. Overview

Aspects of goal-oriented *a posteriori* error estimation and GMR have been the subject of significant practical interest in computational science for a decades [2, 7, 45, 53, 55, 61, 64] and so a significant amount of its theory has been developed for a large variety of specific applications (see [54] and references therein). The lack of any sophisticated GMR framework for DPG methods was only recently identified as a shortcoming in [51]. Accordingly, the principal motivation of this article was to address the present DPG shortcoming in GMR. However, along the way, we uncovered a rich duality theory which has begun to take on a life of its own.

**Main contributions and outline.** The remainder of this article is organized as follows:

- In Section 3, the general Petrov–Galerkin duality theory is outlined in an idealized setting. Here, important notions based on optimality are introduced along with a review of relevant material from the literature on test norms and variational formulations.
- In Section 4, the discrete setting for the primal and dual problems and a fundamental duality theorem (Theorem 4.2) is presented. It is in this section that DPG* methods are first introduced and the QoI error is identified with norms of specific residuals of the DPG and DPG* solutions.
- In Section 5, our GMR strategy is presented.
- In Section 6, conventional energy-based DPG *a posteriori* error estimation is presented. Here, a new version of [13, Theorem 2.1] is given.
- In Section 7, three different *a posteriori* error estimation strategies for DPG* methods are proposed.
- Finally, Section 8 presents numerical experiments with a wide variety of possible quantities of interest.

**Notation.** Throughout this document, $\Omega \subseteq \mathbb{R}^d$ will denote a bounded connected Lipschitz domain and $\mathcal{T}$, referred to as the mesh, will be a finite open disjoint partition of $\Omega$ into Lipschitz subdomains $K \in \mathcal{T}$. Specifically, $\mathcal{T}$ is collection of open subsets $K \subseteq \Omega$, $|\mathcal{T}| < \infty$, $\bigcup_{K\in\mathcal{T}} \overline{K} = \overline{\Omega}$, and $K \cap \widetilde{K} = \varnothing$, for all $K \neq \widetilde{K} \in \mathcal{T}$. Here, each $K \in \mathcal{T}$, referred to as an element, is necessarily Lipschitz.

Occasionally, we will write $A \lesssim B$ or $B \lesssim A$ if there exists a constant $C > 0$, independent of the maximal mesh size, $h$, or the individual element size, $h_K$, such that $A \leq CB$ or $B \leq CA$, respectively. Similarly, $A \eqsim B$ is understood to mean that $A \lesssim B$ and $B \lesssim A$.

We will also often deal with bounded linear operators $\mathcal{L} \in B(\mathcal{X}, \mathcal{Y})$ between Banach spaces $\mathcal{X}$ and $\mathcal{Y}$. From now on, when dealing with quotients, such as in the definition of the norm $\|\mathcal{L}\|_{B(\mathcal{X},\mathcal{Y})} = \sup_{\boldsymbol{x}\in\mathcal{X}\setminus\{0\}} \frac{\|\mathcal{L}\boldsymbol{x}\|_{\mathcal{Y}}}{\|\boldsymbol{x}\|_{\mathcal{X}}}$, we will assume that the infinima and suprema ignore zero, wherein such quotients are not defined. Furthermore, because the meaning can often be understood by context in the case of such linear operators, instead of writing out $\|\mathcal{L}\|_{B(\mathcal{X},\mathcal{Y})}$ in totality, at times we will suppress the subscript above and simply write $\|\mathcal{L}\|$.

## 3. PRELIMINARY THEORY

3.1. **The Riesz operator.** Let $\mathcal{W}$ be a Hilbert space over $\mathbb{R}$. We will denote the norm on $\mathcal{W}$ by $\|\cdot\|_{\mathcal{W}} = (\cdot,\cdot)_{\mathcal{W}}^{1/2}$, where $(\cdot,\cdot)_{\mathcal{W}}$ is the associated inner product. In the case $\mathcal{W} = L^2(K)$, where $K \subseteq \overline{\Omega}$, we will simply just write $\|\cdot\|_K = (\cdot,\cdot)_K^{1/2} = \big(\int_K \cdot^2\big)^{1/2}$.

Denote the duality pairing between $\mathcal{W}$ and $\mathcal{W}'$ by $\langle\cdot,\cdot\rangle_{\mathcal{W}',\mathcal{W}}$. That is, for all bounded linear functionals $F \in \mathcal{W}'$ and all $w \in \mathcal{W}$, $\langle F, w\rangle_{\mathcal{W}',\mathcal{W}} = F(w)$. Now, recalling the Riesz representation theorem [19], the topological dual of $\mathcal{W}$, denoted $\mathcal{W}'$, is isometrically isomorphic to $\mathcal{W}$. Moreover, this isomorphism $\mathcal{R}_{\mathcal{W}} : \mathcal{W} \to \mathcal{W}'$ can be defined explicitly:

$$\langle \mathcal{R}_{\mathcal{W}}\, w, \widetilde{w}\rangle_{\mathcal{W}',\mathcal{W}} = \langle w, \mathcal{R}_{\mathcal{W}}\, \widetilde{w}\rangle_{\mathcal{W},\mathcal{W}'} = (w,\widetilde{w})_{\mathcal{W}} = (\mathcal{R}_{\mathcal{W}}\, w, \mathcal{R}_{\mathcal{W}}\, \widetilde{w})_{\mathcal{W}'}\,, \quad \forall\, w, \widetilde{w} \in \mathcal{W}\,.$$

We call $\mathcal{R}_{\mathcal{W}} : \mathcal{W} \to \mathcal{W}'$, defined above, the *Riesz operator*. Note that $\|\,\mathcal{R}_{\mathcal{W}}\, w\|_{\mathcal{W}'} = \|w\|_{\mathcal{W}}$, for all $w \in \mathcal{W}$, $\|F\|_{\mathcal{W}'} = \|\,\mathcal{R}_{\mathcal{W}}^{-1}\, F\|_{\mathcal{W}}$, for all $F \in \mathcal{W}'$, $\mathcal{R}_{\mathcal{W}} = \mathcal{R}_{\mathcal{W}}'$, and $\mathcal{R}_{\mathcal{W}'} = \mathcal{R}_{\mathcal{W}}^{-1}$ (under the identification $\mathcal{W} \sim \mathcal{W}''$). From now on, we will not explicitly declare the spaces in a duality pairing, $\langle\cdot,\cdot\rangle_{\mathcal{W}',\mathcal{W}}$, because they can be deduced from context.

With the Riesz operator now introduced, it is appropriate to immediately present a key lemma which will get significant mileage throughout this article. Our proof relies upon the key properties established above.

**Lemma 3.1.** *Let $F \in \mathcal{W}'$, where $\mathcal{W}$ is Hilbert, and let $\mathcal{W}_0 \subseteq \mathcal{W}$ be a closed subspace. Define $\mathcal{W}_1 = \mathcal{W}_0^{\perp}$. Then*

$$\sup_{w\in\mathcal{W}} \frac{|F(w)|^2}{\|w\|_{\mathcal{W}}^2} = \sup_{w_0\in\mathcal{W}_0} \frac{|F(w_0)|^2}{\|w_0\|_{\mathcal{W}}^2} + \sup_{w_1\in\mathcal{W}_1} \frac{|F(w_1)|^2}{\|w_1\|_{\mathcal{W}}^2}\,.$$

*Proof.* Let $f = \mathcal{R}_{\mathcal{W}}^{-1}\, F$ and so $(f,w)_{\mathcal{W}} = F(w)$, for all $w \in \mathcal{W}$. Moreover, $\|f\|_{\mathcal{W}} = \|F\|_{\mathcal{W}'}$. If we orthogonally decompose $f = f_0 + f_1$, where $f_0 \in \mathcal{W}_0$ and $f_1 \in \mathcal{W}_1$, then, by orthogonality,

$$F(w_0) = (f_0, w_0)_{\mathcal{W}}\,, \quad \forall\, w_0 \in \mathcal{W}_0\,, \qquad \text{and} \qquad F(w_1) = (f_1, w_1)_{\mathcal{W}}\,, \quad \forall\, w_1 \in \mathcal{W}_1\,.$$

Therefore, $f_0 = \mathcal{R}_{\mathcal{W}_0}^{-1}(F|_{\mathcal{W}_0})$ and $f_1 = \mathcal{R}_{\mathcal{W}_1}^{-1}(F|_{\mathcal{W}_1})$. Finally,

$$\|F\|_{\mathcal{W}'}^2 = \|f\|_{\mathcal{W}}^2 = \|f_0\|_{\mathcal{W}_0}^2 + \|f_1\|_{\mathcal{W}_1}^2 = \|F|_{\mathcal{W}_0}\|_{\mathcal{W}_0'}^2 + \|F|_{\mathcal{W}_1}\|_{\mathcal{W}_1'}^2\,.$$

□

3.2. **Variational boundary value problems with a linear quantity of interest.** Let $\mathcal{U}$ and $\mathcal{V}$ be Hilbert spaces over $\mathbb{R}$. In the abstract setting, the variational boundary value problems we consider are posed over such spaces using a continuous bilinear form $b : \mathcal{U} \times \mathcal{V} \to \mathbb{R}$. Here, $\mathcal{U}$ is called the *trial space* and $\mathcal{V}$ is called the *test space*. The members of both spaces $\mathcal{U}$ and $\mathcal{V}$ may have many components but are routinely called *functions*. For instance, *trial functions* $\boldsymbol{u} \in \mathcal{U}$ and *test functions* $\boldsymbol{v} \in \mathcal{V}$.

For a given functional $\ell \in \mathcal{V}'$, called the *load*, we define the (*primal*) *solution* to be the unique function $\boldsymbol{u}^\star \in \mathcal{U}$ satisfying

$$b(\boldsymbol{u}^\star, \boldsymbol{v}) = \ell(\boldsymbol{v})\,, \quad \forall\, \boldsymbol{v} \in \mathcal{V}\,. \tag{3.1}$$

Note that the bilinear form $b$ naturally generates a continuous linear operator $\mathcal{B} : \mathcal{U} \to \mathcal{V}'$ and, taking into account the reflexivity of $\mathcal{V}$, also generates its dual $\mathcal{B}' : \mathcal{V} \to \mathcal{U}'$:

$$\langle \mathcal{B}\boldsymbol{u}, \boldsymbol{v}\rangle = \langle \boldsymbol{u}, \mathcal{B}'\boldsymbol{v}\rangle = b(\boldsymbol{u}, \boldsymbol{v})\,, \quad \forall\, \boldsymbol{u} \in \mathcal{U}\,,\ \boldsymbol{v} \in \mathcal{V}\,. \tag{3.2}$$

Often, in practical simulations, we are not altogether interested in the global features of the solution $\boldsymbol{u}^\star = \mathcal{B}^{-1}\ell$ of a given problem of type (3.1). Instead, we are interested in a derived quantity $G(\boldsymbol{u}^\star)$, called the *quantity of interest* (QoI). In this context, we choose to call the corresponding functional $G \in \mathcal{U}'$ the *goal*.

Formally, observe that

$$G(\boldsymbol{u}^\star) = \langle G, \mathcal{B}^{-1}\ell\rangle = \langle \ell, (\mathcal{B}')^{-1}\rangle = \ell(\boldsymbol{v}^\star)\,, \tag{3.3}$$

where $\boldsymbol{v}^\star \in (\mathcal{B}')^{-1}G$. We will call any such function $\boldsymbol{v}^\star \in \mathcal{V}$—which may be unique—a (*test space*) *influence function* and note that it satisfies

$$b(\boldsymbol{u},\boldsymbol{v}^\star) = G(\boldsymbol{u})\,, \quad \forall\, \boldsymbol{u} \in \mathcal{U}\,. \tag{3.4}$$

Notice from the definition above that $\boldsymbol{v}^\star$ acts like a generalized Green's function defined for the functional $G \in \mathcal{U}'$. Moreover, notice that the QoI can be calculated from the solution of the primal problem (3.1) or from a solution of the dual problem (3.4).

3.3. **Well-posedness.** Clearly, in order for solutions to (3.1) and (3.4) to exist, we must assume that $\ell \in \mathrm{Range}(\mathcal{B})$ and $G \in \mathrm{Range}(\mathcal{B}')$. In order to simplify our analysis in this section, we will simply assume both that $\mathrm{Range}(\mathcal{B}') = \mathcal{U}'$ and $\mathrm{Range}(\mathcal{B}) = \mathcal{V}'$. These assumptions will make both $\boldsymbol{u}^\star$ and $\boldsymbol{v}^\star$ unique, but only the first assumption will ultimately be necessary for our analysis in later sections.

**Assumption 1.** The bilinear form $b : \mathcal{U}\times\mathcal{V} \to \mathbb{R}$ is bounded with continuity constant $M < +\infty$, where

$$\sup_{\boldsymbol{u}\in\mathcal{U}}\sup_{\boldsymbol{v}\in\mathcal{V}} \frac{b(\boldsymbol{u},\boldsymbol{v})}{\|\boldsymbol{u}\|_{\mathcal{U}}\|\boldsymbol{v}\|_{\mathcal{V}}} = M\,, \tag{3.5}$$

and it satisfies the Banach-Babuška-Nečas inf-sup condition with stability constant $\gamma > 0$, where

$$\inf_{\boldsymbol{u}\in\mathcal{U}}\sup_{\boldsymbol{v}\in\mathcal{V}} \frac{b(\boldsymbol{u},\boldsymbol{v})}{\|\boldsymbol{u}\|_{\mathcal{U}}\|\boldsymbol{v}\|_{\mathcal{V}}} = \gamma\,.$$

Notice that Assumption 1 readily implies

$$\gamma\|\boldsymbol{u}\|_{\mathcal{U}} \le \|\mathcal{B}\boldsymbol{u}\|_{\mathcal{V}'} \le M\|\boldsymbol{u}\|_{\mathcal{U}}\,, \quad \forall\, \boldsymbol{u}\in\mathcal{U}\,. \tag{3.6}$$

3.4. **Minimum residual principles.** Let $\mathcal{U}_h \subseteq \mathcal{U}$ be a finite-dimensional subspace of the trial space. It desirable for us to characterize the *optimal* (i.e. minimal error) solution $\boldsymbol{u}_h^{\mathrm{opt}} \in \mathcal{U}_h$ to the primal problem (3.1). Because the exact solution $\boldsymbol{u}^\star = \mathcal{B}^{-1}\ell$ is inaccessible *a priori*, the optimal solution $\boldsymbol{u}_h^{\mathrm{opt}}$ cannot be defined by explicitly invoking it. Instead, we propose the following minimum residual principle:

$$\boldsymbol{u}_h^{\mathrm{opt}} = \operatorname*{arg\,min}_{\boldsymbol{u}_h\in\mathcal{U}_h} \|\mathcal{B}\boldsymbol{u}_h - \ell\|_{\mathcal{V}'}^2\,. \tag{3.7}$$

Observe that $\|\mathcal{B}\boldsymbol{u}_h - \ell\|_{\mathcal{V}'}^2 = \langle \mathcal{B}\boldsymbol{u}_h - \ell, \mathcal{R}_{\mathcal{V}}^{-1}(\mathcal{B}\boldsymbol{u}_h - \ell)\rangle$, for all $\boldsymbol{u}_h \in \mathcal{U}_h$. Therefore, the first-order optimality condition associated with (3.7) is equivalent to the variational equation

$$\langle \mathcal{B}\boldsymbol{u}_h^{\mathrm{opt}}, \mathcal{R}_{\mathcal{V}}^{-1}\,\mathcal{B}\boldsymbol{u}_h\rangle = \langle \ell, \mathcal{R}_{\mathcal{V}}^{-1}\,\mathcal{B}\boldsymbol{u}_h\rangle\,, \quad \forall\, \boldsymbol{u}_h \in \mathcal{U}_h\,. \tag{3.8}$$

We will now introduce a number of definitions which will become very important soon. Define the *trial-to-test* operator $\Theta = \mathcal{R}_{\mathcal{V}}^{-1}\,\mathcal{B}$ and the *induced (trial space) Riesz operator* $\mathcal{A} = \mathcal{B}'\,\mathcal{R}_{\mathcal{V}}^{-1}\,\mathcal{B} = \mathcal{B}'\Theta = \Theta'\mathcal{B}$. Define the *energy* inner product $a : \mathcal{U}\times\mathcal{U} \to \mathbb{R}$,

$$a(\boldsymbol{u},\widetilde{\boldsymbol{u}}) = \langle \mathcal{A}\boldsymbol{u}, \widetilde{\boldsymbol{u}}\rangle = \langle \boldsymbol{u}, \mathcal{A}\widetilde{\boldsymbol{u}}\rangle\,, \quad \forall\, \boldsymbol{u},\widetilde{\boldsymbol{u}} \in \mathcal{U}\,,$$

and also the *induced load* $F \in \mathcal{U}'$, $F = \mathcal{B}'\,\mathcal{R}_{\mathcal{V}}^{-1}\,\ell = \Theta'\ell$. With these definitions, we can reexpress (3.8) as either

$$b(\boldsymbol{u}_h^{\mathrm{opt}}, \boldsymbol{v}_h) = \ell(\boldsymbol{v}_h)\,, \quad \forall\, \boldsymbol{v}_h \in \Theta(\mathcal{U}_h)\,, \tag{3.9}$$

which is known as the *optimal test function expression* [26, 27], or as

$$a(\boldsymbol{u}_h^{\mathrm{opt}}, \boldsymbol{u}_h) = F(\boldsymbol{u}_h)\,, \quad \forall\, \boldsymbol{u}_h \in \mathcal{U}_h\,, \tag{3.10}$$

which we will refer to as the *coercive expression*. In (3.9), the subspace $\Theta(\mathcal{U}_h) \subseteq \mathcal{V}$ is referred to as the *optimal test space*, and its members are called *optimal test functions*. Recall (3.6) and notice that the energy inner product $a : \mathcal{U}\times\mathcal{U} \to \mathbb{R}$, induces a norm

$$|||\boldsymbol{u}|||_{\mathcal{U}}^2 = a(\boldsymbol{u},\boldsymbol{u}) = \langle \mathcal{A}\boldsymbol{u},\boldsymbol{u}\rangle = \|\mathcal{B}\boldsymbol{u}\|_{\mathcal{V}'}^2\,, \quad \forall\, \boldsymbol{u}\in\mathcal{U}\,.$$

We will refer to $|||\cdot|||_{\mathcal{U}}$ as the *energy norm*.

3.5. **Dual minimum residual principles.** Let us return to (3.1) and (3.4). These problems can be rewritten $\mathcal{B}\boldsymbol{u}^\star = \ell$ and $\mathcal{B}'\boldsymbol{v}^\star = G$, respectively. Applying the dual of trial-to-test operator, $\Theta'$, to both sides of the former equation delivers $\mathcal{A}\boldsymbol{u}^\star = F$. When $\mathcal{U}_h = \mathcal{U}$, comparison with (3.10) then demonstrates that $\boldsymbol{u}_h^{\mathrm{opt}} = \boldsymbol{u}^\star$. Meanwhile, if we define $\boldsymbol{v}^\star = \Theta\boldsymbol{\omega}^\star$, for some $\boldsymbol{\omega}^\star \in \mathcal{U}$, the latter equation delivers $\mathcal{A}\boldsymbol{\omega}^\star = G$. With little ambiguity, we shall also refer to $\boldsymbol{\omega}^\star \in \mathcal{U}$ as the (trial space) influence function. Its distinction from the (test space) influence function $\boldsymbol{v}^\star \in \mathcal{V}$ should always be clear from context.

A Bubnov-Galerkin approximation, $\boldsymbol{\omega}_h$, of the influence function $\boldsymbol{\omega}^\star$ can be readily characterized:

$$a(\boldsymbol{\omega}_h, \boldsymbol{u}_h) = G(\boldsymbol{u}_h)\,, \quad \forall\, \boldsymbol{u}_h \in \mathcal{U}_h\,. \tag{3.11}$$

The discrete solution above, $\boldsymbol{\omega}_h = \boldsymbol{\omega}_h^{\mathrm{en}}$, can easily be seen to come from the Ritz method following from the quadratic energy principle

$$\boldsymbol{\omega}_h^{\mathrm{en}} = \operatorname*{arg\,min}_{\boldsymbol{u}_h \in \mathcal{U}_h} \big( a(\boldsymbol{u}_h, \boldsymbol{u}_h) - 2G(\boldsymbol{u}_h) \big)\,.$$

Equivalently [52, Section 2.3], the approximation can be characterized as $\boldsymbol{\omega}_h = \boldsymbol{\omega}_h^{\mathrm{opt}}$, where $\boldsymbol{\omega}_h^{\mathrm{opt}}$ is defined as the optimal solution coming from the minimum residual problem

$$\boldsymbol{\omega}_h^{\mathrm{opt}} = \operatorname*{arg\,min}_{\boldsymbol{u}_h \in \mathcal{U}_h} |||\mathcal{A}\boldsymbol{u}_h - G|||_{\mathcal{U}'}\,.$$

Because this second perspective will be the most insightful, we will denote the approximation $\boldsymbol{\omega}_h$ as $\boldsymbol{\omega}_h^{\mathrm{opt}}$ from now on.

Observe that an approximation to the influence function $\boldsymbol{v}^\star \in \mathcal{V}$, denoted $\boldsymbol{v}_h^{\mathrm{opt}}$, where

$$b(\boldsymbol{u}_h, \boldsymbol{v}_h^{\mathrm{opt}}) = G(\boldsymbol{u}_h)\,, \quad \forall\, \boldsymbol{u}_h \in \mathcal{U}_h\,,$$

can be recovered by post-processing the computed $\boldsymbol{\omega}_h^{\mathrm{opt}}$. Indeed, $\boldsymbol{v}_h^{\mathrm{opt}} = \Theta\boldsymbol{\omega}_h^{\mathrm{opt}}$. Note that this makes the discrete influence function $\boldsymbol{v}_h^{\mathrm{opt}} \in \Theta(\mathcal{U}_h)$ an optimal test function and

$$\boldsymbol{v}_h^{\mathrm{opt}} = \operatorname*{arg\,min}_{\boldsymbol{v}_h \in \Theta(\mathcal{U}_h)} |||\mathcal{B}'\boldsymbol{v}_h - G|||_{\mathcal{U}'}\,.$$

3.6. **The mixed method interpretation.** Reflect upon (3.10) and (3.11):

$$a(\boldsymbol{u}_h^{\mathrm{opt}}, \boldsymbol{u}_h) = F(\boldsymbol{u}_h) \qquad \text{and} \qquad a(\boldsymbol{\omega}_h^{\mathrm{opt}}, \boldsymbol{u}_h) = G(\boldsymbol{u}_h)\,, \quad \forall\, \boldsymbol{u}_h \in \mathcal{U}_h\,.$$

Both optimal discrete solutions $\boldsymbol{u}_h^{\mathrm{opt}}$ and $\boldsymbol{\omega}_h^{\mathrm{opt}}$ are unique. Indeed, define the Riesz representation of the residual $\boldsymbol{\psi}_h^{\mathrm{opt}} = \mathcal{R}_{\mathcal{V}}^{-1}(\mathcal{B}\boldsymbol{u}_h^{\mathrm{opt}} - \ell)$ and recall that $\boldsymbol{v}_h^{\mathrm{opt}} = \mathcal{R}_{\mathcal{V}}^{-1}\,\mathcal{B}\boldsymbol{\omega}_h^{\mathrm{opt}}$. Now, observe that

$$b(\boldsymbol{u}_h, \boldsymbol{\psi}_h^{\mathrm{opt}}) = \langle \mathcal{B}'\,\mathcal{R}_{\mathcal{V}}^{-1}(\mathcal{B}\boldsymbol{u}_h^{\mathrm{opt}} - \ell), \boldsymbol{u}_h \rangle = a(\boldsymbol{u}_h^{\mathrm{opt}}, \boldsymbol{u}_h) - F(\boldsymbol{u}_h) = 0\,, \quad \forall\, \boldsymbol{u}_h \in \mathcal{U}_h\,,$$

and

$$b(\boldsymbol{u}_h, \boldsymbol{v}_h^{\mathrm{opt}}) = \langle \mathcal{B}'\,\mathcal{R}_{\mathcal{V}}^{-1}\,\mathcal{B}\boldsymbol{\omega}_h^{\mathrm{opt}}, \boldsymbol{u}_h \rangle = a(\boldsymbol{\omega}_h^{\mathrm{opt}}, \boldsymbol{u}_h) = G(\boldsymbol{u}_h)\,, \quad \forall\, \boldsymbol{u}_h \in \mathcal{U}_h\,.$$

Therefore, note the following saddle-point characterizations of the idealized primal and dual problems above:

$$\begin{cases} (\boldsymbol{\psi}_h^{\mathrm{opt}}, \boldsymbol{v})_{\mathcal{V}} - b(\boldsymbol{u}_h^{\mathrm{opt}}, \boldsymbol{v}) = -\ell(\boldsymbol{v})\,, & \forall\, \boldsymbol{v} \in \mathcal{V}\,, \\ b(\boldsymbol{u}_h, \boldsymbol{\psi}_h^{\mathrm{opt}}) \qquad\qquad\quad = 0\,, & \forall\, \boldsymbol{u}_h \in \mathcal{U}_h\,, \end{cases} \tag{3.12}$$

and

$$\begin{cases} (\boldsymbol{v}_h^{\mathrm{opt}}, \boldsymbol{v})_{\mathcal{V}} - b(\boldsymbol{\omega}_h^{\mathrm{opt}}, \boldsymbol{v}) = 0\,, & \forall\, \boldsymbol{v} \in \mathcal{V}\,, \\ b(\boldsymbol{u}_h, \boldsymbol{v}_h^{\mathrm{opt}}) \qquad\qquad\quad = G(\boldsymbol{u}_h)\,, & \forall\, \boldsymbol{u}_h \in \mathcal{U}_h\,, \end{cases} \tag{3.13}$$

respectively. We will refer to these as the *idealized* mixed methods.

Similarly, by introducing the functions $\boldsymbol{\psi}^\star = \mathcal{R}_{\mathcal{V}}^{-1}(\mathcal{B}\boldsymbol{u}^\star - \ell)$ and $\boldsymbol{v}^\star = \mathcal{R}_{\mathcal{V}}^{-1}\,\mathcal{B}\boldsymbol{\omega}^\star$, note that

$$\begin{cases} (\boldsymbol{\psi}^\star, \boldsymbol{v})_{\mathcal{V}} - b(\boldsymbol{u}^\star, \boldsymbol{v}) = -\ell(\boldsymbol{v})\,, & \forall\, \boldsymbol{v} \in \mathcal{V}\,, \\ b(\boldsymbol{u}, \boldsymbol{\psi}^\star) \qquad\qquad = 0\,, & \forall\, \boldsymbol{u} \in \mathcal{U}\,, \end{cases} \tag{3.14}$$

and

$$\begin{cases} (\boldsymbol{v}^\star, \boldsymbol{v})_{\mathcal{V}} - b(\boldsymbol{\omega}^\star, \boldsymbol{v}) = 0\,, & \forall\, \boldsymbol{v} \in \mathcal{V}\,, \\ b(\boldsymbol{u}, \boldsymbol{v}^\star) = G(\boldsymbol{u})\,, & \forall\, \boldsymbol{u} \in \mathcal{U}\,. \end{cases} \tag{3.15}$$

Because the null space of $\mathcal{B}'$ may be non-trivial, an important distinction between (3.4) and (3.15) is that (3.15) always has a unique solution. Indeed, for arbitrary $\ell \in \mathcal{V}'$ and $G \in \mathcal{U}'$, existence and uniqueness can be proven for each saddle point problem (3.12)–(3.15) by simply checking Brezzi's conditions [10]. In fact, for any arbitrary $\ell \in \mathrm{Range}(\mathcal{B})$, $\boldsymbol{\psi}^\star = 0$. Meanwhile, for any arbitrary $G \in \mathcal{U}'$, $\boldsymbol{v}^\star \perp \mathrm{Null}(\mathcal{B}')$ is the unique influence function satisfying (3.4).

**Theorem 3.1.** *Let Assumption 1 hold. Let $\ell \in \mathrm{Range}(\mathcal{B})$ and $G \in \mathcal{U}'$. Then $(\boldsymbol{u}^\star, \boldsymbol{\psi}^\star)$ and $(\boldsymbol{\omega}^\star, \boldsymbol{v}^\star)$ exist and are unique. Furthermore, $\boldsymbol{\psi}^\star = 0$ and $\boldsymbol{v}^\star \perp \mathrm{Null}(\mathcal{B}')$. Similarly, $(\boldsymbol{u}_h^{\mathrm{opt}}, \boldsymbol{\psi}_h^{\mathrm{opt}})$ and $(\boldsymbol{\omega}_h^{\mathrm{opt}}, \boldsymbol{v}_h^{\mathrm{opt}})$ exist and are unique and $\boldsymbol{\psi}_h^{\mathrm{opt}}, \boldsymbol{v}_h^{\mathrm{opt}} \perp \mathrm{Null}(\mathcal{B}')$.*

*Proof.* The operator form of (3.14) is

$$\begin{cases} \mathcal{R}_{\mathcal{V}}\, \boldsymbol{\psi}^\star - \mathcal{B}\boldsymbol{u}^\star = -\ell\,, \\ \mathcal{B}'\boldsymbol{\psi}^\star = 0\,. \end{cases}$$

Applying $\Theta' = \mathcal{B}'\,\mathcal{R}_{\mathcal{V}}^{-1}$ to the first equation and substituting the second delivers $\mathcal{B}'\,\mathcal{R}_{\mathcal{V}}^{-1}\,\mathcal{B}\boldsymbol{u}^\star = \mathcal{B}'\,\mathcal{R}_{\mathcal{V}}^{-1}\,\ell$ or, equivalently, $\mathcal{A}\boldsymbol{u}^\star = F$. By the Riesz Representation Theorem, this equation has a unique solution. Therefore, $\boldsymbol{u}^\star$ exists and is unique. Note that $\ell \in \mathrm{Range}(\mathcal{B}) = \mathrm{Null}(\mathcal{B}')^\perp$ and $\mathcal{B}\boldsymbol{u}^\star = \ell$. Therefore, simultaneously, $\boldsymbol{u}^\star = \mathcal{B}^{-1}\ell$. Moreover, $\boldsymbol{\psi}_1^\star = \mathcal{R}_{\mathcal{V}}^{-1}(\mathcal{B}\boldsymbol{u}^\star - \ell) \perp \mathrm{Null}(\mathcal{B}')$ is at least one candidate for the solution component, $\boldsymbol{\psi}^\star$. Let $\boldsymbol{\psi}^\star = \boldsymbol{\psi}_0^\star + \boldsymbol{\psi}_1^\star$, where $\boldsymbol{\psi}_0^\star \in \mathrm{Null}(\mathcal{B}')$ is arbitrary. Testing the first equation of (3.14) with $\boldsymbol{\psi}_0^\star$, we find $\|\boldsymbol{\psi}_0^\star\|_{\mathcal{V}}^2 = 0$. Therefore, $\boldsymbol{\psi}^\star = \boldsymbol{\psi}_1^\star = \mathcal{R}_{\mathcal{V}}^{-1}(\mathcal{B}\boldsymbol{u}^\star - \ell) = \mathcal{R}_{\mathcal{V}}^{-1}(\mathcal{B}\mathcal{B}^{-1}\ell - \ell) = 0$.

Observe that $\mathrm{Null}((\mathcal{B}')')^\perp = \mathrm{Null}(\mathcal{B})^\perp = \mathcal{U}'$. As was demonstrated above for $\boldsymbol{u}^\star$, the influence function $\boldsymbol{\omega}^\star$ can be shown to be the unique solution to $\mathcal{A}\boldsymbol{\omega}^\star = G$. Consider the decomposition $\boldsymbol{v}^\star = \boldsymbol{v}_0^\star + \boldsymbol{v}_1^\star$, where $\boldsymbol{v}_0^\star \in \mathrm{Null}(\mathcal{B}')$ is arbitrary and $\boldsymbol{v}_1^\star = \Theta\mathcal{A}^{-1}G$. Testing the first equation of (3.15) with $\boldsymbol{v}_0^\star$, it also follows that $\|\boldsymbol{v}_0^\star\|_{\mathcal{V}} = 0$ and so $\boldsymbol{v}^\star \perp \mathrm{Null}(\mathcal{B}')$.

The second pair of statements are proved similarly. □

3.7. **The optimal test norm.** Suppose that a test norm $|||\cdot|||_{\mathcal{V}}$ exists such that $|||\mathcal{B}\cdot|||_{\mathcal{V}'} = \|\cdot\|_{\mathcal{U}}$. With this norm, which we will call the *optimal test norm*,

$$\langle \mathcal{R}_{\mathcal{U}}\,\boldsymbol{u}, \widetilde{\boldsymbol{u}}\rangle = (\boldsymbol{u}, \widetilde{\boldsymbol{u}})_{\mathcal{U}} \overset{\mathrm{opt}}{=} (\mathcal{B}\boldsymbol{u}, \mathcal{B}\widetilde{\boldsymbol{u}})_{\mathcal{V}'} = \langle \mathcal{B}'\,\mathcal{R}_{\mathcal{V}}^{-1}\,\mathcal{B}\boldsymbol{u}, \widetilde{\boldsymbol{u}}\rangle\,, \quad \forall\, \boldsymbol{u}, \widetilde{\boldsymbol{u}} \in \mathcal{U}\,. \tag{3.16}$$

Here, the equality with the opt-notation indicates to the reader it holds only in the special $|||\cdot|||_{\mathcal{V}}$ setting.

Clearly, $\mathcal{R}_{\mathcal{U}} \overset{\mathrm{opt}}{=} \mathcal{B}'\,\mathcal{R}_{\mathcal{V}}^{-1}\,\mathcal{B} = \mathcal{A}$. In the case that $\mathcal{B}$ is an isomorphism, the Riesz operator for the optimal test norm is $\mathcal{R}_{\mathcal{V}} \overset{\mathrm{opt}}{=} \mathcal{B}\,\mathcal{R}_{\mathcal{U}}^{-1}\,\mathcal{B}'$. Here, the definitions of these special Riesz operators are easy to remember by observing that they are the unique isometries where the following diagram commutes:

$$\begin{array}{ccc} \mathcal{U}' & \xleftarrow{\ \mathcal{B}'\ } & \mathcal{V} \\ {\scriptstyle \mathcal{R}_{\mathcal{U}}}\uparrow & & \downarrow{\scriptstyle \mathcal{R}_{\mathcal{V}}} \\ \mathcal{U} & \xrightarrow{\ \mathcal{B}\ } & \mathcal{V}' \end{array}\,.$$

Moreover, if $\|\cdot\|_{\mathcal{V}} \overset{\mathrm{opt}}{=} |||\cdot|||_{\mathcal{V}}$, then $M \overset{\mathrm{opt}}{=} \gamma \overset{\mathrm{opt}}{=} 1$ in (3.6) and $b(\boldsymbol{u}, \boldsymbol{v})$ becomes a duality pairing [11].

Alternatively, if $\mathrm{Null}(\mathcal{B}')^\perp \neq \{0\}$, then $|||\cdot|||_{\mathcal{V}}$ may be treated as a quotient norm or endowed with a Riesz map defined $\mathcal{R}_{\mathcal{V}} \overset{\mathrm{opt}}{=} \mathcal{B}\,\mathcal{R}_{\mathcal{U}}^{-1}\,\mathcal{B}' + \mathcal{P}$, where $\mathcal{P} : \mathcal{V} \to \mathrm{Null}(\mathcal{B}')$ is an orthogonal projection. The general setting we wish to consider is summarized by the following lemma.

**Lemma 3.2.** *Let Assumption 1 hold. Let $\mathcal{P} : \mathcal{V} \to \mathrm{Null}(\mathcal{B}')$ be the orthogonal projection. Define*

$$|||\boldsymbol{u}|||_{\mathcal{U}} = \|\mathcal{B}\boldsymbol{u}\|_{\mathcal{V}'} \qquad \textit{and} \qquad |||\boldsymbol{v}|||_{\mathcal{V}}^2 = \|\mathcal{B}'\boldsymbol{v}\|_{\mathcal{U}'}^2 + \|\mathcal{P}\boldsymbol{v}\|_{\mathcal{V}}^2 \tag{3.17}$$

*Then*

$$\|\boldsymbol{u}\|_{\mathcal{U}} = |||\mathcal{B}\boldsymbol{u}|||_{\mathcal{V}'} \qquad \textit{and} \qquad \|\boldsymbol{v}\|_{\mathcal{V}}^2 = |||\mathcal{B}'\boldsymbol{v}|||_{\mathcal{U}'}^2 + \|\mathcal{P}\boldsymbol{v}\|_{\mathcal{V}}^2\,, \tag{3.18}$$

*for all* $\boldsymbol{u} \in \mathcal{U}$ *and* $\boldsymbol{v} \in \mathcal{V}$.

*Proof.* The first relationship is immediate. To prove the second, let $\mathcal{P}^\perp : \mathcal{V} \to \mathrm{Null}(\mathcal{B}')^\perp$ be the orthogonal projection complementary to $\mathcal{P}$. By the Closed Range Theorem, observe that $\mathrm{Null}(\mathcal{B}')^\perp = \mathrm{Range}(\mathcal{B})$ and both $\mathrm{Null}(\mathcal{B}')$ and $\mathrm{Null}(\mathcal{B}')^\perp$ are closed. Therefore, $\|\boldsymbol{v}\|_{\mathcal{V}}^2 = \|\mathcal{P}\boldsymbol{v}\|_{\mathcal{V}}^2 + \|\mathcal{P}^\perp\boldsymbol{v}\|_{\mathcal{V}}^2$, where

$$\|\mathcal{P}^\perp\boldsymbol{v}\|_{\mathcal{V}} = \sup_{\boldsymbol{v}_1 \in \mathrm{Null}(\mathcal{B}')^\perp} \frac{(\boldsymbol{v},\boldsymbol{v}_1)_{\mathcal{V}}}{\|\boldsymbol{v}_1\|_{\mathcal{V}}} = \sup_{\boldsymbol{u}\in\mathcal{U}} \frac{(\boldsymbol{v},\mathcal{R}_{\mathcal{V}}^{-1}\,\mathcal{B}\boldsymbol{u})_{\mathcal{V}}}{\|\,\mathcal{R}_{\mathcal{V}}^{-1}\,\mathcal{B}\boldsymbol{u}\|_{\mathcal{V}}} = \sup_{\boldsymbol{u}\in\mathcal{U}} \frac{\langle \boldsymbol{v},\mathcal{B}\boldsymbol{u}\rangle}{\|\mathcal{B}\boldsymbol{u}\|_{\mathcal{V}'}} = |||\mathcal{B}'\boldsymbol{v}|||_{\mathcal{U}'}\,.$$

□

**Corollary 3.3.** *Let Assumption 1 hold. Let* $\ell \in \mathrm{Range}(\mathcal{B})$ *and* $G \in \mathcal{U}'$. *Then, with the notation of Theorem 3.1 and Lemma 3.2,*

$$\|\boldsymbol{u}^\star - \boldsymbol{u}\|_{\mathcal{U}} = |||\mathcal{B}\boldsymbol{u} - \ell|||_{\mathcal{V}'} \qquad \textit{and} \qquad \|\boldsymbol{v}^\star - \boldsymbol{v}\|_{\mathcal{V}}^2 = \|\mathcal{P}\boldsymbol{v}\|_{\mathcal{V}}^2 + |||\mathcal{B}'\boldsymbol{v} - G|||_{\mathcal{U}'}^2\,, \tag{3.19}$$

*for all* $\boldsymbol{u} \in \mathcal{U}$ *and* $\boldsymbol{v} \in \mathcal{V}$.

*Proof.* This immediately follows from (3.18) and Theorem 3.1, namely $\boldsymbol{v}^\star \perp \mathrm{Null}(\mathcal{B}')$. □

Unfortunately, actually working directly with the optimal test norm is not computationally feasible. This is, in part, due to presence of the inverse of the trial space Riesz operator in the term $\mathcal{B}\,\mathcal{R}_{\mathcal{U}}^{-1}\,\mathcal{B}'$. Nevertheless, for special bilinear forms, we can often still come suitably close as we will motivate in the next subsection.

3.8. **Ultraweak variational formulations.** Many modeling problems are originally derived in a strong form:

$$\mathcal{L}\boldsymbol{u} = f\,, \tag{3.20}$$

where $\mathcal{L}$ is a linear operator and $f$ is a prescribed function. Multiplying this equation with test functions and formally integrating by parts over the domain $\Omega$ delivers the so-called *ultraweak variation formulation*:

$$(\boldsymbol{u},\mathcal{L}^*\boldsymbol{v})_\Omega = (f,\boldsymbol{v})_\Omega\,, \quad \forall\, \boldsymbol{v} \in \mathcal{V}\,.$$

This formulation of the abstract modeling problem can be proved to be well-posed and have a solution consistent with (3.20), for all compatible functions $f \in L^2(\Omega)$, assuming that $\mathcal{L}^* : \mathcal{V} \to L^2(\Omega)$ is a closed linear operator where $\mathcal{V} = \mathrm{Dom}(\mathcal{L}^*) \subseteq L^2(\Omega)$ is dense [23].

The specific variational formulations we will be concerned with involve hybridized variables. With the definition $\mathcal{V}_K = \{\boldsymbol{v}|_K : \boldsymbol{v} \in \mathcal{V}\}$, they reduce to problems which can be posed in the following functional setting:

$$\begin{cases} \text{Find } \boldsymbol{u}^\star = (u^\star, \hat{u}^\star) \in \mathcal{U} = L^2(\Omega) \times \hat{\mathcal{U}} : \\ (u^\star, \mathcal{L}^*\boldsymbol{v})_\Omega + \langle \hat{u}^\star, \boldsymbol{v}\rangle_{\partial\mathcal{T}} = (f,\boldsymbol{v})_\Omega\,, \quad \forall\, \boldsymbol{v} \in \mathcal{V} = \displaystyle\prod_{K\in\mathcal{T}} \mathcal{V}_K\,. \end{cases} \tag{3.21}$$

Here, $\langle\cdot,\cdot\rangle_{\partial\mathcal{T}}$ denotes the mesh-dependent functional pairing:

$$\langle \hat{u}, \boldsymbol{v}\rangle_{\partial\mathcal{T}} = \sum_{K\in\mathcal{T}} \langle \hat{u}, \mathrm{tr}_{\mathcal{V}}^K\, \boldsymbol{v}\rangle\,, \quad \forall\, \hat{u} \in \hat{\mathcal{U}}, \boldsymbol{v} \in \mathcal{V}\,,$$

where $\mathrm{tr}_{\mathcal{V}}^K$ is a continuous trace operator. Usually, such variational problems are derived by multiplying (3.20) with so-called *broken test functions* $\boldsymbol{v}$, formally integrating by parts *element-wise*, and then deriving the appropriate interface space $\hat{\mathcal{U}}$ for well-posedness. Generally, this space $\hat{\mathcal{U}}$ may be thought of a space of Lagrange multipliers, involving only functions associated to the mesh skeleton, which complements the *broken* or discontinuous nature of the test space, $\mathcal{V} = \prod_{K\in\mathcal{T}} \mathcal{V}_K$ [14, 28].

From (3.21), define the broken ultraweak bilinear form $b(\boldsymbol{u},\boldsymbol{v}) = (u,\mathcal{L}^*\boldsymbol{v})_\Omega + \langle \hat{u},\boldsymbol{v}\rangle_{\partial\mathcal{T}}$ and the load functional $\ell(\boldsymbol{v}) = (f,\boldsymbol{v})_\Omega$. Recall definition (3.2). For all $\boldsymbol{v}_1 \perp \mathrm{Null}(\mathcal{B}')$, observe that the corresponding optimal test norm (3.17) can be expressed as

$$\tag{3.22} |||\boldsymbol{v}_1|||_{\mathcal{V}}^2 = \|\mathcal{B}'\boldsymbol{v}_1\|_{\mathcal{U}}^2 = \left( \sup_{\boldsymbol{u}\in\mathcal{U}} \frac{(u,\mathcal{L}^*\boldsymbol{v}_1)_\Omega + \langle \hat{u},\boldsymbol{v}_1\rangle_{\partial\mathcal{T}}}{\left(\|u\|_\Omega^2 + \|\hat{u}\|_{\hat{\mathcal{U}}}^2\right)^{1/2}} \right)^2 = \|\mathcal{L}^*\boldsymbol{v}_1\|_\Omega^2 + \sup_{\hat{u}\in\hat{\mathcal{U}}} \frac{\langle \hat{u},\boldsymbol{v}_1\rangle_{\partial\mathcal{T}}^2}{\|\hat{u}\|_{\hat{\mathcal{U}}}^2}\,,$$

by Lemma 3.1. Here, the first term is readily computable but the second term is generally not. As a compromise, we may instead endow the full test space $\mathcal{V}$ with the quasi-optimal test norm $\|\boldsymbol{v}\|_{\mathcal{L}^*}^2 = \|\mathcal{L}^*\boldsymbol{v}\|_\Omega^2 + \|\mathcal{P}_{\mathcal{L}^*}\boldsymbol{v}\|_\Omega$, where $\mathcal{P}_{\mathcal{L}^*} : \mathcal{V} \to \mathrm{Null}(\mathcal{L}^*) \supseteq \mathrm{Null}(\mathcal{B}')$ is the orthogonal projection. Alternatively, a related norm called the (*adjoint*) *graph norm* is often used in practice:

$$\tag{3.23} \|\boldsymbol{v}\|_{\mathcal{L}^*,\alpha}^2 = \|\mathcal{L}^*\boldsymbol{v}\|_\Omega^2 + \alpha^2\|\boldsymbol{v}\|_\Omega^2\,.$$

Computations with this graph norm, for sufficiently small values of the parameter $\alpha$, give quasi-optimal results for a wide variety of problems. For a deeper account of test norm choices, as well as improvements on the graph norm for some singular perturbation problems, see [16,47] and references therein. Note that the graph norm easily decomposes into the following sum: $\|\boldsymbol{v}\|_{\mathcal{L}^*,\alpha}^2 = \sum_{K\in\mathcal{T}} \left(\|\mathcal{L}^*\boldsymbol{v}\|_K^2 + \alpha^2\|\boldsymbol{v}\|_K^2\right)$, and each respective term can be viewed as an independent norm on $\mathcal{V}_K$. Because other norms and test spaces could be used in practice, we formulate the key properties of the graph norm above as an independent assumption.

**Assumption 2.** The test space $\mathcal{V}$ is *broken*; that is, $\mathcal{V} = \prod_{K\in\mathcal{T}} \mathcal{V}_K$, where $\mathcal{V}_K = \{\boldsymbol{v}|_K : \boldsymbol{v} \in \mathcal{V}\}$, for all $K \in \mathcal{T}$. Moreover, the corresponding test space norm $\|\cdot\|_{\mathcal{V}}$ is *localizable*; that is, $\|\cdot|_K\|_{\mathcal{V}}$ is also a norm, for all $K \in \mathcal{T}$.

## 4. DPG* METHODS AND DUALITY OF THE ERROR IN THE QOI

All of the expressions in the previous section were derived under the assumption that the test space Riesz map $\mathcal{R}_{\mathcal{V}}$ can be inverted over the entire space $\mathcal{V}$. Although this assumption is a helpful guide for designing DPG methods, it cannot be made during careful numerical analysis. In practice, the inversion of the Riesz map is carried out only approximately, on a finite-dimensional subspace $\mathcal{V}_r \subsetneq \mathcal{V}$.

### 4.1. **DPG and DPG* methods.**

Let Assumption 2 hold. A *practical* discontinuous Petrov–Galerkin (DPG) method with optimal test functions is a finite element method defined by the problem

$$\tag{4.1} \begin{cases} \text{Find } \boldsymbol{u}_{h,r} \in \mathcal{U}_h : \\ b(\boldsymbol{u}_{h,r},\Theta_r\boldsymbol{w}_h) = \ell(\Theta_r\boldsymbol{w}_h)\,, \qquad \forall\,\boldsymbol{w}_h \in \mathcal{U}_h\,, \end{cases}$$

where $\Theta_r : \mathcal{U}_h \to \mathcal{V}_r$ is an approximate trial-to-test operator defined by the chosen inner product $(\cdot,\cdot)_{\mathcal{V}}$ on the test space $\mathcal{V}$:

$$\tag{4.2} (\Theta_r\boldsymbol{w}_h,\boldsymbol{v}_r)_{\mathcal{V}} = b(\boldsymbol{w}_h,\boldsymbol{v}_r)\,, \quad \forall\,\boldsymbol{w}_h \in \mathcal{U}_h\,,\ \boldsymbol{v}_r \in \mathcal{V}_r\,.$$

Alternatively, a practical DPG* method is defined by the fully discrete dual problem

$$\tag{4.3} \begin{cases} \text{Find } \boldsymbol{v}_{h,r} = \Theta_r\boldsymbol{\omega}_{h,r} \in \mathcal{V}_r : \\ b(\boldsymbol{u}_h,\Theta_r\boldsymbol{\omega}_{h,r}) = G(\boldsymbol{u}_h)\,, \qquad \forall\,\boldsymbol{u}_h \in \mathcal{U}_h\,. \end{cases}$$

Define the *approximate* Riesz representation of the residual $\boldsymbol{\psi}_{h,r} = \Theta_r\boldsymbol{u}_{h,r}$. Generally, the range of a trial-to-test operator $\Theta_r$ is a proper closed subspace of the test space, $\mathrm{Range}(\Theta_r) \subsetneq \mathcal{V}_r$. We call any function in the range of $\Theta_r$ an *approximate optimal test function*. For general test norms, computation of the trial-to-test operator $\Theta_r$ in (4.2) requires the inversion of a large symmetric Gram matrix coming from the inner product $(\cdot,\cdot)_{\mathcal{V}}$. This is made feasible by using a broken test space and localizable test norms, given by Assumption 2. In this case, the Gram matrix becomes block-diagonal and the inverse can be computed locally and in parallel [52].

Under the following assumption of the existence of a Fortin operator, $\Pi_r : \mathcal{V} \to \mathcal{V}_r$, both (4.1) and (4.3) are well-posed [46].

**Assumption 3.** For the closed subspaces $\mathcal{U}_h \subseteq \mathcal{U}$ and $\mathcal{V}_r \subseteq \mathcal{V}$, there exists a bounded linear operator $\Pi_r : \mathcal{V} \to \mathcal{V}_r$ such that

$$b(\boldsymbol{u}_h, \boldsymbol{v} - \Pi_r \boldsymbol{v}) = 0 \,, \quad \forall\, \boldsymbol{u}_h \in \mathcal{U}_h \,,\, \boldsymbol{v} \in \mathcal{V} \,. \tag{4.4}$$

Although only Assumption 3 is actually required for well-posedness, when $\Pi_r$ is constructed analytically, it is often additionally designed to be a projection as well [57]. This additional structure can be exploited, so we also include the following stronger assumption.

**Assumption 4.** In addition to Assumption 3, $\Pi_r$ is projection, $\Pi_r \circ \Pi_r = \Pi_r$.

**Theorem 4.1.** *Let Assumption 1 and 3 hold. Let $\ell \in \mathrm{Range}(\mathcal{B})$ and $G \in \mathcal{U}'$. Then $(\boldsymbol{u}_{h,r}, \boldsymbol{\psi}_{h,r})$ and $(\boldsymbol{\omega}_{h,r}, \boldsymbol{v}_{h,r})$ exist and are unique. Moreover, $\boldsymbol{\psi}_{h,r}, \boldsymbol{v}_{h,r} \perp \mathrm{Null}(\mathcal{B}')$.*

*Proof.* We only produce the argument for $\boldsymbol{v}_{h,r}$. By definition, $\boldsymbol{v}_{h,r} = \Theta_r \boldsymbol{\omega}_{h,r} \in \mathcal{V}_r$. By Assumption 3, $\boldsymbol{\omega}_{h,r}$ exists and is unique. Therefore, $\boldsymbol{v}_{h,r}$ also exists and is also unique. Moreover, by the definition of $\Theta_r$,

$$(\boldsymbol{v}_{h,r}, \boldsymbol{v}_r)_{\mathcal{V}} = b(\boldsymbol{\omega}_{h,r}, \boldsymbol{v}_r) \,, \quad \forall\, \boldsymbol{v}_r \in \mathcal{V}_r \,. \tag{4.5}$$

Testing (4.5) with any $\boldsymbol{v}_r \in \mathcal{V}_r \cap \mathrm{Null}(\mathcal{B}')$, observe that $(\boldsymbol{v}_{h,r}, \boldsymbol{v}_r)_{\mathcal{V}} = 0$. Therefore, $\boldsymbol{v}_{h,r} \perp \mathcal{V}_r \cap \mathrm{Null}(\mathcal{B}')$. Clearly, $\boldsymbol{v}_r \perp \mathcal{V}_r^{\perp} \cap \mathrm{Null}(\mathcal{B}')$. The result follows immediately. □

4.2. **Error in the quantity of interest.** Here, we demonstrate that the duality seen in (3.3) also holds for the errors in the DPG and DPG* problems, even though the two discrete solutions exist in spaces of *different* dimensions. This is a fundamental result and may be viewed as a new generalization of [7, Equation 1.8] or [45, Theorem 3.1] which holds for discrete trial and test spaces, $\mathcal{U}_h$ and $\mathcal{V}_r$, of unequal dimension.

**Theorem 4.2.** *The following identity holds for all $\boldsymbol{w}_h \in \mathcal{U}_h$:*

$$G(\boldsymbol{u}^\star - \boldsymbol{u}_{h,r}) = b(\boldsymbol{u}^\star - \boldsymbol{u}_{h,r}, \boldsymbol{v}^\star - \Theta_r \boldsymbol{w}_h) = b(\boldsymbol{u}^\star - \boldsymbol{w}_h, \boldsymbol{v}^\star - \boldsymbol{v}_{h,r}) = \ell(\boldsymbol{v}^\star - \boldsymbol{v}_{h,r}) \,. \tag{4.6}$$

*Proof.* Due to Galerkin orthogonality in (4.1), $b(\boldsymbol{u}^\star - \boldsymbol{u}_{h,r}, \Theta_r \boldsymbol{w}_h) = 0$, for any $\boldsymbol{w}_h \in \mathcal{U}_h$. Likewise,

$$G(\boldsymbol{u}^\star - \boldsymbol{u}_{h,r}) = b(\boldsymbol{u}^\star - \boldsymbol{u}_{h,r}, \boldsymbol{v}^\star) = b(\boldsymbol{u}^\star - \boldsymbol{u}_{h,r}, \boldsymbol{v}^\star - \Theta_r \boldsymbol{w}_h) = b(\boldsymbol{u}^\star - \boldsymbol{u}_{h,r}, \boldsymbol{v}^\star - \boldsymbol{v}_{h,r}) \,.$$

Alternatively, due to Galerkin orthogonality in (4.3), $b(\boldsymbol{w}_h, \boldsymbol{v}^\star - \boldsymbol{v}_{h,r}) = 0$, for any $\boldsymbol{w}_h \in \mathcal{U}_h$. Thus,

$$\ell(\boldsymbol{v}^\star - \boldsymbol{v}_{h,r}) = b(\boldsymbol{u}^\star, \boldsymbol{v}^\star - \boldsymbol{v}_{h,r}) = b(\boldsymbol{u}^\star - \boldsymbol{w}_h, \boldsymbol{v}^\star - \boldsymbol{v}_{h,r}) = b(\boldsymbol{u}^\star - \boldsymbol{u}_{h,r}, \boldsymbol{v}^\star - \boldsymbol{v}_{h,r}) \,.$$

□

The following corollary, which follows immediately from the special case of Theorem 4.2 where $\mathcal{V}_r = \mathcal{V}$, delivers the well-known duality in the error in the symmetric coercive setting [60] and, in particular, in the least-squares setting, $a(\boldsymbol{u}, \boldsymbol{w}) = (\mathcal{L}\boldsymbol{u}, \mathcal{L}\boldsymbol{w})_\Omega$ and $F(\boldsymbol{w}) = (f, \boldsymbol{w})_\Omega$ [18].

**Corollary 4.3.** *The following identity holds:*

$$G(\boldsymbol{u}^\star - \boldsymbol{u}_h^{\mathrm{opt}}) = a(\boldsymbol{u}^\star - \boldsymbol{u}_h^{\mathrm{opt}}, \boldsymbol{\omega}^\star - \boldsymbol{\omega}_h^{\mathrm{opt}}) = F(\boldsymbol{\omega}^\star - \boldsymbol{\omega}_h^{\mathrm{opt}}) \,.$$

Note that $b(\boldsymbol{u}^\star - \boldsymbol{u}_{h,r}, \boldsymbol{v}^\star - \boldsymbol{v}_{h,r}) \leq |||\boldsymbol{u}^\star - \boldsymbol{u}_{h,r}|||_{\mathcal{U}} \|\boldsymbol{v}^\star - \boldsymbol{v}_{h,r}\|_{\mathcal{V}}$. Invoking (3.5), there is also the immediate corollary—a notable generalization of [45, Corollary 3.2]—which may be used in *a priori* analysis to demonstrate the accelerated convergence rate of the QoI error, $e_{\mathrm{QoI}}$. Because the *a priori* analysis of DPG* methods would be a considerable detour, it has been postponed to a follow-up article [49].

**Corollary 4.4.** *Define $e_{\mathrm{QoI}} = b(\boldsymbol{u}^\star - \boldsymbol{u}_{h,r}, \boldsymbol{v}^\star - \boldsymbol{v}_{h,r})$. The following crude upper bounds hold:*

$$|e_{\mathrm{QoI}}| \leq M \cdot \begin{cases} \|\boldsymbol{u}^\star - \boldsymbol{u}_{h,r}\|_{\mathcal{U}} \displaystyle\inf_{\boldsymbol{w}_h \in \mathcal{U}_h} \|\boldsymbol{v}^\star - \Theta_r \boldsymbol{w}_h\|_{\mathcal{V}} \\ \|\boldsymbol{v}^\star - \boldsymbol{v}_{h,r}\|_{\mathcal{V}} \displaystyle\inf_{\boldsymbol{w}_h \in \mathcal{U}_h} \|\boldsymbol{u}^\star - \boldsymbol{w}_h\|_{\mathcal{U}} \end{cases} \tag{4.7}$$

We also pose the following theorem which is the motivating factor behind the adaptive strategy introduced in the next section.

**Theorem 4.5.** *Define* $e_{\mathrm{QoI}} = b(\boldsymbol{u}^\star - \boldsymbol{u}_{h,r}, \boldsymbol{v}^\star - \boldsymbol{v}_{h,r})$. *The following crude upper bound holds:*

$$|e_{\mathrm{QoI}}| \leq \gamma^{-1} \|\mathcal{B}\boldsymbol{u}_{h,r} - \ell\|_{\mathcal{V}'} \|\mathcal{B}'\boldsymbol{v}_{h,r} - G\|_{\mathcal{U}'} \,. \tag{4.8}$$

**Remark 4.1.** The simplicity of (4.7) and (4.8) is actually surprising when compared to analogous upper bounds for other mixed methods [59]. Notably, these bounds do not involve the auxiliary functions $\boldsymbol{\psi}_{h,r}$ and $\boldsymbol{\omega}_{h,r}$.

*Proof of Theorem 4.5.* Observe that $|e_{\mathrm{QoI}}| = b(\boldsymbol{u}^\star - \boldsymbol{u}_{h,r}, \boldsymbol{v}^\star - \boldsymbol{v}_{h,r}) \leq |||\boldsymbol{u}^\star - \boldsymbol{u}_{h,r}|||_{\mathcal{U}} \|\boldsymbol{v}^\star - \boldsymbol{v}_{h,r}\|_{\mathcal{V}}$. Therefore, by Lemma 3.2,

$$|e_{\mathrm{QoI}}| \leq \|\mathcal{B}\boldsymbol{u}_{h,r} - \ell\|_{\mathcal{V}'} (\|\mathcal{P}(\boldsymbol{v}^\star - \boldsymbol{v}_{h,r})\|_{\mathcal{V}}^2 + |||\mathcal{B}'\boldsymbol{v}_{h,r} - G|||_{\mathcal{U}'}^2)^{1/2} \,.$$

Recall that $\boldsymbol{v}^\star$,$\boldsymbol{v}_{h,r} \perp \mathrm{Null}(\mathcal{B}')$, by Theorems 3.1 and 4.1. Therefore, $\mathcal{P}(\boldsymbol{v}^\star - \boldsymbol{v}_{h,r}) = 0$ and the result follows by observing that $|||\cdot|||_{\mathcal{U}'} \leq \gamma^{-1} \|\cdot\|_{\mathcal{U}'}$. □

## 5. GOAL-ORIENTED ADAPTIVE MESH REFINEMENT

In the literature, there are a number of established goal-oriented adaptive mesh refinement (GMR) and error estimation strategies. The most common general approach to goal-oriented error estimation is the dual-weighted residual (DWR) method [7, Section 5], which requires the solution of a dual problem on an enriched test space or additional post-processing of the discrete dual solution.[1] This method relies on a functional relationship to the error like (4.6) which can be localized to the element level and directly used in element marking. Most other competitive strategies instead deliver very useful upper and lower bounds on the QoI error [54, 60, 62]. In these strategies, marking is either directly driven by the bounding error estimates or by another estimate entirely.

Indeed, some such marking strategies, which happen to use two independent error estimators to bound a very crude QoI error estimate similar to (4.8), can be proved to deliver optimal GMR convergence rates [6, 37, 48, 55]. Meanwhile, even though many of the other GMR strategies work extremely well in practice, constructing rigorous proofs of optimal convergence rates with them has remained largely unprofitable. The approach to GMR we have taken is similar to these recent strategies and also [60] in that it also employs two independent error estimators, $\eta$ and $\eta^*$, for both the primal and dual (i.e. DPG and DPG*) methods, respectively.

### 5.1. Refinement indicators and adaptive mesh refinement algorithms.

Recall theorem 4.5. In this manuscript, GMR involves the following repetitive sequence of actions: (1) compute the approximate solutions $\boldsymbol{u}_{h,r}$ and $\boldsymbol{v}_{h,r}$; (2) construct residual-based *a posteriori* error estimators

$$\eta \approx \|\mathcal{B}\boldsymbol{u}_{h,r} - \ell\|_{\mathcal{V}'} \qquad \text{and} \qquad \eta^* \approx \|\mathcal{B}'\boldsymbol{v}_{h,r} - G\|_{\mathcal{U}'} \,,$$

and decompose them into element-wise components, called *refinement indicators*,

$$\eta^2 = \sum_{K \in \mathcal{T}} (\eta_K)^2 \qquad \text{and} \qquad (\eta^*)^2 = \sum_{K \in \mathcal{T}} (\eta_K^*)^2 \,;$$

(3) mark suitable elements for refinement using a marking convention influenced by contributions of both sets of refinement indicators $\{\eta_K\}_{K\in\mathcal{T}}$ and $\{\eta_K^*\}_{K\in\mathcal{T}}$; (4) refine all marked elements with the intention of driving the upper bound (4.8) to zero with an exceptional rate.[2] This general AMR strategy is summarized in Algorithm 1.

The key difference between GMR and solution-adaptive mesh refinement (SMR), which has predominately driven DPG AMR past studies, is in the marking strategy (step (4)). We have used a so-called "greedy" strategy in our experiments (see Section 8) instead of a more mathematically sound Dörfler-influenced strategy [33]. Upon the selection of a refinement parameter $0 < \theta < 1$, this marking strategy is defined as follows: calculate $\eta_{\max} = \max_{K\in\mathcal{T}} \widetilde{\eta}_K$ and mark all elements $K \in \mathcal{T}$ such that

$$\theta\, \eta_{\max} \leq \widetilde{\eta}_K \,. \tag{5.1}$$

Here, if SMR or GMR is being considered, then $\widetilde{\eta}_K := \eta_K$ or $\widetilde{\eta}_K := \eta_K \cdot \eta_K^*$, respectively.

[1]Although the dual problem (4.3) is also posed on an enriched test space $\dim(\mathcal{V}_r) > \dim(\mathcal{U}_h)$, the same stiffness matrix as in the primal problem (4.1) is ultimately used for the global solve. Therefore, the complexity of both problems is identical. This is not the case for the DWR method when the test space is enriched.

[2]At this step, a minimal set of additional elements are also refined, beyond the set of marked elements, to ensure that the new mesh has only one level of hanging nodes [31, Chapter 3].

**Algorithm 1** Adaptive mesh refinement

**Input:** initial mesh $\mathcal{T}$, marking strategy, complexity tolerance $\texttt{TOL}$.
**while** $|\mathcal{T}| < \texttt{TOL}$ **do**
(1) Solve for $\boldsymbol{u}_{h,r}$ and $\boldsymbol{v}_{h,r}$ on $\mathcal{T}$.
(2) Compute the two sets of refinement indicators $\{\eta_K\}_{K\in\mathcal{T}}$ and $\{\eta_K^*\}_{K\in\mathcal{T}}$.
(3) Mark elements for refinement $\mathcal{M}\subseteq\mathcal{T}$, as dictated by *the marking strategy*.
(4) Refine all marked elements $K\in\mathcal{M}$ and construct a new mesh $\mathcal{T}$.
**return** approximate solution $\boldsymbol{u}_{h,r}$, approximate QoI $G(\boldsymbol{u}_{h,r})$.

We will only consider one class of refinement indicators $\{\eta_K\}_{K\in\mathcal{T}}$ for the DPG problem (see Section 6). This class coincides with the same well-established, standard energy norm refinement indicators that have been used for SMR in many previous DPG studies [13, 16, 29, 38, 51, 65, 73]. Alternatively, we will consider three different classes of refinement indicators $\{\eta_K^*\}_{K\in\mathcal{T}}$ for the DPG* problem. Each of these are derived in Section 7.

## 6. A POSTERIORI ERROR ESTIMATION FOR DPG METHODS

In this section, we consider a specific well-studied implicit error estimator for the term

$$|||\boldsymbol{u}^\star - \boldsymbol{u}_{h,r}|||_{\mathcal{U}} = \|\mathcal{B}\boldsymbol{u}_{h,r} - \ell\|_{\mathcal{V}'}$$

in (4.8). Namely, we consider $\eta = \eta_r(\boldsymbol{u}_{h,r})$, where

$$\eta_r(\boldsymbol{u}) = \sup_{\boldsymbol{v}_r\in\mathcal{V}_r} \frac{b(\boldsymbol{u},\boldsymbol{v}_r) - \ell(\boldsymbol{v}_r)}{\|\boldsymbol{v}_r\|_{\mathcal{V}}}\,, \quad \forall\, \boldsymbol{u}\in\mathcal{U}\,. \tag{6.1}$$

From now on, we will also denote this error estimator simply as $\eta_r(\boldsymbol{u}) = \|\mathcal{B}\boldsymbol{u} - \ell\|_{\mathcal{V}_r'}$.

The DPG error estimator (6.1) has been well-studied and analyzed in the literature. For expanded discussions on it, we refer the interested reader to [13, 29]. Before we present the main result of this section, we summarize the most important properties of $\eta_r(\boldsymbol{u})$.

6.1. **Properties of the error estimator.** Observe that $\eta_r(\boldsymbol{u}_{h,r}) = \|\mathcal{B}\boldsymbol{u}_{h,r} - \ell\|_{\mathcal{V}_r'} = \|\boldsymbol{\psi}_{h,r}\|_{\mathcal{V}}$. Let Assumption 2 hold for $\mathcal{V}$ and the subspace $\mathcal{V}_r \subseteq \mathcal{V}$. Specifically, $\mathcal{V} = \prod_{K\in\mathcal{T}} \mathcal{V}_K$ and $\mathcal{V}_r = \prod_{K\in\mathcal{T}} \mathcal{V}_{K,r}$, where $\mathcal{V}_{K,r}\subseteq\mathcal{V}_K$, for all $K\in\mathcal{T}$. Then, as a consequence of Lemma 3.1,

$$\eta_r(\boldsymbol{u})^2 = \sum_{K\in\mathcal{T}} \|\mathcal{B}\boldsymbol{u} - \ell\|^2_{\mathcal{V}'_{K,r}}\,.$$

Define $\eta_K = \|\mathcal{B}\boldsymbol{u}_{h,r} - \ell\|_{\mathcal{V}'_{K,r}} = \|\boldsymbol{\psi}_{h,r}\|_{\mathcal{V}_K}$. If $\mathcal{V}$ is localizable, then each refinement indicator $\eta_K$ can be computed locally [66].

Recall that $\boldsymbol{\psi}^\star = 0$, by Theorem 3.1. Therefore,

$$|||\boldsymbol{u}^\star - \boldsymbol{u}_h^{\text{opt}}|||_{\mathcal{U}} = \|\mathcal{B}(\boldsymbol{u}^\star - \boldsymbol{u}_h^{\text{opt}})\|_{\mathcal{V}'} = \|\boldsymbol{\psi}_h^{\text{opt}}\|_{\mathcal{V}}\,.$$

When $\boldsymbol{u}_h^{\text{opt}}$ is replaced by $\boldsymbol{u}_{h,r}$ and $\boldsymbol{\psi}_h^{\text{opt}}$ is replaced by $\boldsymbol{\psi}_{h,r}$, the expression noticeably changes:

$$|||\boldsymbol{u}^\star - \boldsymbol{u}_{h,r}|||^2_{\mathcal{U}} = \|\,\mathcal{R}_{\mathcal{V}}\,\boldsymbol{\psi}_{h,r} - \mathcal{B}\boldsymbol{u}_{h,r} + \ell\|^2_{\mathcal{V}'} + \eta_r(\boldsymbol{u}_{h,r})^2\,.$$

This is an immediate consequence of the following lemma.

**Lemma 6.1.** *Let Assumption1 hold. Let $\mathcal{V}_r$ be any closed subspace of $\mathcal{V}$. Define $\boldsymbol{\psi}_{h,r}\in\mathcal{V}_r$ to be the unique solution of $\langle\mathcal{R}_{\mathcal{V}}\,\boldsymbol{\psi}_{h,r},\boldsymbol{v}_r\rangle = \langle\mathcal{B}\boldsymbol{u}_h - \ell,\boldsymbol{v}_r\rangle$, for all $\boldsymbol{v}_r\in\mathcal{V}_r$. Then*

$$\|\mathcal{B}(\boldsymbol{u}^\star - \boldsymbol{u}_h)\|^2_{\mathcal{V}'} = \|\,\mathcal{R}_{\mathcal{V}}\,\boldsymbol{\psi}_{h,r} - \mathcal{B}\boldsymbol{u}_h + \ell\|^2_{\mathcal{V}'} + \|\boldsymbol{\psi}_{h,r}\|^2_{\mathcal{V}}\,.$$

*Proof.* Observe that

$$\begin{aligned}\|\mathcal{B}(\boldsymbol{u}^\star - \boldsymbol{u}_h)\|_{\mathcal{V}'}^2 &= \|(\mathcal{R}_\mathcal{V}\,\boldsymbol{\psi}_{h,r} - \mathcal{R}_\mathcal{V}\,\boldsymbol{\psi}_{h,r}) - \mathcal{B}\boldsymbol{u}_h + \ell\|_{\mathcal{V}'}^2 \\ &= \|\,\mathcal{R}_\mathcal{V}\,\boldsymbol{\psi}_{h,r} - \mathcal{B}\boldsymbol{u}_h + \ell\|_{\mathcal{V}'}^2 - 2\langle \mathcal{R}_\mathcal{V}\,\boldsymbol{\psi}_{h,r} - \mathcal{B}\boldsymbol{u}_h + \ell, \boldsymbol{\psi}_{h,r}\rangle + \|\,\mathcal{R}_\mathcal{V}\,\boldsymbol{\psi}_{h,r}\|_{\mathcal{V}'}^2 \\ &= \|\,\mathcal{R}_\mathcal{V}\,\boldsymbol{\psi}_{h,r} - \mathcal{B}\boldsymbol{u}_h + \ell\|_{\mathcal{V}'}^2 + \|\boldsymbol{\psi}_{h,r}\|_{\mathcal{V}}^2 \,.\end{aligned}$$

□

6.2. **Reliability and efficiency.** The primary result in this section is an improvement on the main theorem in [13] in the case that the Fortin operator $\Pi_r$ is a bounded projection (see Assumption 4) To prove the result, we will require Theorem 6.2 and Theorem 6.3. The former has become reasonably well-known in the literature and many different proofs for it are given in [72]. The latter is perhaps much less well known. We provide a proof of this result separately, in Appendix A.

**Theorem 6.2** (Complementary projections)**.** *Let* $\Pi$ *be any bounded projection on a Hilbert space* $\mathcal{W}$, $\Pi\circ\Pi = \Pi$. *Then*

$$\|\Pi\| = \|1 - \Pi\| \,.$$

**Theorem 6.3** (Pythagoras)**.** *Let* $\mathcal{W}$ *be a Hilbert space and* $\mathcal{W}_0 \subseteq \mathcal{W}$ *be a nontrivial closed subspace. Let* $\mathcal{P} : \mathcal{W} \to \mathcal{W}_0$ *be the orthogonal projection onto* $\mathcal{W}_0$ *and let* $\Pi : \mathcal{W} \to \mathcal{W}_0$ *be any other bounded projection onto* $\mathcal{W}_0 = \Pi(\mathcal{W})$. *Then*

$$\|\Pi - \mathcal{P}\|^2 + 1 = \|\Pi\|^2 \,.$$

**Theorem 6.4** (Improved *a posteriori* error estimates for DPG)**.** *Let Assumptions 1 and 4 hold. Let* $\ell \in \mathrm{Range}(\mathcal{B})$, $\boldsymbol{u}^\star = \mathcal{B}^{-1}\ell$, *and* $\boldsymbol{u}_h \in \mathcal{U}_h$ *be arbitrary. Denote the best approximation of* $\boldsymbol{u}^\star$ *in* $\mathcal{U}_h$ *as*

$$\boldsymbol{u}_h^\star = \operatorname*{arg\,min}_{\boldsymbol{u}_h\in\mathcal{U}_h} \|\boldsymbol{u}^\star - \boldsymbol{u}_h\|_{\mathcal{U}}.$$

*Then the computable residual* $\eta_r(\boldsymbol{u}) = \|\ell - \mathcal{B}\boldsymbol{u}\|_{\mathcal{V}_r'}$ *and the data approximation error* $\mathrm{osc}(\ell) = \|\ell\circ(1-\Pi_r)\|_{\mathcal{V}'}$ *satisfy*

$$\gamma^2\|\boldsymbol{u}^\star - \boldsymbol{u}_h\|_{\mathcal{U}}^2 \leq \eta_r(\boldsymbol{u}_h)^2 + \Big(\eta_r(\boldsymbol{u}_h)\sqrt{\|\Pi_r\|^2 - 1} + \mathrm{osc}(\ell)\Big)^2. \tag{6.2}$$

$$\eta_r(\boldsymbol{u}_h) \leq M\|\boldsymbol{u}^\star - \boldsymbol{u}_h\|_{\mathcal{U}}\,, \tag{6.3}$$

$$\textit{and} \qquad \mathrm{osc}(\ell) \leq M\|\Pi_r\|\|\boldsymbol{u}^\star - \boldsymbol{u}_h^\star\|_{\mathcal{U}}\,. \tag{6.4}$$

*Remark* 6.1. The reliability bound in (6.2) is a new version of that found in [13], with the additional assumption that $\Pi_r \circ \Pi_r = \Pi_r$. Note that if $\mathrm{osc}(\ell) = 0$ then

$$\gamma\|\boldsymbol{u}^\star - \boldsymbol{u}_h\|_{\mathcal{U}} \leq |\!|\!|\boldsymbol{u}^\star - \boldsymbol{u}_h|\!|\!|_{\mathcal{U}} \leq \|\Pi_r\|\eta_r(\boldsymbol{u}_h)\,.$$

Moreover, if $\|\Pi_r\| = 1$—that is, if $\Pi_r$ is an orthogonal projection—then

$$\gamma^2\|\boldsymbol{u}^\star - \boldsymbol{u}_h\|_{\mathcal{U}}^2 \leq |\!|\!|\boldsymbol{u}^\star - \boldsymbol{u}_h|\!|\!|_{\mathcal{U}}^2 \leq \eta_r(\boldsymbol{u}_h)^2 + \mathrm{osc}(\ell)^2.$$

*Remark* 6.2. Reproducing the remarks in [13, Therorem 2.1], each of the inequalities in theorem 6.4 can be demonstrated to be sharp.

*Proof of Theorem 6.4.* To arrive at (6.3), simply observe that

$$\eta_r(\boldsymbol{u}_h) = \|\mathcal{B}(\boldsymbol{u}^\star - \boldsymbol{u}_h)\|_{\mathcal{V}_r'} \leq \|\mathcal{B}(\boldsymbol{u}^\star - \boldsymbol{u}_h)\|_{\mathcal{V}'} \leq \|\mathcal{B}\|\|\boldsymbol{u}^\star - \boldsymbol{u}_h\|_{\mathcal{U}}\,.$$

To arrive at (6.4), first let $\boldsymbol{v} \in \mathcal{V}$ where $\|\boldsymbol{v}\| = 1$ be arbitrary. Then, by (4.4) and Galerkin orthogonality,

$$\ell(\boldsymbol{v} - \Pi_r\boldsymbol{v}) = b(\boldsymbol{u}^\star, \boldsymbol{v} - \Pi_r\boldsymbol{v}) = b(\boldsymbol{u}^\star - \boldsymbol{u}_h, \boldsymbol{v} - \Pi_r\boldsymbol{v}) \leq \|\mathcal{B}\|\|1 - \Pi_r\|\|\boldsymbol{u}^\star - \boldsymbol{u}_h\|_{\mathcal{U}}\,,$$

for every $\boldsymbol{u}_h \in \mathcal{U}_h$. The result then follows from Assumption 4 and Theorem 6.2: $\|1 - \Pi_r\| = \|\Pi_r\|$.

The reliability bound is much more subtle. Begin by defining $\widetilde{\mathcal{V}}_r = \mathrm{Range}(\Pi_r) \subseteq \mathcal{V}_r$. Now, define $\widetilde{\boldsymbol{\psi}}_{h,r} \in \widetilde{\mathcal{V}}_r$ to be the unique solution of $\langle \mathcal{R}_\mathcal{V}\, \widetilde{\boldsymbol{\psi}}_{h,r}, \widetilde{\boldsymbol{v}}_r\rangle = \langle \ell - \mathcal{B}\boldsymbol{u}_h, \widetilde{\boldsymbol{v}}_r\rangle$, for all $\widetilde{\boldsymbol{v}}_r \in \widetilde{\mathcal{V}}_r$. Recall that $\gamma\|\boldsymbol{u}\|_\mathcal{U} \leq \|\mathcal{B}\boldsymbol{u}\|_{\mathcal{V}'}$, for all $\boldsymbol{u} \in \mathcal{U}$. Therefore, because $\widetilde{\mathcal{V}}_r$ is a closed subspace of $\mathcal{V}$, note that

$$\gamma^2\|\boldsymbol{u}^\star - \boldsymbol{u}_h\|_\mathcal{U}^2 \leq \|\mathcal{B}(\boldsymbol{u}^\star - \boldsymbol{u}_h)\|_{\mathcal{V}'}^2 = \|\,\mathcal{R}_\mathcal{V}\,\widetilde{\boldsymbol{\psi}}_{h,r} + \mathcal{B}\boldsymbol{u}_h - \ell\|_{\mathcal{V}'}^2 + \|\widetilde{\boldsymbol{\psi}}_{h,r}\|_\mathcal{V}^2\,, \tag{6.5}$$

by Lemma 6.1.

Define $\mathcal{P}_r : \mathcal{V} \to \widetilde{\mathcal{V}}_r$ to be the orthogonal projection onto the range of $\Pi_r$. Because $\widetilde{\boldsymbol{\psi}}_{h,r} \in \widetilde{\mathcal{V}}_r$, notice that

$$(\widetilde{\boldsymbol{\psi}}_{h,r}, \boldsymbol{v} - \Pi_r\boldsymbol{v})_\mathcal{V} = (\widetilde{\boldsymbol{\psi}}_{h,r}, \mathcal{P}_r(\boldsymbol{v} - \Pi_r\boldsymbol{v}))_\mathcal{V} = (\widetilde{\boldsymbol{\psi}}_{h,r}, \mathcal{P}_r\boldsymbol{v} - \Pi_r\boldsymbol{v})_\mathcal{V}\,, \quad \forall\, \boldsymbol{v} \in \mathcal{V}\,. \tag{6.6}$$

Now, observe that

$$\begin{aligned}
\|\,\mathcal{R}_\mathcal{V}\,\widetilde{\boldsymbol{\psi}}_{h,r} + \mathcal{B}\boldsymbol{u}_h - \ell\|_{\mathcal{V}'} &= \sup_{\boldsymbol{v}\in\mathcal{V}} \frac{(\widetilde{\boldsymbol{\psi}}_{h,r}, \boldsymbol{v})_\mathcal{V} + b(\boldsymbol{u}_h, \boldsymbol{v}) - \ell(\boldsymbol{v})}{\|\boldsymbol{v}\|_\mathcal{V}} \\
&= \sup_{\boldsymbol{v}\in\mathcal{V}} \frac{(\widetilde{\boldsymbol{\psi}}_{h,r}, \boldsymbol{v} - \Pi_r\boldsymbol{v})_\mathcal{V} + b(\boldsymbol{u}_h, \boldsymbol{v} - \Pi_r\boldsymbol{v}) - \ell(\boldsymbol{v} - \Pi_r\boldsymbol{v})}{\|\boldsymbol{v}\|_\mathcal{V}} \\
&= \sup_{\boldsymbol{v}\in\mathcal{V}} \frac{(\widetilde{\boldsymbol{\psi}}_{h,r}, \mathcal{P}_r\boldsymbol{v} - \Pi_r\boldsymbol{v})_\mathcal{V} - \ell(\boldsymbol{v} - \Pi_r\boldsymbol{v})}{\|\boldsymbol{v}\|_\mathcal{V}} \\
&\leq \|\widetilde{\boldsymbol{\psi}}_{h,r}\|_\mathcal{V}\|\mathcal{P}_r - \Pi_r\| + \|\ell \circ (1 - \Pi_r)\|_{\mathcal{V}'} \\
&= \|\widetilde{\boldsymbol{\psi}}_{h,r}\|_\mathcal{V}\sqrt{\|\Pi_r\|^2 - 1} + \mathrm{osc}(\ell)\,.
\end{aligned} \tag{6.7}$$

In the third line we have used (6.6) and (4.4), meanwhile, in the final line we have used Theorem 6.3. Finally, observe that

$$\|\widetilde{\boldsymbol{\psi}}_{h,r}\|_\mathcal{V} = \sup_{\widetilde{\boldsymbol{v}}_r\in\widetilde{\mathcal{V}}_r} \frac{b(\boldsymbol{u}, \widetilde{\boldsymbol{v}}_r) - \ell(\widetilde{\boldsymbol{v}}_r)}{\|\widetilde{\boldsymbol{v}}_r\|_\mathcal{V}} \leq \sup_{\boldsymbol{v}_r\in\mathcal{V}_r} \frac{b(\boldsymbol{u}, \boldsymbol{v}_r) - \ell(\boldsymbol{v}_r)}{\|\boldsymbol{v}_r\|_\mathcal{V}} = \eta_r(\boldsymbol{u}_h)\,.$$

With this observation in hand, (6.5) and (6.7) complete the proof. □

## 7. A POSTERIORI ERROR ESTIMATION FOR DPG* METHODS

In this section, we present three different strategies to approximate the term

$$|||\boldsymbol{v}^\star - \boldsymbol{v}_{h,r}|||_\mathcal{V} = \|\mathcal{B}'\boldsymbol{v}_{h,r} - G\|_{\mathcal{U}'}$$

in (4.8). In attempt at brevity, only the explicit error estimator $\eta^{*,\mathrm{expl}}(\boldsymbol{v}_{h,r})$ in Section 7.2 is rigorously justified.

7.1. **Functional setting and quantities of interest.** With the exception of Section 7.3, this section only considers ultraweak variational formulations, $b(\boldsymbol{u}, \boldsymbol{v}) = (u, \mathcal{L}^*\boldsymbol{v})_\Omega + \langle \hat{u}, \boldsymbol{v}\rangle_{\partial\mathcal{T}}$. In turn, we may consider general quantities of interest defined by goal functionals

$$G(\boldsymbol{u}) = G_\Omega(u) + \hat{G}(\hat{u})\,, \quad \forall\, \boldsymbol{u} = (u, \hat{u}) \in \mathcal{U} = L^2(\Omega) \times \hat{\mathcal{U}}\,, \tag{7.1}$$

where $G_\Omega \in L^2(\Omega)'$ and $\hat{G} \in \hat{\mathcal{U}}'$. Clearly, the first term can always be expressed as $G_\Omega(u) = (u, g)_\Omega$, for some fixed $g \in L^2(\Omega)$. Additionally, we may consider $\hat{G}(\hat{u}) = \langle \hat{u}, \hat{g}\rangle_{\partial\mathcal{T}}$, for some fixed $\hat{g} \in \mathcal{V}$.

In our experiments, we focused only on smooth linear functionals which are also well-posed on the underlying *unbroken* trial space. In particular, goal functionals

$$G(\boldsymbol{u}) = (u, g)_\Omega + \langle \hat{g}_{\partial\Omega}, \hat{u}|_{\partial\Omega}\rangle\,, \quad \forall\, \boldsymbol{u} = (u, \hat{u}) \in L^2(\Omega) \times \hat{\mathcal{U}}\,, \tag{7.2}$$

where $g \in L^2(\Omega)$ and $\hat{g}_{\partial\Omega} \in L^2(\partial\Omega)$ sufficiently regular that the second term $\langle \hat{g}_{\partial\Omega}, \hat{u}|_{\partial\Omega}\rangle$ was well-defined. In such scenarios, any well-defined functional $\langle \hat{g}_{\partial\Omega}, \hat{u}|_{\partial\Omega}\rangle$ can be identified with a functional $\hat{G}(\hat{u}) = \langle \hat{u}, \hat{g}\rangle_{\partial\mathcal{T}}$, where $\hat{g}$ has additional (unbroken) regularity, (see Remark 7.3).

7.2. **Explicit error estimators.** Let $G(\boldsymbol{u}) = (u, g)_\Omega + \langle \hat{u}, \hat{g} \rangle_{\partial\mathcal{T}}$. By Lemma 3.1, observe that

$$\|\mathcal{B}'\boldsymbol{v} - G\|^2_{\mathcal{U}'} = \left( \sup_{\boldsymbol{u}\in\mathcal{U}} \frac{(u, \mathcal{L}^*\boldsymbol{v} - g)_\Omega + \langle \hat{u}, \boldsymbol{v} - \hat{g} \rangle_{\partial\mathcal{T}}}{\left(\|u\|^2_\Omega + \|\hat{u}\|^2_{\hat{\mathcal{U}}}\right)^{1/2}} \right)^2 = \|\mathcal{L}^*\boldsymbol{v} - g\|^2_\Omega + \sup_{\hat{u}\in\hat{\mathcal{U}}} \frac{\langle \hat{u}, \boldsymbol{v} - \hat{g} \rangle^2_{\partial\mathcal{T}}}{\|\hat{u}\|^2_{\hat{\mathcal{U}}}} . \tag{7.3}$$

The first term on the right-hand side corresponds to the local $L^2$-residuals, while the second measures the jumps of $\boldsymbol{v}$ across the element interfaces.

Let $\mathcal{F}$ be the set of all faces in the mesh skeleton and $h_F$ be the diameter of each face $F \in \mathcal{F}$. For appropriate functions $\boldsymbol{v}_r$, if $\hat{g}$ is smooth enough, then it can be shown that the jump term above is equivalent to a sum of weighted norms over the element faces:

$$\sup_{\hat{u}\in\hat{\mathcal{U}}} \frac{\langle \hat{u}, \boldsymbol{v}_r - \hat{g} \rangle^2_{\partial\mathcal{T}}}{\|\hat{u}\|^2_{\hat{\mathcal{U}}}} \approx \sum_{F\in\mathcal{F}} w_F(h_F) \| [\![ \boldsymbol{v}_r - \hat{g} ]\!] \|^2_{\hat{\mathcal{H}}(F)} , \tag{7.4}$$

where, relative to any fixed element $K^+ \in \mathcal{T}$, $[\![ \cdot ]\!]|_F$ denotes the inter-element jump across the face $F \in \mathcal{F}$, $F \subseteq \partial K^+$. Generally, $w_F : \mathbb{R}_{>0} \to \mathbb{R}_{\geq 0}$ and the face-dependent norm $\| \cdot \|_{\hat{\mathcal{H}}(F)}$ in (7.4) is stronger than the natural norm induced by the pairing $\langle \cdot, \cdot \rangle_{\partial\mathcal{T}}$. Both this norm and the weighting function $w_F$ must be determined in a case-by-case basis. Theorem 7.5 illustrates a common scenario where $\mathcal{V} = H^1(\mathcal{T}) \times \boldsymbol{H}(\mathrm{div}, \mathcal{T})$ in $\mathbb{R}^3$.

Depending on the location of edge, often either $w_F(h_F) = h_F^\theta$ or $w_F(h_F) = 0$, where $\theta \in \mathbb{R}$ is constant throughout all faces in the mesh. In this scenario, (7.3), (7.4), and the shape-regularity of the mesh motivate the error estimator

$$(\eta^{*,\mathrm{expl}})^2 = \sum_{K\in\mathcal{T}} \Big( \|\mathcal{L}^*\boldsymbol{v}_{h,r} - g\|^2_K + w_F(h_K) \sum_{F\in\mathcal{F}: F\subseteq\partial K} \| [\![ \boldsymbol{v}_{h,r} - \hat{g} ]\!] \|^2_{\hat{\mathcal{H}}(F)} \Big) . \tag{7.5}$$

7.3. **Implicit error estimators.** Recall (3.13). Namely,

$$\begin{cases} (\boldsymbol{v}_h^{\mathrm{opt}}, \boldsymbol{v})_{\mathcal{V}} - b(\boldsymbol{\omega}_h^{\mathrm{opt}}, \boldsymbol{v}) = 0 , & \forall\, \boldsymbol{v} \in \mathcal{V} , \\ b(\boldsymbol{u}_h, \boldsymbol{v}_h^{\mathrm{opt}}) = G(\boldsymbol{u}_h) , & \forall\, \boldsymbol{u}_h \in \mathcal{U}_h . \end{cases}$$

Define a particular representation of the residual in the ideal DPG* solution, $\boldsymbol{\varepsilon}_h^{\mathrm{opt}} \in \mathcal{U}$, as the unique solution of the saddle-point problem

$$a(\boldsymbol{\varepsilon}_h^{\mathrm{opt}}, \boldsymbol{u}) = G(\boldsymbol{u}) - b(\boldsymbol{u}, \boldsymbol{v}_h^{\mathrm{opt}}) , \quad \forall\, \boldsymbol{u} \in \mathcal{U} .$$

Clearly, $b(\boldsymbol{\varepsilon}_h^{\mathrm{opt}}, \Theta\boldsymbol{u}_h) = a(\boldsymbol{\varepsilon}_h^{\mathrm{opt}}, \boldsymbol{u}_h) = 0$, for all $\boldsymbol{u}_h \in \mathcal{U}_h$ and therefore $(\boldsymbol{\varepsilon}_h^{\mathrm{opt}}, \boldsymbol{v}_h^{\mathrm{opt}}) \in \mathcal{U} \times \Theta(\mathcal{U}_h)$ is the unique solution of

$$\begin{cases} a(\boldsymbol{\varepsilon}_h^{\mathrm{opt}}, \boldsymbol{u}) + b(\boldsymbol{u}, \boldsymbol{v}_h^{\mathrm{opt}}) = G(\boldsymbol{u}) , & \forall\, \boldsymbol{u} \in \mathcal{U} , \\ b(\boldsymbol{\varepsilon}_h^{\mathrm{opt}}, \boldsymbol{v}_h) = 0 , & \forall\, \boldsymbol{v}_h \in \Theta(\mathcal{U}_h) . \end{cases} \tag{7.6}$$

Moreover, $|||\boldsymbol{\varepsilon}_h^{\mathrm{opt}}|||_{\mathcal{U}} = |||\mathcal{A}^{-1}(\mathcal{B}'\boldsymbol{v}_h^{\mathrm{opt}} - G)|||_{\mathcal{U}} = |||\mathcal{B}'\boldsymbol{v}_h^{\mathrm{opt}} - G|||_{\mathcal{U}'}$. The saddle-point problem (7.6) is similar to the system analyzed in the saddle-point least-squares method [4,5]. Likewise, define a representation of the residual in the practical DPG* solution, $\boldsymbol{\varepsilon}_{h,r}$, as the unique solution of

$$a(\boldsymbol{\varepsilon}_{h,r}, \boldsymbol{u}) = G(\boldsymbol{u}) - b(\boldsymbol{u}, \boldsymbol{v}_{h,r}) , \quad \forall\, \boldsymbol{u} \in \mathcal{U} , \tag{7.7}$$

and deduce that $|||\boldsymbol{\varepsilon}_{h,r}|||_{\mathcal{U}} = |||\mathcal{B}'\boldsymbol{v}_{h,r} - G|||_{\mathcal{U}'}$.[3]

Notice that $\boldsymbol{\varepsilon}_{h,r}$ cannot be computed exactly. In order to estimate $\boldsymbol{\varepsilon}_{h,r}$, we may pose (7.7) on an enriched trial space $\mathcal{U}_H \subseteq \mathcal{U}$ containing the original trial space $\mathcal{U}_h \subsetneq \mathcal{U}_H$. Furthermore, because the inner product $a : \mathcal{U} \times \mathcal{U} \to \mathbb{R}$ is defined using the inverse of the full Riesz map $\mathcal{R}_{\mathcal{V}}^{-1}$, we use a further enriched test space

[3]Notice that $(\boldsymbol{\varepsilon}_{h,r}, \boldsymbol{u})_{\mathcal{U}} = G(\boldsymbol{u}) - b(\boldsymbol{u}, \boldsymbol{v}_{h,r})$ is also a convenient defining equation. However, the trial norm will be somewhat formulation dependent and, in the ultraweak setting, $\|(u, \hat{u})\|^2_{\mathcal{U}} = \|u\|^2_\Omega + \|\hat{u}\|_{\hat{\mathcal{U}}}$ presents its own natural challenges.

$\mathcal{V}_R \subseteq \mathcal{V}$ to generate a natural approximation, $a_R : \mathcal{U}_H \times \mathcal{U}_H \to \mathbb{R}$. In other words, in order to compute a practical estimate the solution $\boldsymbol{\varepsilon}_{h,r}$ to (7.7), we propose solving a discrete variational problem

$$a_R(\boldsymbol{\varepsilon}, \boldsymbol{u}_H) = G(\boldsymbol{u}_H) - b(\boldsymbol{u}_H, \boldsymbol{v}_{h,r}), \quad \forall \boldsymbol{u}_H \in \mathcal{U}_H , \tag{7.8}$$

where $a_R(\boldsymbol{u}_H, \widetilde{\boldsymbol{u}}_H) = \langle \mathcal{B}\boldsymbol{u}_H, \mathcal{R}_{\mathcal{V}_R}^{-1} \mathcal{B}\widetilde{\boldsymbol{u}}_H \rangle$, for all $\boldsymbol{u}_H, \widetilde{\boldsymbol{u}}_H \in \mathcal{U}_H$. The function $\boldsymbol{\varepsilon}$ defined above is a computationally feasible estimate of $\boldsymbol{\varepsilon}_{h,r}$ and

$$\eta^{*,\mathrm{impl}} = a_R(\boldsymbol{\varepsilon}, \boldsymbol{\varepsilon})^{1/2} ,$$

can be used as a (global) implicit error estimator.

In general, (7.8) a global problem. For instance, in the broken test space setting where $\mathcal{V}_R = \prod_{K\in\mathcal{T}} \mathcal{V}_{R,K}$,

$$a_R(\boldsymbol{u}_H, \widetilde{\boldsymbol{u}}_H) = \sum_{K\in\mathcal{T}} \langle \mathcal{B}\boldsymbol{u}_H, \mathcal{R}_{\mathcal{V}_{R,K}}^{-1} \mathcal{B}\widetilde{\boldsymbol{u}}_H \rangle , \quad \forall \boldsymbol{u}_H, \widetilde{\boldsymbol{u}}_H \in \mathcal{U}_H ,$$

where each term $a_{R,K}(\boldsymbol{u}_H, \widetilde{\boldsymbol{u}}_H) = \langle \mathcal{B}\boldsymbol{u}_H, \mathcal{R}_{\mathcal{V}_{R,K}}^{-1} \mathcal{B}\widetilde{\boldsymbol{u}}_H \rangle$ has the symmetric structure of an inner product but is not necessary definite. Define $\mathcal{U}_{H,K} = \{\boldsymbol{u}|_K \mid \boldsymbol{u} \in \mathcal{U}_H\}$. Consequently, the auxiliary local problem for the restricted residual estimate, $\boldsymbol{\varepsilon}_K \in \mathcal{U}_{H,K}$, namely

$$a_{R,K}(\boldsymbol{\varepsilon}_K, \boldsymbol{u}_H) = G(\boldsymbol{u}_H) - b(\boldsymbol{u}_H, \boldsymbol{v}_{h,r}), \quad \forall \boldsymbol{u}_H \in \mathcal{U}_{H,K} , \tag{7.9}$$

cannot be expected to have a unique solution. However, intentionally removing the kernel of $\mathcal{B}|_{\mathcal{U}_{H,K}}$ from the solution to (7.9) will introduce uniqueness and allow for the efficient local construction of the necessary refinement indicators:

$$\eta_K^{*,\mathrm{impl}} = a_{R,K}(\boldsymbol{\varepsilon}_K, \boldsymbol{\varepsilon}_K)^{1/2} , \quad \forall K \in \mathcal{T} . \tag{7.10}$$

The precise method for removing the kernel of $\mathcal{B}|_{\mathcal{U}_{H,K}}$ from (7.9) differentiates each category of possible implicit DPG* error estimator. There are several different classical strategies to choose from [74], each involving their own subtle analysis. In our experiments, we simply introduced artificial homogeneous boundary conditions (see Section 8.2). This decision resulted in a polluted estimate of $\boldsymbol{\varepsilon}_K$, however, it did not appear to adversely affect the performance of resulting adaptive mesh refinement algorithm.

7.4. **An additional ad hoc class of refinement indicators.** As a bona fide *a posteriori* error estimator for DPG* methods, the ad hoc error estimator proposed in this section, $\eta^{*,\mathrm{a.h.}}$, has very poor properties. Therefore, we do *not* recommend it for actual error estimation in DPG* problems. Moreover, we only recommend it for the limited set of goal functionals (7.1) where $\hat{G} = 0$. Nevertheless, we have observed that this "error estimator" can deliver very favorable pre-asymptotic results with the GMR algorithm in Section 5.

We begin with a number of observations. First, for every goal $G \in \mathcal{U}'$, there exists a unique function $\boldsymbol{g} \in \mathcal{U}$ such that

$$(\boldsymbol{g}, \boldsymbol{u})_{\mathcal{U}} = G(\boldsymbol{u}), \quad \forall \boldsymbol{u} \in \mathcal{U}.$$

Namely, $\boldsymbol{g} = \mathcal{R}_{\mathcal{U}}^{-1} G$ is the Riesz representation of $G$ induced by the trial norm $\|\cdot\|_{\mathcal{U}}$. Here, recall that the opt-notation exists to remind the reader that the equality only holds in the idealized setting, $\|\cdot\|_{\mathcal{V}} \overset{\mathrm{opt}}{=} |||\cdot|||_{\mathcal{V}}$. As we saw in (3.16), $\mathcal{R}_{\mathcal{U}} \overset{\mathrm{opt}}{=} \mathcal{B}' \mathcal{R}_{\mathcal{V}}^{-1} \mathcal{B} = \mathcal{A}$. Recall that $\mathcal{A}\boldsymbol{\omega}^\star = G$. Therefore, $\boldsymbol{\omega}^\star \overset{\mathrm{opt}}{=} \boldsymbol{g}$. Moreover, by theorem 3.1 and (3.19),

$$\|\boldsymbol{v}^\star - \boldsymbol{v}_h^{\mathrm{opt}}\|_{\mathcal{V}} = |||\mathcal{B}'\boldsymbol{v}_h^{\mathrm{opt}} - G|||_{\mathcal{U}'} = |||\mathcal{B}' \mathcal{R}_{\mathcal{V}}^{-1} \mathcal{B}\boldsymbol{\omega}_h^{\mathrm{opt}} - G|||_{\mathcal{U}'} \overset{\mathrm{opt}}{=} \| \mathcal{R}_{\mathcal{U}} \boldsymbol{\omega}_h^{\mathrm{opt}} - G\|_{\mathcal{U}'} \overset{\mathrm{opt}}{=} \|\boldsymbol{g} - \boldsymbol{\omega}_h^{\mathrm{opt}}\|_{\mathcal{U}} .$$

We also note another identity. Let $\boldsymbol{u} \in \mathcal{U}$ be arbitrary and define $\boldsymbol{\psi} = \mathcal{R}_{\mathcal{V}}^{-1}(\mathcal{B}\boldsymbol{u} - \ell)$. With the optimal test norm, note that

$$G(\boldsymbol{u}^\star - \boldsymbol{u}) = G(\mathcal{B}^{-1} \mathcal{R}_{\mathcal{V}} \mathcal{R}_{\mathcal{V}}^{-1} \mathcal{B}(\boldsymbol{u}^\star - \boldsymbol{u})) \overset{\mathrm{opt}}{=} -G(\mathcal{R}_{\mathcal{U}}^{-1} \mathcal{B}'\boldsymbol{\psi}) \overset{\mathrm{opt}}{=} -b(\boldsymbol{g}, \boldsymbol{\psi}) .$$

This indicates an explicit relationship between error in the QoI and the residual of the primal problem.

To now begin to derive the error estimator, first consider the broken ultraweak setting, $b(\boldsymbol{u}, \boldsymbol{v}) = (u, \mathcal{L}^*\boldsymbol{v})_\Omega + \langle \hat{u}, \boldsymbol{v}\rangle_{\partial\mathcal{T}}$, and assume that $G(\boldsymbol{u}) = G_\Omega(u)$, where $G_\Omega(u) = (g, u)_\Omega$. In this case, observe that $\boldsymbol{g} = (g, 0)$, where $g = \mathcal{R}_{L^2(\Omega)}^{-1} G_\Omega$. With the corresponding adjoint graph norm (3.23), instead of the optimal test norm (3.22), it can be shown that $\boldsymbol{\omega}^\star = \boldsymbol{g} + \alpha^2\widetilde{\boldsymbol{\omega}}$, for a fixed function $\widetilde{\boldsymbol{\omega}}$, and so can be made arbitrarily close to $\boldsymbol{g}$ by adjusting the constant $\alpha > 0$. Lemma 7.1, adopted from [42, Lemma 7], makes this fact precise in the present setting.

Recall (4.3) and define $\boldsymbol{\omega}_{h,r} = (\omega_{h,r}, \hat{\omega}_{h,r})$. Now, observe that $\|\boldsymbol{g} - \boldsymbol{\omega}_{h,r}\|^2_{\mathcal{U}} = \|g - \omega_{h,r}\|^2_\Omega + \|\hat{\omega}_{h,r}\|^2_{\hat{\mathcal{U}}}$. Because the norm $\|\cdot\|_{\hat{\mathcal{U}}}$ is difficult to estimate numerically, we ignore the second and define the ad hoc DPG* "error estimator" $(\eta^{*,\text{a.h.}})^2 = \sum_{K\in\mathcal{T}}(\eta_K^{*,\text{a.h.}})^2$, where

$$\eta_K^{*,\text{a.h.}} = \|g - \omega_{h,r}\|_K\,, \quad \forall\, K \in \mathcal{T}\,. \tag{7.11}$$

*Remark* 7.1. Reasonable behavior with this error estimator will only be found when the influence function $\boldsymbol{\omega}^\star$ coming from the solution of the dual problem (3.4) is approximately equal to $\boldsymbol{g} = \mathcal{R}_{\mathcal{U}}^{-1}\,G$. This is possible when the adjoint graph norm is used, $G(\boldsymbol{u}) = G_\Omega(u)$, and $\alpha$ is sufficiently close to 0. Additionally, it is important to note that

$$\big|\|g - \omega^\star\|_\Omega - \|\omega^\star - \omega_{h,r}\|_\Omega\big| \leq \|g - \omega_{h,r}\|_\Omega\,.$$

Therefore, even though an *a priori* estimate can show that $\|\omega^\star - \omega_{h,r}\|_\Omega \to 0$ as the mesh is refined, $\eta^{*,\text{a.h.}} \not\to 0$.

**Lemma 7.1.** *Let $b(\boldsymbol{u},\boldsymbol{v}) = (u, \mathcal{L}^*\boldsymbol{v})_\Omega + \langle \hat{u}, \boldsymbol{v}\rangle_{\partial\mathcal{T}}$, and assume that $G(\boldsymbol{u}) = G_\Omega(u) + \hat{G}(\hat{u})$, where $G_\Omega(u) = (g,u)$ and $\hat{G} \in \hat{\mathcal{U}}'$. Let Assumption 1 hold. Endow $\mathcal{V}$ with the adjoint graph norm $\|\boldsymbol{v}\|^2_{\mathcal{L}^*,\alpha} = \|\mathcal{L}^*\boldsymbol{v}\|^2_\Omega + \alpha^2\|\boldsymbol{v}\|^2_\Omega$. Then $\boldsymbol{\omega}^\star = (g,0) + \alpha^2\widetilde{\boldsymbol{\omega}}$, where and $\widetilde{\boldsymbol{\omega}} = (\widetilde{\omega}, \hat{\widetilde{\omega}}) \in L^2(\Omega)\times\hat{\mathcal{U}}$, is the unique solution of*

$$(\widetilde{\omega}, \mathcal{L}^*\boldsymbol{v})_\Omega + \langle \hat{\widetilde{\omega}}, \boldsymbol{v}\rangle_{\partial\mathcal{T}} = (\boldsymbol{v}^\star, \boldsymbol{v})_\Omega\,, \quad \forall\, \boldsymbol{v} \in \mathcal{V}\,.$$

*Proof.* With the adjoint graph norm, (3.15) can be rewritten

$$\begin{cases} (\mathcal{L}^*\boldsymbol{v}^\star, \mathcal{L}^*\boldsymbol{v})_\Omega + \alpha^2(\boldsymbol{v}^\star,\boldsymbol{v})_\Omega - (\omega^\star, \mathcal{L}^*\boldsymbol{v})_\Omega - \langle\hat{\omega}^\star, \boldsymbol{v}\rangle_{\partial\mathcal{T}} = 0\,, & \forall\, \boldsymbol{v} \in \mathcal{V}\,,\\ (u, \mathcal{L}^*\boldsymbol{v}^\star)_\Omega = (g,u)_\Omega\,, & \forall\, u \in L^2(\Omega)\,,\\ \langle \hat{u}, \boldsymbol{v}^\star\rangle_{\partial\mathcal{T}} = \hat{G}(\hat{u})\,, & \forall\, \hat{u} \in \hat{\mathcal{U}}\,. \end{cases}$$

Recall that $\boldsymbol{v}^\star \perp \mathrm{Null}(\mathcal{B}')$, by theorem 3.1. Therefore, the bottom two equations, together expressible as $\mathcal{B}'\boldsymbol{v}^\star = G$, uniquely determine $\boldsymbol{v}^\star$. Testing the second equation with $\mathcal{L}^*\boldsymbol{v}$ and substituting into the first, we then find that

$$(\omega^\star, \mathcal{L}^*\boldsymbol{v})_\Omega + \langle\hat{\omega}^\star, \boldsymbol{v}\rangle_{\partial\mathcal{T}} = \alpha^2(\boldsymbol{v}^\star,\boldsymbol{v})_\Omega + (g, \mathcal{L}^*\boldsymbol{v})_\Omega\,, \quad \forall\, \boldsymbol{v} \in \mathcal{V}\,,$$

which uniquely determines $\boldsymbol{\omega}^\star$, as necessary. □

7.5. **Reliability of the explicit error estimator.** In this section, we give a specific example of the explicit error estimator (7.5) which applies to a large class of ultraweak variational formulations of second-order elliptic boundary value problems and immediately generalizes to many ultraweak formulations of more complicated problems, such as many coming from continuum models such as linearized elasticity and Stokes flow. We then prove that this error estimator is *reliable* in the sense that it bounds the DPG* residual $\|\mathcal{B}'\boldsymbol{v}_h^{\text{opt}} - G\|_{\mathcal{U}'}$ from above by a mesh-independent constant. We note that our analysis here only considers the case of tetrahedral meshes even though our numerical experiments are carried out with hexahedral elements.

For all Lipschitz subdomains $K \subseteq \Omega$, let $\mathrm{tr}^K : H^1(K) \to H^{1/2}(\partial K)$, $\mathrm{tr}_n^K : \boldsymbol{H}(\mathrm{div},K) \to H^{-1/2}(\partial K)$, $\mathrm{tr}_t^K : \boldsymbol{H}(\mathbf{curl},K) \to H^{-1/2}(\mathrm{div},\partial K)$ and $\mathrm{tr}_T^K : \boldsymbol{H}(\mathbf{curl},K) \to H^{-1/2}(\mathbf{Curl},\partial K)$ be the canonical trace operators, where $\mathrm{tr}^K u = u|_{\partial K}$, $\mathrm{tr}_n^K \boldsymbol{\sigma} = \boldsymbol{\sigma}|_{\partial K}\cdot\boldsymbol{n}_K$, $\mathrm{tr}_T^K\boldsymbol{\sigma} = \boldsymbol{n}_K \times \boldsymbol{\sigma}|_{\partial K}$, and $\mathrm{tr}_t^K\boldsymbol{\sigma} = (\boldsymbol{n}_K\times\boldsymbol{\sigma}|_{\partial K})\times\boldsymbol{n}_K$ for smooth functions. Here, by $\boldsymbol{n}_K$, we mean the outer unit normal vector on $\partial K$.

Let $\Gamma_{\mathrm{D}}, \Gamma_{\mathrm{N}} \subseteq \partial\Omega$ be relatively open. Using the trace operators given above, define the Hilbert spaces $H^1_{\mathrm{D}}(\Omega) = \{u \in H^1(\Omega) : (\mathrm{tr}^\Omega u)|_{\Gamma_{\mathrm{D}}} = 0\}$, $\boldsymbol{H}_{\mathrm{N}}(\mathrm{div},\Omega) = \{\boldsymbol{\sigma} \in \boldsymbol{H}(\mathrm{div},\Omega) : (\mathrm{tr}_n^\Omega\boldsymbol{\sigma})|_{\Gamma_{\mathrm{N}}} = 0\}$, and $\boldsymbol{H}_{\mathrm{N}}(\mathbf{curl},\Omega) = \{\boldsymbol{\sigma} \in \boldsymbol{H}(\mathbf{curl},\Omega) : (\mathrm{tr}_t^\Omega\boldsymbol{\sigma})|_{\Gamma_{\mathrm{N}}} = 0\}$. Define the mesh-dependent trace operators $\mathrm{tr} = \prod_{K\in\mathcal{T}}\mathrm{tr}^K$, $\mathrm{tr}_n = \prod_{K\in\mathcal{T}}\mathrm{tr}_n^K$, and $\mathrm{tr}_T = \prod_{K\in\mathcal{T}}\mathrm{tr}_T^K$. Then define the interface spaces $H^{1/2}_{\mathrm{D}}(\mathcal{S}) = \mathrm{tr}(H^1_{\mathrm{D}}(\Omega))$ and $H^{-1/2}_{\mathrm{N}}(\mathcal{S}) = \mathrm{tr}_n(\boldsymbol{H}_{\mathrm{N}}(\mathrm{div},\Omega))$ and the broken test spaces $H^1(\mathcal{T}) = \prod_{K\in\mathcal{T}} H^1(K)$ and $\boldsymbol{H}(\mathrm{div},\mathcal{T}) =$ $\prod_{K\in\mathcal{T}} \boldsymbol{H}(\mathrm{div}, K)$. Note that an interface function $\hat{u} \in H^{1/2}_{\mathrm{D}}(\mathcal{S})$ is entirely single-valued, meanwhile, an interface function $\hat{\sigma} \in H^{-1/2}_{\mathrm{N}}(\mathcal{S})$ is only single-valued *up to its sign*. Finally, using the $H^{1/2}(\partial K)$–$H^{-1/2}(\partial K)$ duality pairing $\langle\cdot,\cdot\rangle_{\partial K}$, define

$$\langle \hat{u}, \boldsymbol{\tau}\cdot\boldsymbol{n}\rangle_{\partial\mathcal{T}} = \sum_{K\in\mathcal{T}} \langle \hat{u}|_{\partial K}, \mathrm{tr}^K_n\, \boldsymbol{\tau}\rangle_{\partial K} \qquad \text{and} \qquad \langle \hat{\sigma}, v\rangle_{\partial\mathcal{T}} = \sum_{K\in\mathcal{T}} \langle \hat{\sigma}|_{\partial K}, \mathrm{tr}^K v\rangle_{\partial K}\,,$$

for all $\hat{u} \in H^{1/2}_{\mathrm{D}}(\mathcal{S})$, $\hat{\sigma} \in H^{-1/2}_{\mathrm{N}}(\mathcal{S})$, $v \in H^1(\mathcal{T})$, and $\boldsymbol{\tau} \in \boldsymbol{H}(\mathrm{div},\mathcal{T})$. Here, $\boldsymbol{n}$ is understood as the double-valued unit normal vector field on the mesh skeleton. It is important to remark that the interface spaces $H^{1/2}_{\mathrm{D}}(\mathcal{S})$ and $H^{-1/2}_{\mathrm{N}}(\mathcal{S})$ are equipped with special quotient norms which will be denoted $\|\hat{u}\|_{H^{1/2}(\partial\mathcal{T})}$ and $\|\hat{\sigma}\|_{H^{-1/2}(\partial\mathcal{T})}$, respectively. These and other similar norms are intensively analyzed in [14].

We will naturally consider scalar- and vector-valued functions which are inherently double-valued on any interior face $F = \partial K^+ \cap \partial K^-$. Reducing to the vector case first, consider a locally smooth, discontinuous function $\boldsymbol{\tau}$ with the double-valued restriction $\boldsymbol{\tau}|_F$. Now, identify the two individual branches of this restriction; $(\boldsymbol{\tau}|_F)_{K^+}$ and $(\boldsymbol{\tau}|_F)_{K^-}$. For any such function $\boldsymbol{\tau}$, define $[\![\boldsymbol{\tau}\cdot\boldsymbol{n}]\!]$ face-by-face as the normal jump across the mesh skeleton:

$$[\![\boldsymbol{\tau}\cdot\boldsymbol{n}]\!]\big|_F = \begin{cases} (\boldsymbol{\tau}|_F)_{K^+}\cdot\boldsymbol{n}_{K^+} + (\boldsymbol{\tau}|_F)_{K^-}\cdot\boldsymbol{n}_{K^-}\,, & \text{if } \exists K^-\in\mathcal{T} : \partial K^+\cap\partial K^- = F\,,\\ \boldsymbol{\tau}\cdot\boldsymbol{n}_{K^+}\,, & \text{otherwise}\,. \end{cases}$$

Notably, another vector trace we will require at one step in the proof is the tangential jump across the mesh skeleton:

$$[\![\boldsymbol{n}\times\boldsymbol{\tau}]\!]\big|_F = \begin{cases} \boldsymbol{n}_{K^+}\times(\boldsymbol{\tau}|_F)_{K^+} + \boldsymbol{n}_{K^-}\times(\boldsymbol{\tau}|_F)_{K^-}\,, & \text{if } \exists K^-\in\mathcal{T} : \partial K^+\cap\partial K^- = F\,,\\ \boldsymbol{n}_{K^+}\times\boldsymbol{\tau}\,, & \text{otherwise}\,. \end{cases}$$

Alternatively, for scalar-valued functions, we will only require the unsigned jump function

$$[\![v]\!]\big|_F = \begin{cases} |(v|_F)_{K^+} - (v|_F)_{K^-}|\,, & \text{if } \exists K^-\in\mathcal{T} : \partial K^+\cap\partial K^- = F\,,\\ |v|\,, & \text{otherwise}\,. \end{cases}$$

When posing error estimates using the jump functions above, it is natural to introduce terms like $[\![(\boldsymbol{\tau}-\widetilde{\boldsymbol{g}})\cdot\boldsymbol{n}]\!]$, where $\widetilde{\boldsymbol{g}}$ is an $\boldsymbol{H}(\mathrm{div},\Omega)$ extension of sufficiently regular boundary data $\hat{g}_n$. In such cases, we instead write the shorthand $[\![\boldsymbol{\tau}\cdot\boldsymbol{n}-\hat{g}_n]\!]$, which reminds the reader that the interior jump terms can be neglected:

$$[\![\boldsymbol{\tau}\cdot\boldsymbol{n}-\hat{g}_n]\!]\big|_F = \begin{cases} (\boldsymbol{\tau}|_F)_{K^+}\cdot\boldsymbol{n}_{K^+} + (\boldsymbol{\tau}|_F)_{K^-}\cdot\boldsymbol{n}_{K^-}\,, & \text{if } \exists K^-\in\mathcal{T} : \partial K^+\cap\partial K^- = F\,,\\ \boldsymbol{\tau}\cdot\boldsymbol{n}_{K^+} - \hat{g}_n\,, & \text{otherwise}\,. \end{cases}$$

A similar meaning is attributed to the jump function $[\![v_r - \hat{g}]\!]$.

Before proving Theorem 7.5, we require the following identities and theorems. The proof of the identities in Lemma 7.2 can be found in [25, Section 3.2]. Theorem 7.3 is a special case of [32, Lemma 5] for $n = 3$ and $k = 2$. Theorem 7.4 is a special case of [32, Lemma 6] which has important connections to the work many people, including Clément [20], Scott & Zhang [71], Schöberl [70], and Falk & Winther [36].

**Lemma 7.2.** *Let* $\boldsymbol{\tau}\in\boldsymbol{H}(\mathrm{div},\mathcal{T})$ *and* $v\in H^1(\mathcal{T})$. *Then*

$$\sup_{\hat{u}\in H^{1/2}_{\mathrm{D}}(\mathcal{S})} \frac{\langle \hat{u}, \boldsymbol{\tau}\cdot\boldsymbol{n}\rangle_{\partial\mathcal{T}}}{\|\hat{u}\|_{H^{1/2}(\partial\mathcal{T})}} = \sup_{u\in H^1_{\mathrm{D}}(\Omega)} \frac{\langle \mathrm{tr}\, u, \boldsymbol{\tau}\cdot\boldsymbol{n}\rangle_{\partial\mathcal{T}}}{\|u\|_{H^1(\Omega)}}$$

*and*

$$\sup_{\hat{\sigma}\in H^{-1/2}_{\mathrm{N}}(\mathcal{S})} \frac{\langle \hat{\sigma}, v\rangle_{\partial\mathcal{T}}}{\|\hat{\sigma}\|_{H^{-1/2}(\partial\mathcal{T})}} = \sup_{\boldsymbol{\sigma}\in\boldsymbol{H}_{\mathrm{N}}(\mathrm{div},\Omega)} \frac{\langle \mathrm{tr}_n\,\boldsymbol{\sigma}, v\rangle_{\partial\mathcal{T}}}{\|\boldsymbol{\sigma}\|_{\boldsymbol{H}(\mathrm{div},\Omega)}}\,.$$

**Theorem 7.3.** *Given any $\boldsymbol{\sigma} \in \boldsymbol{H}(\mathrm{div},\Omega)$, there exist $\boldsymbol{\varphi} \in \boldsymbol{H}^1(\Omega)$ and $\boldsymbol{\Phi} \in \boldsymbol{H}^1(\Omega)$ such that $\boldsymbol{\sigma} = \mathbf{curl}\,\boldsymbol{\varphi} + \boldsymbol{\Phi}$ and*

$$\|\boldsymbol{\varphi}\|_{\boldsymbol{H}^1(\Omega)} + \|\boldsymbol{\Phi}\|_{\boldsymbol{H}^1(\Omega)} \lesssim \|\boldsymbol{\sigma}\|_{\boldsymbol{H}(\mathrm{div},\Omega)}\,.$$

*In this situation, $\boldsymbol{\sigma} = \mathbf{curl}\,\boldsymbol{\varphi} + \boldsymbol{\Phi}$ is called a* regular decomposition *of $\boldsymbol{\sigma}$.*

**Theorem 7.4.** *For any face $F \in \mathcal{F}$, define $h_F$ to be its diameter and define $\Omega_F \subseteq \Omega$ to be the patch of elements $K \in \mathcal{T}$ neighboring it. Let $p \in \mathbb{N}_0$ and $\mathcal{P}^1_{\mathrm{c,D}}(\mathcal{T}) = \mathcal{P}^1(\mathcal{T}) \cap H^1_{\mathrm{D}}(\Omega)$, where $\mathcal{P}^1(\mathcal{T})$ denotes the space of $\mathcal{T}$-piecewise polynomials of degree less than or equal to $1$. Let $\mathcal{RT}^0_{\mathrm{N}}(\mathcal{T}) = \mathcal{RT}^0(\mathcal{T}) \cap \boldsymbol{H}_{\mathrm{N}}(\mathrm{div},\Omega)$, where $\mathcal{RT}^0(\mathcal{T})$ is the lowest-order Raviart-Thomas space over $\mathcal{T}$ [8]. Likewise, let $\mathcal{N}^0_{\mathrm{N}}(\mathcal{T}) = \mathcal{N}^0(\mathcal{T}) \cap \boldsymbol{H}_{\mathrm{N}}(\mathrm{div},\Omega)$, where $\mathcal{N}^0(\mathcal{T})$ is the lowest-order Nédélec space over $\mathcal{T}$.*

*There exist commuting operators $\Pi_{\mathbf{grad}} : H^1_{\mathrm{D}}(\Omega) \to \mathcal{P}^1_{\mathrm{c,D}}(\mathcal{T})$, $\Pi_{\mathrm{div}} : \boldsymbol{H}_{\mathrm{N}}(\mathrm{div},\Omega) \to \mathcal{RT}^0_{\mathrm{N}}(\mathcal{T})$, and $\Pi_{\mathbf{curl}} : \boldsymbol{H}_{\mathrm{N}}(\mathbf{curl},\Omega) \to \mathcal{N}^0_{\mathrm{N}}(\mathcal{T})$ such that $\Pi_{\mathrm{div}}\,\mathbf{curl}\,\boldsymbol{\varphi} = \mathbf{curl}\,\Pi_{\mathbf{curl}}\boldsymbol{\varphi}$,*

$$\|\,\mathrm{tr}(u - \Pi_{\mathbf{grad}}u)\|_F \lesssim h_F^{1/2}\|u\|_{H^1(\Omega_F)}\,, \quad \forall u \in H^1_{\mathrm{D}}(\Omega)\,,$$

*and, for all regular decompositions of $\boldsymbol{\sigma} \in \boldsymbol{H}(\mathrm{div},\Omega)$, $\boldsymbol{\sigma} = \mathbf{curl}\,\boldsymbol{\varphi} + \boldsymbol{\Phi}$,*

$$\|\,\mathrm{tr}_T(\boldsymbol{\varphi} - \Pi_{\mathbf{curl}}\boldsymbol{\varphi})\|_F + \|\,\mathrm{tr}_n(\boldsymbol{\Phi} - \Pi_{\mathrm{div}}\boldsymbol{\Phi})\|_F \lesssim h_F^{1/2}\|\boldsymbol{\sigma}\|_{\boldsymbol{H}(\mathrm{div},\Omega_F)}\,.$$

**Theorem 7.5.** *Assume that $\Omega \subseteq \mathbb{R}^3$ is a bounded Lipschitz domain and let $\mathcal{T}$ be a regular subdivision of $\Omega$. Assume that $\overline{\Gamma_{\mathrm{D}} \cup \Gamma_{\mathrm{N}}} = \partial\Omega$ are both relatively open in $\partial\Omega$, $\Gamma_{\mathrm{D}} \cap \Gamma_{\mathrm{N}} = \varnothing$, $\overline{\Gamma_{\mathrm{D}}} \cap \overline{\Gamma_{\mathrm{N}}}$ is Lipschitz, and $\Gamma_{\mathrm{D}} \neq \varnothing$. Define $\mathcal{F}_{\mathrm{N}}$ and $\mathcal{F}_{\mathrm{D}}$ to be the sets of interior faces in $\mathcal{T}$ union the faces in common with $\Gamma_{\mathrm{N}}$ and $\Gamma_{\mathrm{D}}$, respectively.*

*For all $\boldsymbol{u} = (u,\hat{u},\hat{\sigma}) \in \mathcal{U}$ and $\boldsymbol{v} = (v,\boldsymbol{\tau}) \in \mathcal{V}$, where $u \in L^2(\Omega)$, $\hat{u} \in H^{1/2}_{\mathrm{D}}(\mathcal{S})$, $\hat{\sigma} \in H^{-1/2}_{\mathrm{N}}(\mathcal{S})$, $v \in H^1(\mathcal{T})$, and $\boldsymbol{\tau} \in \boldsymbol{H}(\mathrm{div},\mathcal{T})$ define $\mathcal{U}$ and $\mathcal{V}$, let*

$$b(\boldsymbol{u},\boldsymbol{v}) = (u,\mathcal{L}^*\boldsymbol{v})_\Omega - \langle \hat{u}, \boldsymbol{\tau}\cdot\boldsymbol{n}\rangle_{\partial\mathcal{T}} - \langle \hat{\sigma}, v\rangle_{\partial\mathcal{T}} \qquad \text{and} \qquad G(\boldsymbol{u}) = (g,u)_\Omega + \langle \hat{g}_n, \hat{u}\rangle_{\partial\Omega} + \langle \hat{g}, \hat{\sigma}\rangle_{\partial\Omega}\,,$$

*where $g \in L^2(\Omega)$, $\hat{g}_n|_{\Gamma_{\mathrm{N}}} \in L^2(\Gamma_{\mathrm{N}})$, and $\hat{g}|_{\Gamma_{\mathrm{D}}} \in H^1_0(\Gamma_{\mathrm{D}})$. Here, without loss of generality, we implicitly assume that $\hat{g}_n$ and $\hat{g}$ are extended by zero throughout $\Gamma_{\mathrm{D}}$ and $\Gamma_{\mathrm{N}}$, respectively. If $\boldsymbol{v}_r = (v_r,\boldsymbol{\tau}_r) \in \mathcal{V}_r \subseteq \mathcal{V}$ satisfies*

$$b(\boldsymbol{v}_r,\boldsymbol{u}_h) = g(\boldsymbol{u}_h)\,, \quad \forall\, \boldsymbol{u}_h = (0,\hat{u}_h,\hat{\sigma}_h)\,, \tag{7.12}$$

*where $\hat{u}_h \in \mathrm{tr}\left(\mathcal{P}^1_{\mathrm{c,D}}(\mathcal{T})\right)$ and $\hat{\sigma}_h \in \mathrm{tr}\left(\mathcal{RT}^0_{\mathrm{N}}(\mathcal{T})\right)$, then $\|\mathcal{B}'\boldsymbol{v}_r - G\|^2_{\mathcal{U}'} \lesssim \eta^{*,\mathrm{expl}}$, where $\eta^{*,\mathrm{expl}} = \sum_{K\in\mathcal{T}} \eta_K^{*,\mathrm{expl}}$ and*

(7.13)

$$(\eta_K^{*,\mathrm{expl}})^2 = \|\mathcal{L}^*\boldsymbol{v}_r - g\|_K^2 + h_K\left(\sum_{F\in\mathcal{F}_{\mathrm{N}}:F\subseteq\partial K} \|[\![\boldsymbol{\tau}_r\cdot\boldsymbol{n} - \hat{g}_n]\!]\|_F^2 + \sum_{F\in\mathcal{F}_{\mathrm{D}}:F\subseteq\partial K} \|[\![v_r - \hat{g}]\!]\|^2_{H(\mathbf{Curl}_F,F)}\right),$$

*for each $K \in \mathcal{T}$. Here, for all faces $F$, $\|\cdot\|^2_{H(\mathbf{Curl}_F,F)} = \|\cdot\|_F^2 + \|\mathbf{Curl}_F\cdot\|_F^2$, where $\mathbf{Curl}_F = \boldsymbol{n}_F \times \mathbf{Grad}_F$ is the surface curl operator for the face $F$ and $\mathbf{Grad}_F$ is the surface gradient operator restricted to the face $F$ [56].*

*Remark* 7.2. Efficiency of the explicit estimator defined in Theorem 7.5—that is, a complementary lower bound on the residual of the dual problem, $\|\mathcal{B}'\boldsymbol{v}_r - G\|_{\mathcal{U}'}$—can be proven with Verfürth's bubble function technique [74]. We have chosen not to include the proof here because it is not absolutely necessary for our analysis.

*Remark* 7.3. Let $\widetilde{g} \in H^1(\Omega)$ be any $H^1$-extension of the function $\hat{g} \in H^1(\partial\Omega)$ into the domain $\Omega$. Redefine $\hat{g} = \mathrm{tr}(\widetilde{g})$. In this case, for any $\boldsymbol{H}(\mathrm{div})$-extension of $\hat{\sigma}$ into $\Omega$, say $\widetilde{\boldsymbol{\sigma}}$, observe that

$$\langle \hat{g},\hat{\sigma}\rangle_{\partial\Omega} = (\mathbf{grad}\,\widetilde{g},\widetilde{\boldsymbol{\sigma}})_\Omega + (\widetilde{g},\mathrm{div}\,\widetilde{\boldsymbol{\sigma}})_\Omega = \sum_{K\in\mathcal{T}} \left((\mathbf{grad}\,\widetilde{g},\widetilde{\boldsymbol{\sigma}})_K + (\widetilde{g},\mathrm{div}\,\widetilde{\boldsymbol{\sigma}})_K\right) = \langle \hat{g},\hat{\sigma}\rangle_{\partial\mathcal{T}}\,.$$

Therefore, without loss of generality, we may always substitute $\langle \hat{g},\hat{\sigma}\rangle_{\partial\Omega}$ with $\langle \hat{g},\hat{\sigma}\rangle_{\partial\mathcal{T}}$, where $\hat{g}$ is redefined as above. Similarly, we may always substitute $\langle \hat{g}_n,\hat{u}\rangle_{\partial\Omega}$ with $\langle \hat{g}_n,\hat{u}\rangle_{\partial\mathcal{T}}$.

*Proof of Theorem 7.5.* First of all, note that $\hat{g}_n \in \big(H_{\rm D}^{1/2}(\partial\Omega)\big)'$, and $\hat{g} \in \big(H_{\rm N}^{-1/2}(\partial\Omega)\big)'$. This is enough to demonstrate that $G \in \mathcal{U}'$.

Recall (7.3) and Remark 7.3. By Lemmas 3.1 and 7.2,

$$\|\mathcal{B}'\boldsymbol{v}_r - G\|_{\mathcal{U}'}^2 = \|\mathcal{L}^*\boldsymbol{v}_r - g\|_\Omega^2 + \sup_{u\in H^1_{\rm D}(\Omega)} \frac{\langle \operatorname{tr} u, \boldsymbol{\tau}_r\cdot\boldsymbol{n} - \hat{g}_n\rangle_{\partial\mathcal{T}}^2}{\|u\|_{H^1(\Omega)}^2} + \sup_{\boldsymbol{\sigma}\in\boldsymbol{H}_{\rm N}(\operatorname{div},\Omega)} \frac{\langle \operatorname{tr}_n \boldsymbol{\sigma}, v_r - \hat{g}\rangle_{\partial\mathcal{T}}^2}{\|\boldsymbol{\sigma}\|_{\boldsymbol{H}(\operatorname{div},\Omega)}^2}\,.$$

The first term in the right-hand side above can also be expressed $\|\mathcal{L}^*\boldsymbol{v}_r - g\|_\Omega^2 = \sum_{K\in\mathcal{T}} \|\mathcal{L}^*\boldsymbol{v}_r - g\|_K^2$. Therefore, we only have to handle the last two terms.

Invoking (7.12), observe that

$$\begin{aligned}
\langle \operatorname{tr} u, \boldsymbol{\tau}_r\cdot\boldsymbol{n} - \hat{g}_n\rangle_{\partial\mathcal{T}} &= \langle \operatorname{tr}(u - \Pi_{\mathbf{grad}}u), \boldsymbol{\tau}_r\cdot\boldsymbol{n} - \hat{g}_n\rangle_{\partial\mathcal{T}} \\
&= \sum_{F\in\mathcal{F}_{\rm N}} (\operatorname{tr}(u - \Pi_{\mathbf{grad}}u), [\![\boldsymbol{\tau}_r\cdot\boldsymbol{n} - \hat{g}_n]\!])_F \\
&\le \sum_{F\in\mathcal{F}_{\rm N}} \|\operatorname{tr}(u - \Pi_{\mathbf{grad}}u)\|_F \|[\![\boldsymbol{\tau}_r\cdot\boldsymbol{n} - \hat{g}_n]\!]\|_F \\
&\lesssim \sum_{F\in\mathcal{F}_{\rm N}} h_F^{1/2} \|u\|_{H^1(\Omega_F)} \|[\![\boldsymbol{\tau}_r\cdot\boldsymbol{n} - \hat{g}_n]\!]\|_F \\
&\le \Big(\sum_{F\in\mathcal{F}_{\rm N}} \|u\|_{H^1(\Omega_F)}^2\Big)^{1/2} \Big(\sum_{F\in\mathcal{F}_{\rm N}} h_F \|[\![\boldsymbol{\tau}_r\cdot\boldsymbol{n} - \hat{g}_n]\!]\|_F^2\Big)^{1/2}.
\end{aligned} \tag{7.14}$$

By the shape regularity of the mesh, $\sum_{F\in\mathcal{F}} \|u\|_{H^1(\Omega_F)}^2 \le \|u\|_{H^1(\Omega)}^2$, and so

$$\sup_{u\in H^1_{\rm D}(\Omega)} \frac{\langle \operatorname{tr} u, \boldsymbol{\tau}_r\cdot\boldsymbol{n} - \hat{g}_n\rangle_{\partial\mathcal{T}}^2}{\|u\|_{H^1(\Omega)}^2} \lesssim \sum_{F\in\mathcal{F}_{\rm N}} h_F \|[\![\boldsymbol{\tau}_r\cdot\boldsymbol{n} - \hat{g}_n]\!]\|_F^2 \lesssim \sum_{K\in\mathcal{T}} h_K \sum_{F\in\mathcal{F}_{\rm N}:F\subseteq\partial K} \|[\![\boldsymbol{\tau}_r\cdot\boldsymbol{n} - \hat{g}_n]\!]\|_F^2\,. \tag{7.15}$$

Let $\mathbf{curl}\,\boldsymbol{\varphi} + \boldsymbol{\Phi} = \boldsymbol{\sigma}$ be a regular decomposition; $\boldsymbol{\varphi}, \boldsymbol{\Phi} \in \boldsymbol{H}^1(\Omega)$ and $\mathbf{curl}\,\boldsymbol{\varphi}, \boldsymbol{\Phi} \in \boldsymbol{H}_{\rm N}(\operatorname{div},\Omega)$. Therefore,

$$\begin{aligned}
\sup_{\boldsymbol{\sigma}\in\boldsymbol{H}_{\rm N}(\operatorname{div},\Omega)} \frac{\langle \operatorname{tr}_n \boldsymbol{\sigma} - \hat{\sigma}_h, v_r - \hat{g}\rangle_{\partial\mathcal{T}}^2}{\|\boldsymbol{\sigma}\|_{\boldsymbol{H}(\operatorname{div},\Omega)}^2} &\lesssim \sup_{\substack{\boldsymbol{\varphi}\in\boldsymbol{H}^1(\Omega)\\ \mathbf{curl}\,\boldsymbol{\varphi}\in\boldsymbol{H}_{\rm N}(\operatorname{div},\Omega)}} \frac{\langle \operatorname{tr}_n \mathbf{curl}\,\boldsymbol{\varphi} - \hat{\sigma}_h, v_r - \hat{g}\rangle_{\partial\mathcal{T}}^2}{\|\boldsymbol{\varphi}\|_{\boldsymbol{H}^1(\Omega)}^2} \\
&\quad + \sup_{\boldsymbol{\Phi}\in\boldsymbol{H}^1(\Omega)\cap\boldsymbol{H}_{\rm N}(\operatorname{div},\Omega)} \frac{\langle \operatorname{tr}_n \boldsymbol{\Phi} - \hat{\sigma}_h, v_r - \hat{g}\rangle_{\partial\mathcal{T}}^2}{\|\boldsymbol{\Phi}\|_{\boldsymbol{H}^1(\Omega)}^2}\,.
\end{aligned}$$

Consider the second supremum first and set $\hat{\sigma}_h = \operatorname{tr}_n(\Pi_{\operatorname{div}}\boldsymbol{\Phi}) \in \operatorname{tr}\big(\mathcal{RT}^0_{\rm N}(\mathcal{T})\big)$. As in (7.14), we readily determine that

$$\langle \operatorname{tr}_n \boldsymbol{\Phi} - \hat{\sigma}_h, v_r - \hat{g}\rangle_{\partial\mathcal{T}} \le \sum_{F\in\mathcal{F}_{\rm D}} \|\operatorname{tr}_n(\boldsymbol{\Phi} - \Pi_{\operatorname{div}}\boldsymbol{\Phi})\|_F \|[\![v_r - \hat{g}]\!]\|_F \lesssim \|\boldsymbol{\Phi}\|_{\boldsymbol{H}^1(\Omega)} \Big(\sum_{F\in\mathcal{F}_{\rm D}} h_F \|[\![v_r - \hat{g}]\!]\|_F^2\Big)^{1/2}.$$

In the last inequality above, we have again invoked shape regularity. Therefore,

$$\sup_{\boldsymbol{\Phi}\in\boldsymbol{H}^1(\Omega)\cap\boldsymbol{H}_{\rm N}(\operatorname{div},\Omega)} \frac{\langle \operatorname{tr}_n(\boldsymbol{\Phi} - \Pi_{\operatorname{div}}\boldsymbol{\Phi}), v_r - \hat{g}\rangle_{\partial\mathcal{T}}^2}{\|\boldsymbol{\Phi}\|_{\boldsymbol{H}^1(\Omega)}^2} \lesssim \sum_{K\in\mathcal{T}} h_K \sum_{F\in\mathcal{F}_{\rm D}:F\subseteq\partial K} \|[\![v_r - \hat{g}]\!]\|_F^2\,. \tag{7.16}$$

Let $\widetilde{g} \in H^1(\Omega)$ be any $H^1$-extension of the function $\hat{g} \in H^1(\partial\Omega)$ into the domain $\Omega$, $\mathrm{tr}(\widetilde{g}) = \hat{g}$, and set $\hat{\sigma}_h = \mathrm{tr}_n(\Pi_{\mathrm{div}}\,\mathbf{curl}\,\boldsymbol{\varphi}) \in \mathrm{tr}\left(\mathcal{RT}^0_{\mathrm{N}}(\mathcal{T})\right)$. Recall that $\Pi_{\mathrm{div}}\,\mathbf{curl}\,\boldsymbol{\varphi} = \mathbf{curl}\,\Pi_{\mathbf{curl}}\boldsymbol{\varphi}$. Therefore,

$$\begin{aligned}
\langle \mathrm{tr}_n\,\mathbf{curl}\,\boldsymbol{\varphi} - \hat{\sigma}_h, v_r - \hat{g}\rangle_{\partial\mathcal{T}} &= \sum_{K\in\mathcal{T}} (\mathbf{curl}(\boldsymbol{\varphi} - \Pi_{\mathbf{curl}}\boldsymbol{\varphi}), \mathbf{grad}(v_r - \widetilde{g}))_K \\
&= \sum_{K\in\mathcal{T}} \langle \mathrm{tr}^K_T(\boldsymbol{\varphi} - \Pi_{\mathbf{curl}}\boldsymbol{\varphi}), \mathrm{tr}^K_t\,\mathbf{grad}(v_r - \widetilde{g})\rangle_{\partial K} \\
&= \sum_{F\in\mathcal{F}_{\mathrm{D}}} (\mathrm{tr}_T(\boldsymbol{\varphi} - \Pi_{\mathbf{curl}}\boldsymbol{\varphi}), [\![\boldsymbol{n}\times\mathbf{grad}(v_r - \widetilde{g})]\!])_F \\
&\leq \sum_{F\in\mathcal{F}_{\mathrm{D}}} \|\mathrm{tr}_T(\boldsymbol{\varphi} - \Pi_{\mathbf{curl}}\boldsymbol{\varphi})\|_F \|[\![\boldsymbol{n}\times\mathbf{grad}(v_r - \widetilde{g})]\!]\|_F \\
&\lesssim \|\boldsymbol{\varphi}\|_{\boldsymbol{H}^1(\Omega)} \bigg( \sum_{F\in\mathcal{F}_{\mathrm{D}}} h_F \|\mathbf{Curl}_F [\![v_r - \hat{g}]\!]\|_F^2 \bigg)^{1/2},
\end{aligned}$$

where, in the final line, we have used the fact that $\|[\![\boldsymbol{n}\times\mathbf{grad}(v_r - \widetilde{g})]\!]\|_F = \|\mathbf{Curl}_F[\![v_r - \hat{g}]\!]\|_F$. Finally,

$$\sup_{\substack{\boldsymbol{\varphi}\in\boldsymbol{H}^1(\Omega)\\ \mathbf{curl}\,\boldsymbol{\varphi}\in\boldsymbol{H}_{\mathrm{N}}(\mathrm{div},\Omega)}} \frac{\langle \mathrm{tr}_n\,\mathbf{curl}\,\boldsymbol{\varphi} - \hat{\sigma}_h, v_r - \hat{g}\rangle^2_{\partial\mathcal{T}}}{\|\boldsymbol{\varphi}\|^2_{\boldsymbol{H}^1(\Omega)}} \lesssim \sum_{K\in\mathcal{T}} h_K \sum_{F\in\mathcal{F}_{\mathrm{D}}:F\subseteq\partial K} \|\mathbf{Curl}_F[\![v_r - \hat{g}]\!]\|_F^2\,. \tag{7.17}$$

Altogether, (7.15)–(7.17) complete the proof. □

## 8. Numerical experiments

Throughout this section, we attempt to thoroughly demonstrate the efficacy of goal-oriented adaptive mesh refinement (GMR) using Algorithm 1 and the DPG and DPG* refinement indicators introduced in Sections 6 and 7, respectively. All of the following experiments pertain to the broken ultraweak formulation of Poisson's boundary value problem, which is analyzed in detail in [25]. Here, the corresponding ultraweak bilinear form is $b(\boldsymbol{u},\boldsymbol{v}) = (u, \mathrm{div}\,\boldsymbol{\tau})_\Omega + (\boldsymbol{\sigma}, \boldsymbol{\tau} + \mathbf{grad}\,v)_\Omega - \langle \hat{u}, \boldsymbol{\tau}\cdot\boldsymbol{n}\rangle_{\partial\mathcal{T}} - \langle\hat{\sigma}, v\rangle_{\partial\mathcal{T}}$, where $\boldsymbol{u} = (u,\boldsymbol{\sigma},\hat{u},\hat{\sigma}) \in \mathcal{U} = L^2(\Omega)\times \boldsymbol{L}^2(\Omega)\times H^{1/2}_{\mathrm{D}}(\mathcal{S})\times H^{-1/2}_{\mathrm{N}}(\mathcal{S})$ and $\boldsymbol{v} = (v,\boldsymbol{\tau}) \in \mathcal{V} = H^1(\mathcal{T})\times\boldsymbol{H}(\mathrm{div},\mathcal{T})$. Recall that $\eta_K = \|\mathcal{B}\boldsymbol{u}_{h,r} - \ell\|_{\mathcal{V}'_{K,r}}$, and $\eta^*_K$ are defined in (7.13), (7.10), and (7.11), depending on the specific refinement indicator type.

8.1. **Set-up.** To compare the GMR algorithm with the conventional solution-adaptive mesh refinement (SMR) algorithm (see Section 5), we used a manufactured solution with two regions of isolated steep and shallow gradients in the convex domain $\Omega = [0,4]\times[0,1]\times[0,1]$ (see Figure 8.2). The exact expression for this manufactured solution, $u^\star = u^{\mathrm{man}} \in C^\infty_0(\Omega)$, is given by

$$u^{\mathrm{man}}(x,y,z) = f(x/4)f(y)f(z)\,, \quad \text{where} \quad f(x) = x(1-x)\big((x/4) + (1-4x)^2\big)\,. \tag{8.1}$$

The remaining components of $\boldsymbol{u}^\star = (u^\star, \boldsymbol{\sigma}^\star, \hat{u}^\star, \hat{\sigma}^\star)$ can easily be derived from the expression above.

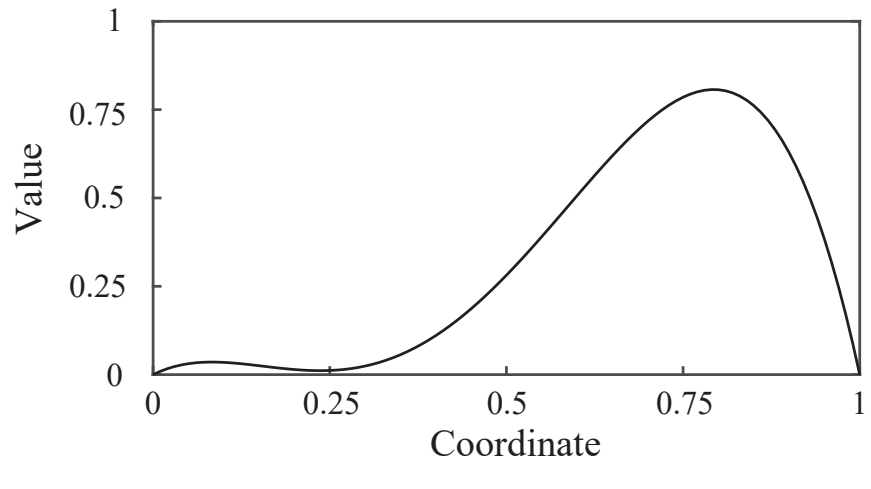

FIGURE 8.1. Graph of the function $f(x)$ in (8.1).

With such a solution, the expected behavior of an SMR strategy is to induce the majority of mesh refinements in the region with the highest gradients in the solution. That is, where the length scale of the solution is the smallest and the solution is the most difficult to resolve. Using the manufactured solution given in (8.1), this behavior is demonstrated in Figure 8.2 (c) from DPG SMR. From now on, SMR specifically refers to Algorithm 1 with the marking strategy (5.1), where $\widetilde{\eta}_K := \eta_K = \|\mathcal{B}\boldsymbol{u}_{h,r} - \ell\|_{\mathcal{V}'_{K,r}}$.

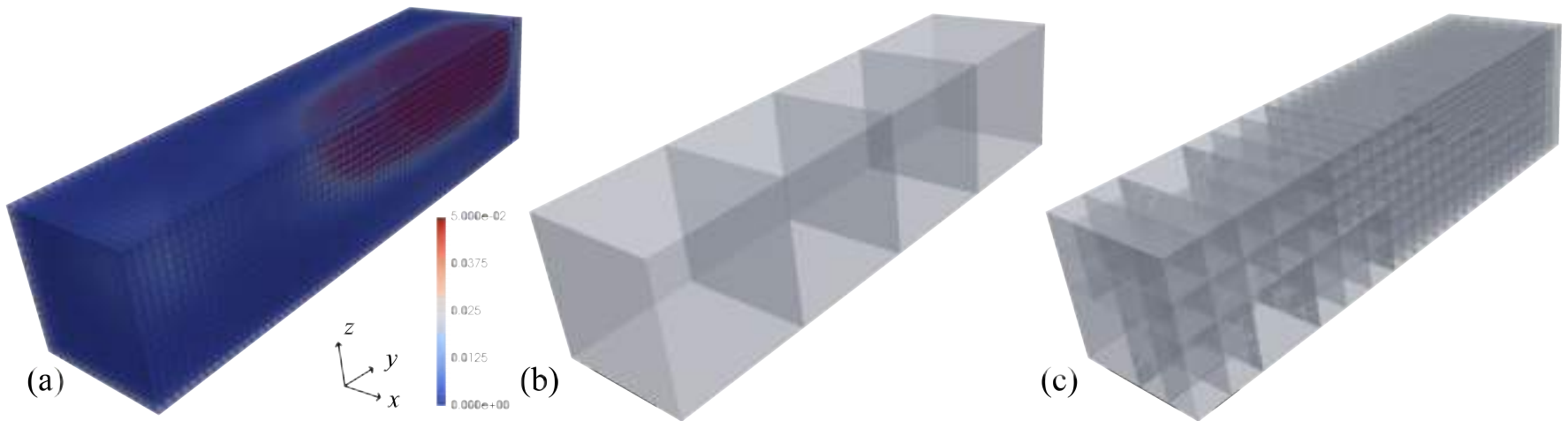

FIGURE 8.2. (a) Converged solution; (b) initial mesh; (c) mesh after twelve SMRs.

Now, consider a goal functional $G \in \mathcal{U}'$ defined in terms of the solution in a region far away from the largest solution gradients. In this circumstance, it is conceivable that the best possible QoI estimate for a fixed computational expenditure, $G(\boldsymbol{u}_h)$, would require a mesh with a refinement pattern very different than one coming from the standard series of SMRs. The following results will clearly verify this conjecture.

8.2. **Experimental design.** Beginning with the four-element mesh depicted in Figure 8.2 (b), we analyzed four different QoIs with goal functionals of the type considered in (7.2) and one pointwise-value QoI (see Section 8.7). Here, is it convenient to express $g \in \mathcal{U}$ as

$$G(\boldsymbol{u}) = \int_\Omega g_1 \cdot u + \int_\Omega \boldsymbol{g}_2 \cdot \boldsymbol{\sigma} + \int_{\Gamma_\mathrm{N}} \hat{g}_3 \cdot \hat{u} + \int_{\Gamma_\mathrm{D}} \hat{g}_4 \cdot \hat{\sigma}\,, \quad \forall \boldsymbol{u} = (u, \boldsymbol{\sigma}, \hat{u}, \hat{\sigma}) \in \mathcal{U}\,. \tag{8.2}$$

In our experiments, $g_1 \in L^2(\Omega)$, $\boldsymbol{g}_2 \in \boldsymbol{L}^2(\Omega)$, $\hat{g}_3 \in L^2(\Gamma_\mathrm{N})$, and $\hat{g}_4 \in H^1_0(\Gamma_\mathrm{D})$ are each piecewise polynomial. Because the ad hoc indicator $\eta_K^{*,\mathrm{a.h.}}$ is not suitable for goal functionals with nonzero boundary contributions (see Remark 7.1), when $G(\boldsymbol{u})$ involved nonzero $\hat{g}_3$ or nonzero $\hat{g}_4$ and $\eta_K^{*,\mathrm{a.h.}}$ was employed, we instead analyzed a sequence of modified goal functionals (see Sections 8.5 and 8.6). In Section 8.7, an extension of this technique is presented for a pointwise-value QoI, $G_{\widetilde{\mathbf{x}}}(\boldsymbol{u}) = u(\widetilde{\mathbf{x}})$.

All of our computations were performed with the finite element software $hp$3D which has complete 3D support for local hierarchical anisotropic $h$ and $p$ refinements with one level of hanging nodes [21, 31] and shape functions for all standard elements conforming in each of the canonical 3D exact sequence energy spaces [41]:

$$H^1(K) \xrightarrow{\mathbf{grad}} \boldsymbol{H}(\mathbf{curl}, K) \xrightarrow{\mathbf{curl}} \boldsymbol{H}(\mathrm{div}, K) \xrightarrow{\mathrm{div}} L^2(K)\,.$$

In our third experiment (Section 8.5), non-homogeneous Neumann boundary conditions were applied to a non-trivial subset of the boundary $\Gamma_\mathrm{N} \subseteq \partial\Omega$. In order to apply this essential boundary condition to the $\hat{\sigma}$-variable, we used projection-based interpolation [22]. For each of the energy spaces above, this is a fully-supported feature of the $hp$3D software.

To implement the practical DPG and DPG* methods, for the ultraweak form of Poisson's boundary value problem, polynomial discretizations of $\mathcal{U}_h$ and $\mathcal{V}_r$ were inferred from previous studies [38, 40, 50, 52]. For the discrete trial space $\mathcal{U}_h$, $L^2(\Omega)$ and $\boldsymbol{L}^2(\Omega)$ were taken from a fixed $p$-order exact sequence. Meanwhile, the interface spaces $H^{1/2}(\mathcal{S})$ and $H^{-1/2}(\mathcal{S})$, were discretized by the $H^1$- and $\boldsymbol{H}(\mathrm{div})$-elementwise traces of the $H^1(\Omega)$- and $\boldsymbol{H}(\mathrm{div}, \Omega)$-conforming shape functions from the same fixed $p$-order exact sequence. Roughly speaking, our assembly of the trace variables $\hat{u}$ and $\hat{\sigma}$ follows established procedure in their respective energy spaces, with the additional extra step of removing the interior bubbles from the final stiffness matrix (see [27, Section 5]). Finally, for the discrete test space $\mathcal{V}_r$, discretizations of $H^1(\mathcal{T})$ and $\boldsymbol{H}(\mathrm{div}, \mathcal{T})$ were taken from a fixed, non-conforming and enriched $(p + \mathrm{d}p)$-order exact sequence.

For a careful account of DPG assembly algorithms, see [66] and [52, Section 4]. Note that we always assembled the full normal equation [52, Section 4.1]. In the DPG* setting, the stiffness matrix is identical and so, on any one mesh, only a single assembly and factorization of the DPG stiffness matrix was ever performed to solve for both $\boldsymbol{u}_{h,r}$ and $\boldsymbol{v}_{h,r}$.

*Parameters.* Before finally presenting the results of our experiments, we now list the outstanding parameters in our algorithms and our choices for them in all experiments:

- Discretizations of the trial space $\mathcal{U}_h$ and the test space $\mathcal{V}_r$ came from an exact sequence of polynomial order $p = 2$ and $p + \mathrm{d}p = 3$, respectively [39]. Note that the polynomial order of the manufactured solution (8.1) is too high for it to be fully recovered with this trial space discretization.
- In the implicit refinement indicator $\eta_K^{*,\mathrm{impl}}$, the local problems (7.9) were solved on individual elements $K$ from the same mesh $\mathcal{T}$ as the global DPG and DPG* problems. Here, the enriched trial space $\mathcal{U}_H$ and further enriched test space $\mathcal{V}_R$ were constructed as previously described for $\mathcal{U}_h$ and $\mathcal{V}_r$, but with basis functions taken from exact sequences of polynomial order $P = p + 1 = 3$ and $P + \mathrm{d}p = 4$, respectively.
- The adjoint graph norm (3.23) was used, $\|\cdot\|_{\mathcal{V}} = \|\cdot\|_{\mathcal{L}^*,\alpha}$, with $\alpha = 1$.
- The refinement factor of $\theta = 0.5$ was set for both the SMR and GMR marking strategies (see Section 5.1).

8.3. **Temperature in a subdomain.** In this subsection, we consider the goal functional given in (8.2) where

$$g_1(x,y,z) = \begin{cases} 1\,, & x \leq 1\,, \\ 0\,, & \text{otherwise}, \end{cases} \qquad \boldsymbol{g}_2 = \boldsymbol{0}\,, \qquad \hat{g}_3 = 0\,, \qquad \text{and} \qquad \hat{g}_4 = 0\,. \tag{8.3}$$

Physically, this corresponds to a QoI which is the average value of the temperature $u$ in the subdomain $0 \leq x \leq 1$. In these experiments, we used homogeneous Dirichlet boundary conditions, $u|_{\partial\Omega} = 0$.

8.3.1. *Results.* Define the relative error in the QoI to be $|G(\boldsymbol{u}^\star - \boldsymbol{u}_{h,r})|/|G(\boldsymbol{u}^\star)|$. In Figure 8.3 (a), we present the relative error vs. the degrees of freedom in each successive solution as the mesh was refined using each AMR strategy. From now on, explicit, implicit, and ad hoc GMR refers to $\widetilde{\eta}_K := \eta_K \cdot \eta_K^*$ with $\eta_K^* = \eta_K^{*,\mathrm{expl}}$ defined in (7.13), $\eta_K^* = \eta_K^{*,\mathrm{impl}}$ defined in (7.10), and $\eta_K^* = \eta_K^{*,\mathrm{a.h.}}$ defined in (7.11), respectively. It is immediately evident that each GMR step was far more efficient at reducing the relative error in the QoI than each SMR step. Moreover, taking into account their entire sequence of refinements, each GMR strategy performed nearly equally, until nearly the final refinement.

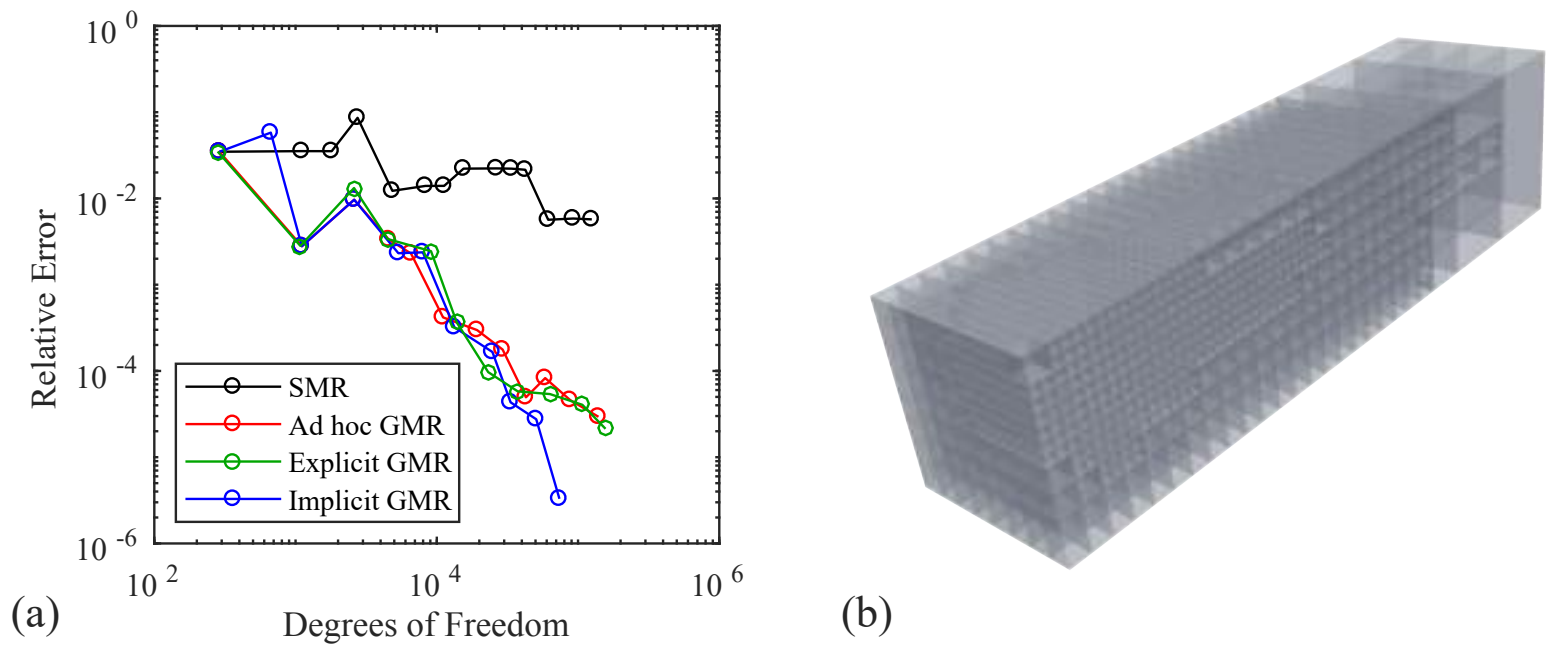

FIGURE 8.3. (a) The error in the QoI for the first example: average temperature $u$; (b) the final adaptively refined mesh using the ad hoc GMR approach.

In Figure 8.3 (b), we see the final refined mesh after twelve adaptive mesh refinements with the GMR marking strategy and the ad hoc refinement indicator $\eta_K^{*,\mathrm{a.h.}}$. However, because there are two other classes of DPG* refinement indicators which preformed well for this problem and QoI, Figure 8.4 is provided to compare all three corresponding final solution and meshes. Here, it is visibly evident that the final meshes are extremely similar, but significantly different from the SMR mesh in Figure 8.2 (c). A strong visual similarity in the final

GMR meshes was also exhibited in each of our studies. Therefore, from now on, we will only provide one representative GMR mesh for illustration.

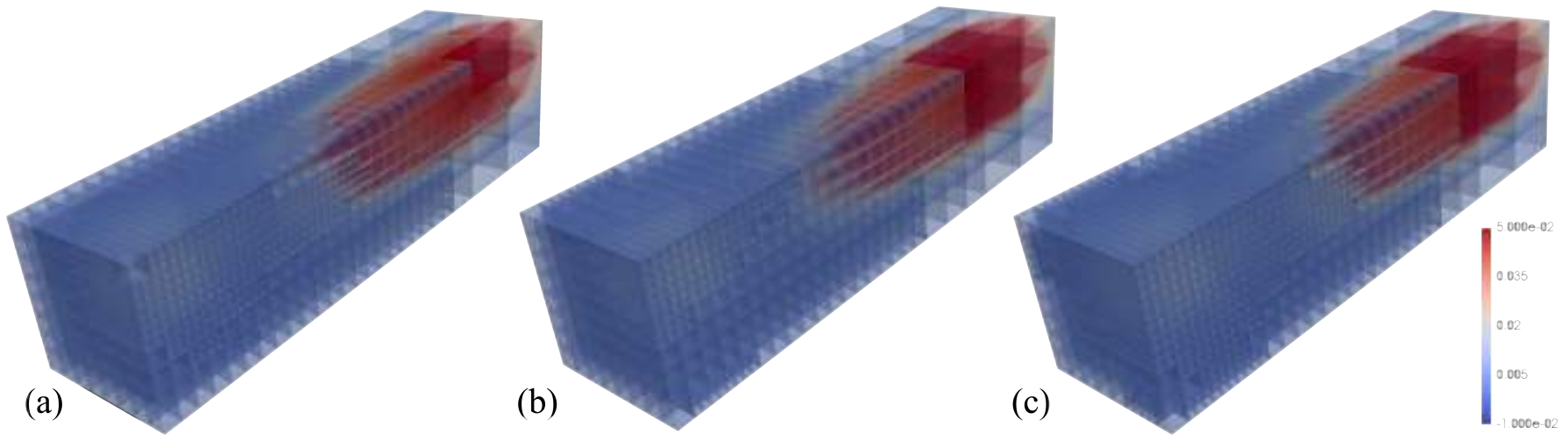

FIGURE 8.4. The solution after twelve refinement steps using the: (a) ad hoc GMR; (b) explicit GMR; (c) implicit GMR.

Finally, we provide Figure 8.5 for a visual depiction of the local temperature error in the region $0 \leq x \leq 1$. In some contexts, a goal functional of the form (8.3) is chosen to drive adaptivity with the intention of significantly reducing the error in a particular solution variable—in this case, it is the temperature $u$—in a *region of interest*. Although this can also be done more accurately by considering a nonlinear goal functional [60], simply using GMR with a closely related linear QoI often provides a sufficient improvement. With this understanding, Figure 8.5 clearly demonstrates that, the total error in the temperature variable $u$ in the region of interest, is far lower as a result of the GMRs as opposed to the conventional SMR for a similar number of degrees of freedom. In Figure 8.5 (b), we have only visualized the error from the ad hoc approach, however, the results from the other two approaches were nearly indistinguishable in comparison.

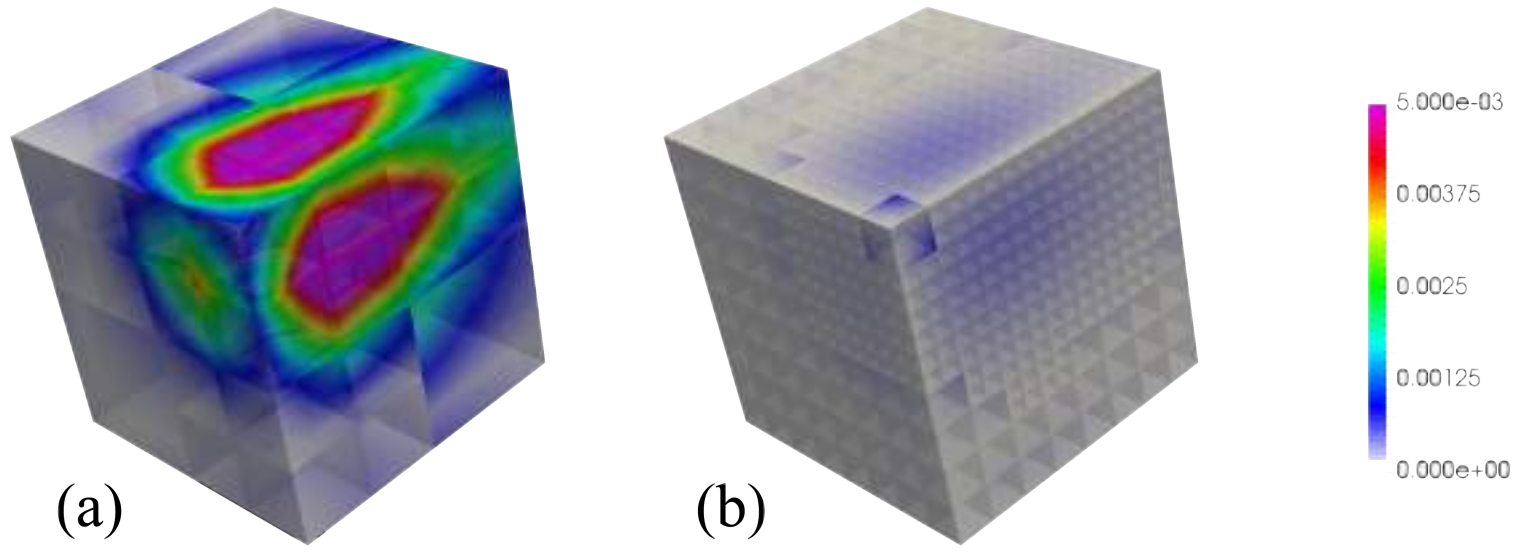

FIGURE 8.5. The error in the solution component $u$ in $0 \leq x \leq 1$ at the final mesh using: (a) conventional DPG SMR; (b) ad hoc DPG GMR.

8.3.2. *The influence function.* It is worth remarking, once again, that the ad hoc refinement indicator $\eta^{*,\mathrm{a.h.}}$ will not converge to zero asymptotically. Indeed, the success of these refinement indicators in a GMR strategy highly depends upon the quality of resolution of the optimal test functions and the difference between the optimal test norm and the norm chosen for computation.

Recall that $g_1 = 1$ if $x \leq 1$ and $g_1 = 0$ otherwise. Meanwhile, $\boldsymbol{g}_2 = \mathbf{0}$ everywhere. Our driving assumption in Section 7.4 is that $(\omega_1, \boldsymbol{\omega}_2) \approx (g_1, \boldsymbol{g}_2)$ when using the adjoint graph norm with small enough parameter $\alpha$. The affirmative results in Figure 8.3 (a) indicate that this assumption was reasonably valid for this problem, even though we only used a moderate $\alpha = 1$. Figure 8.6 is also provided to further justify this conclusion. Here, the $\omega_1$-component of the approximated influence function $\boldsymbol{\omega}_{h,r}$ is visualized on both the initial and final mesh. Clearly, even when computing on the initial mesh, which had only four elements, the difference between the computed $\omega_1$ and its idealized value, $\omega_1 \overset{\mathrm{opt}}{=} g_1$, was only marginal.

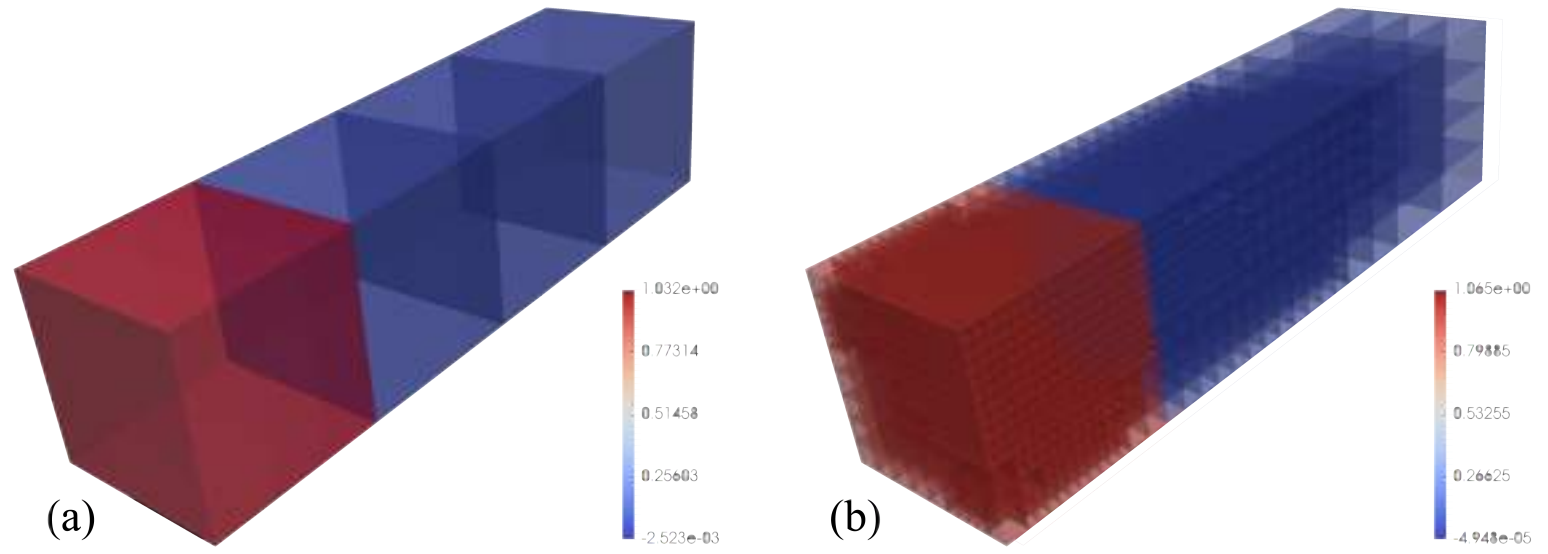

FIGURE 8.6. The influence function $\omega_1$ for: (a) the initial mesh; (b) final adaptively refined mesh using the ad hoc GMR approach.

8.4. **Flux in a subdomain.** In this subsection, we consider the goal functional given in (8.2) where

$$g_1 = 0\,, \qquad \boldsymbol{g}_2(x,y,z) = \begin{cases} (1,0,0)^\mathsf{T}\,, & x \leq 1\,, \\ \mathbf{0}\,, & \text{otherwise}, \end{cases} \qquad \hat{g}_3 = 0\,, \qquad \text{and} \qquad \hat{g}_4 = 0\,. \tag{8.4}$$

Physically, this corresponds to a QoI which is the average value of the $x$-component of the flux, $\sigma_x$, in the subdomain $0 \leq x \leq 1$. In these experiments, homogeneous Dirichlet boundary conditions, $u|_{\partial\Omega} = 0$, were used.

Define the relative error in the QoI to be $|G(\boldsymbol{u}^\star - \boldsymbol{u}_{h,r})|/|G(\boldsymbol{u}^\star)|$. As with the previous experiment, we present the relative error in this QoI with each of the AMR strategies. Again, by inspecting Figure 8.7 (a), it is clear that each of the GMR strategies are far more efficient than conventional DPG SMR.

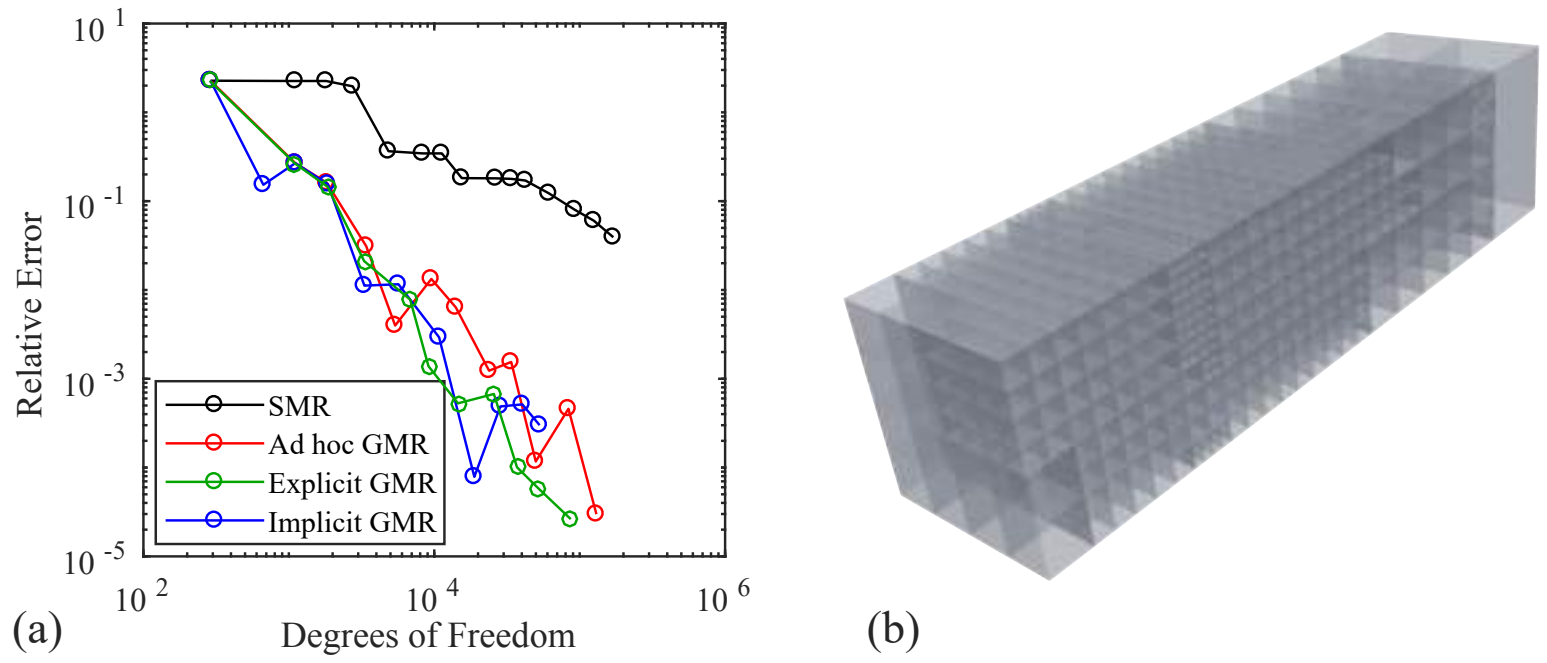

FIGURE 8.7. (a) The error in the QoI for the second example: average flux $\sigma_x$; (b) the final adaptively refined mesh using the explicit GMR approach.

For a visual comparison of the error in $\sigma_x$ in the region of interest $0 \leq x \leq 1$, we provide Figure 8.8. Even though they both come from meshes inducing a similar number of degrees of freedom, notice that the local error from the GMR strategy is at least two orders of magnitude smaller than with the local error from the SMR strategy.

8.5. **Temperature on the boundary.** In this subsection, consider the goal functional $g$ given in (8.2) where

$$g_1 = 0\,, \qquad \boldsymbol{g}_2 = \mathbf{0}\,, \qquad \hat{g}_3(x,y,z) = \begin{cases} 1\,, & x = 0\,, \\ 0\,, & \text{otherwise}, \end{cases} \qquad \text{and} \qquad \hat{g}_4 = 0\,. \tag{8.5}$$

Physically, this corresponds to a QoI which is the average value of the temperature $u$ on the subboundary $x = 0$. Here, homogeneous Dirichlet and homogeneous Neumann boundary conditions were used on disjoint regions of the boundary, $u|_{\Gamma_\mathrm{D}} = 0$ and $\frac{\partial u}{\partial n}|_{\Gamma_\mathrm{N}} = 0$, where $\Gamma_\mathrm{D} = \{(x,y,z) \in \partial\Omega \mid x > 0\}$ and $\overline{\Gamma_\mathrm{D} \cup \Gamma_\mathrm{N}} = \partial\Omega$.

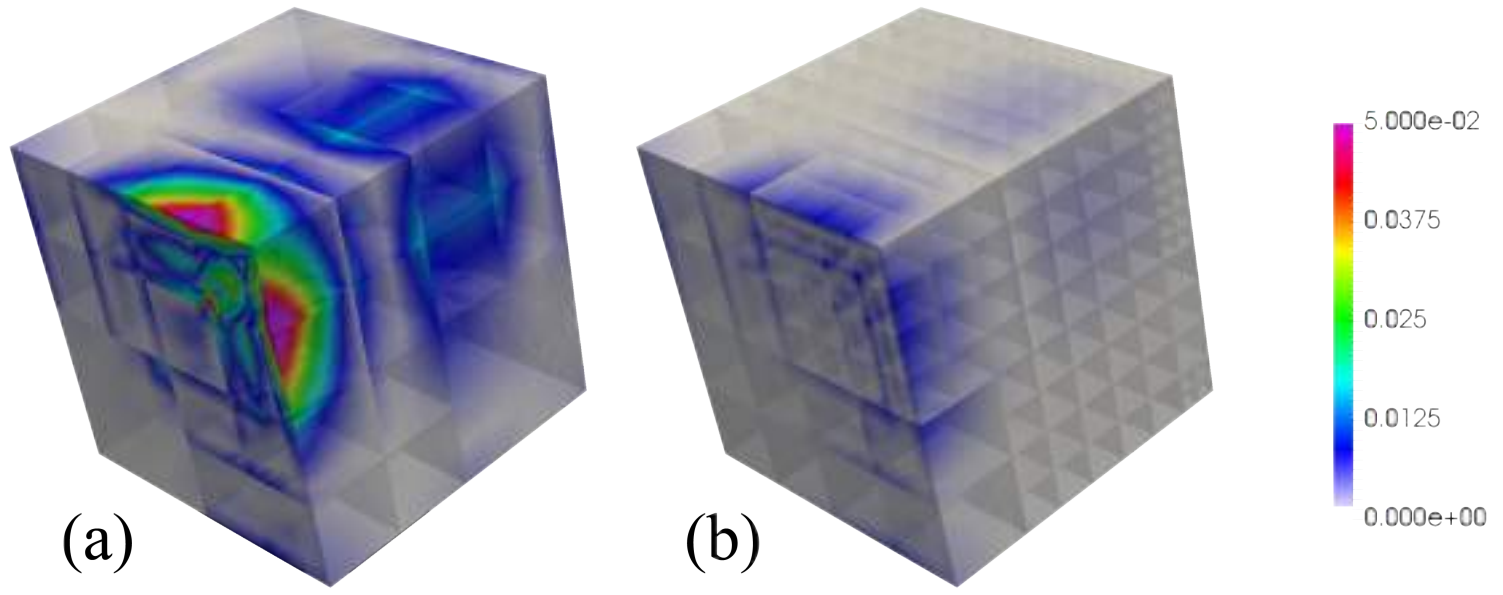

FIGURE 8.8. The final error in the solution component $\sigma_x$ in the subdomain $0 \leq x \leq 1$ using: (a) SMR; (b) explicit GMR.

8.5.1. *An alternative functional for ad hoc refinement indicators.* Given that the ad hoc refinement indicators $\eta_K^{*,\mathrm{a.h.}}$ were not derived for goal functionals involving nonzero $\hat{g}_3$ or $\hat{g}_4$, we actually used a different goal functional (8.6) for this indicator *only*.

Let $\Gamma = \{(x, y, z) \in \partial\Omega \mid x = 0\}$ and define the set of all elements neighboring $\Gamma$ as $\mathcal{T}_\Gamma = \{K \in \mathcal{T} \mid \overline{K} \cap \Gamma \neq \varnothing\}$. Lastly, define the subdomain occupied by this set of elements as $\Omega_\Gamma = \bigcup_{K \in \mathcal{T}_\Gamma} \overline{K}$. Now, instead of (8.5), consider the modified goal functional

$$g_1(\mathbf{x}) = \begin{cases} h_{K,x}^{-1}, & \mathbf{x} \in \Omega_\Gamma, \\ 0, & \text{otherwise,} \end{cases} \qquad \boldsymbol{g}_2 = \mathbf{0}, \qquad \hat{g}_3 = 0, \qquad \text{and} \qquad \hat{g}_4 = 0, \tag{8.6}$$

where, $h_{K,x}$ is the length, in the $x$-dimension, of the element $K \in \mathcal{T}_\Gamma$ enclosing the point $\mathbf{x} \in K$. Notice that because $g_1$ operates on the $u$-component of the solution in (8.2), as the mesh becomes finer near $\Gamma$, this functional will also limit to a characterization of the average temperature on the boundary. The primary novelty of (8.6) is that the definition is mesh dependent. This is demonstrated in Figure 8.9, where $\Omega_\Gamma$ is highlighted in red on different adaptively refined meshes.

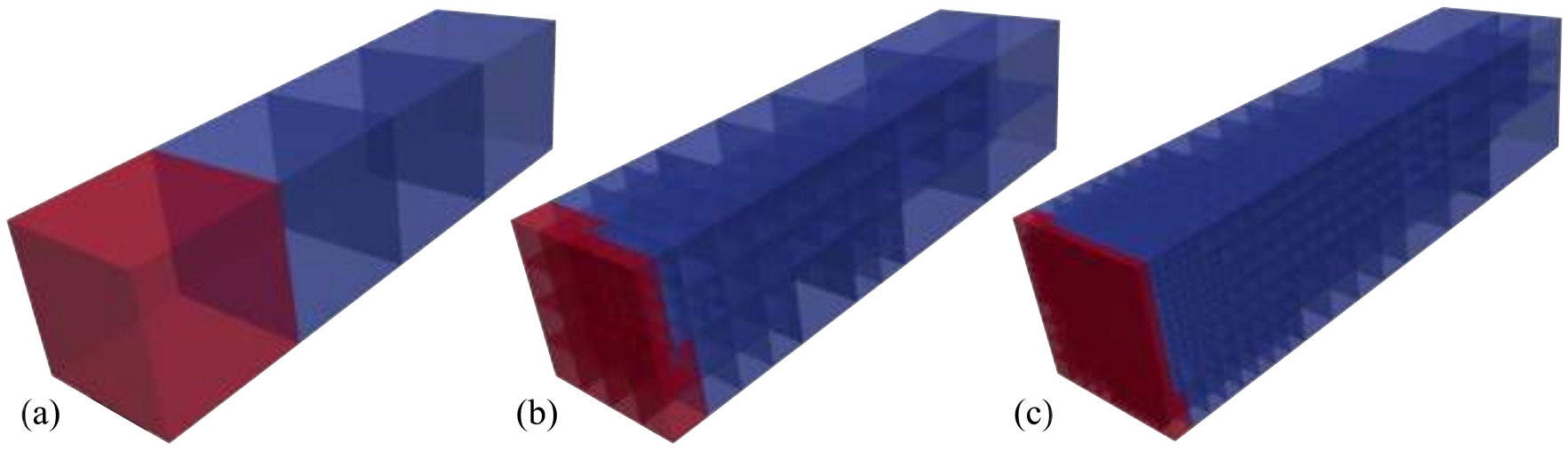

FIGURE 8.9. The region of interest $\Omega_\Gamma$, marked in red, for: (a) the initial mesh; (b) the mesh after six GMR steps with the ad hoc approach; (c) the mesh after twelve ad hoc GMR steps.

8.5.2. *Results.* Notice that $u^\star|_{\partial\Omega} = 0$ from (8.1). Therefore, from the definition of $G$ given in (8.5), $G(\boldsymbol{u}^\star) = \int_{\Gamma_{\mathrm{N}}} u^\star = \int_\Gamma u^\star = 0$. Define the relative error in the QoI to be $\left|\int_\Gamma \hat{u}_{h,r}\right| / \left|\int_\Gamma \hat{u}_{h_{\mathrm{init}},r}\right|$, where $\hat{u}_{h_{\mathrm{init}},r}$ is the third component of approximate solution $\boldsymbol{u}_{h,r}$ computed on the initial mesh, which all approaches having in common (see Figure 8.2 (a)). From only a cursory inspection of Figure 8.10, it is again evident that each of the GMR strategies were far more efficient than conventional SMR.

In Figure 8.11, the visual comparison given of the error in the temperature variable on the region of interest, $x = 0$, demonstrates a substantial improvement over the conventional SMR strategy, even with the ad hoc GMR approach employing the modified goal functional (8.6).

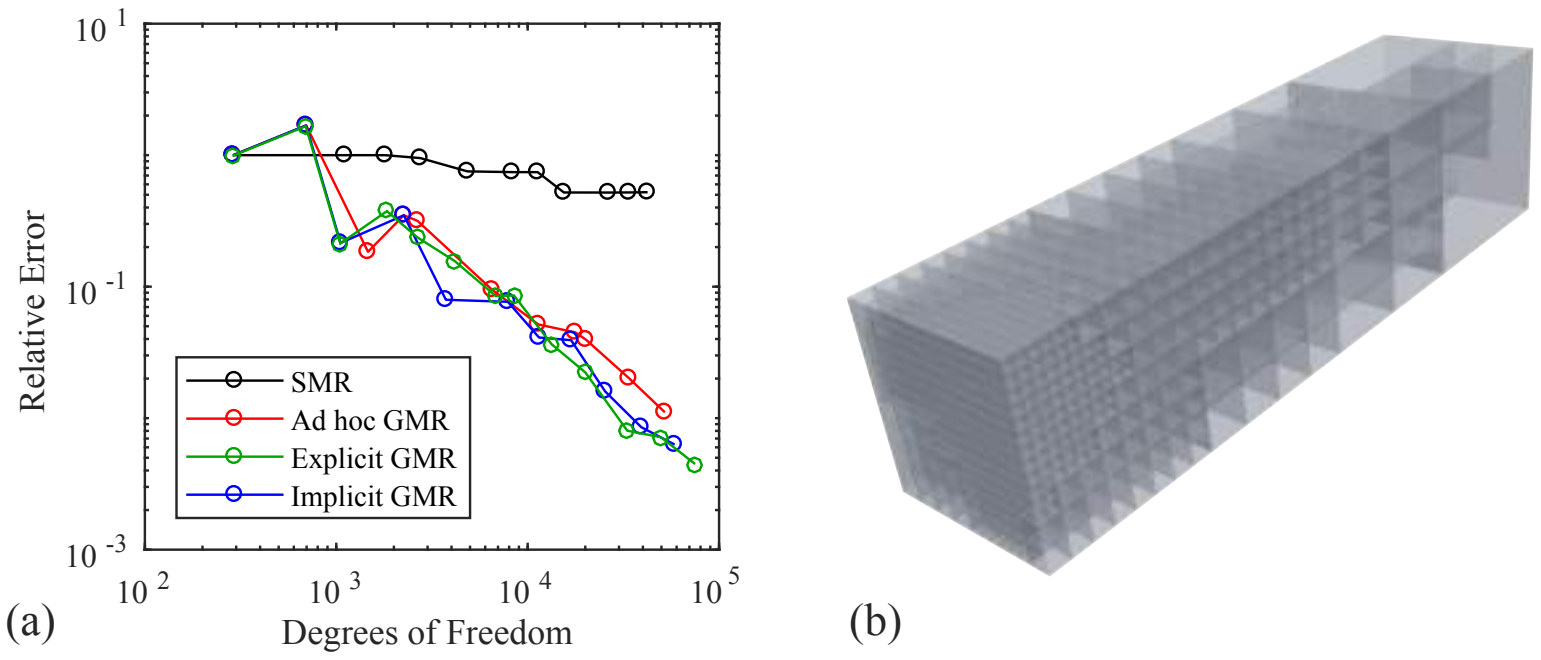

FIGURE 8.10. (a) The error in the QoI for the third example: average temperature $\hat{u}$; (b) the final adaptively refined mesh using the ad hoc GMR approach.

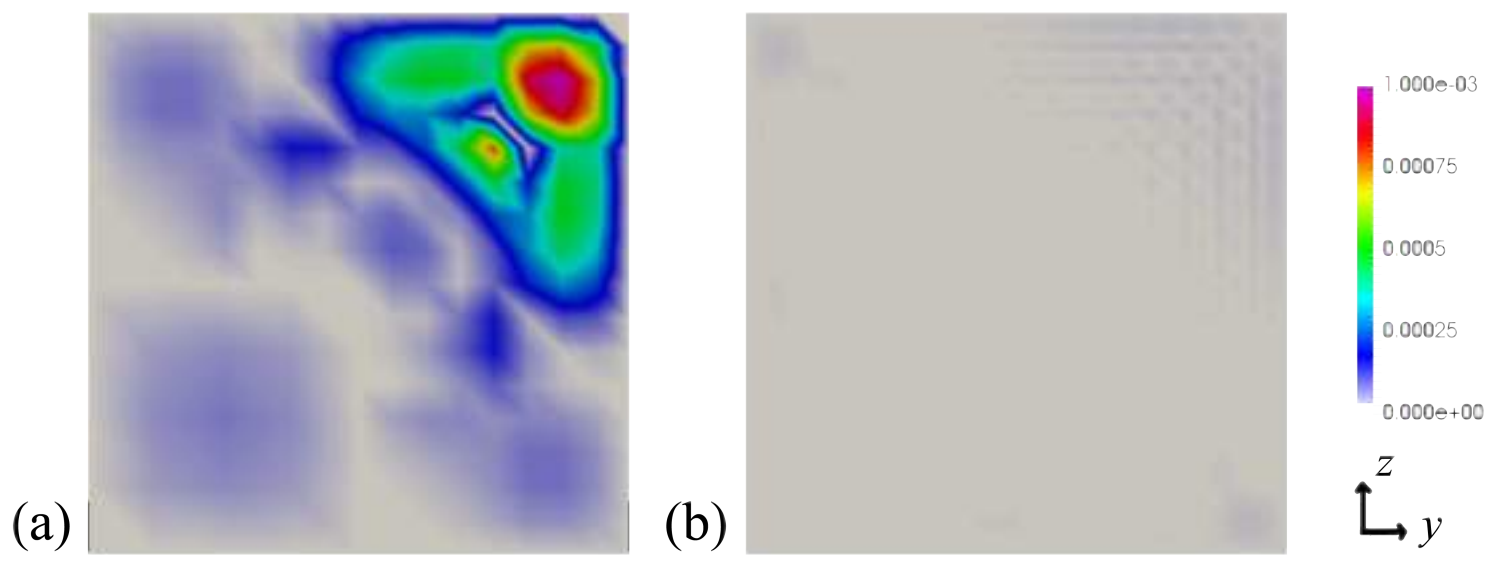

FIGURE 8.11. The final error in the solution $\hat{u}$ on $x = 0$ from: (a) conventional SMR; (b) the ad hoc GMR approach.

8.6. **Flux on the boundary.** Consider the goal functional given in (8.2) where

$$g_1 = 0\,, \qquad \boldsymbol{g}_2 = \boldsymbol{0}\,, \qquad \hat{g}_3 = 0\,, \qquad \text{and} \qquad \hat{g}_4(x,y,z) = \begin{cases} 1\,, & x = 0\,, \\ 0\,, & \text{otherwise.} \end{cases} \tag{8.7}$$

Physically, this corresponds to the QoI being the average value of the flux $\boldsymbol{\sigma}\cdot\boldsymbol{n}$ through the subboundary $x = 0$. In our experiments, we used homogeneous Dirichlet boundary conditions everywhere, $u|_{\partial\Omega} = 0$.

8.6.1. *Energy space considerations.* In this experiment, $\Gamma_{\mathrm{D}} = \partial\Omega$ and so $\big(H_{\mathrm{N}}^{-1/2}(\partial\Omega)\big)' = H^{1/2}(\partial\Omega)$ and $H^1_0(\Gamma_D) = H^1(\partial\Omega)$. Notice that $\hat{g}_4 \in L^2(\partial\Omega)$ but $\hat{g}_4 \notin H^1(\partial\Omega)$. Therefore, the assumptions of theorem 7.5 are not met. In fact, there is no prerequisite reason why $g$, as defined by (8.7), should even be a bounded linear functional on $\mathcal{U}$. Indeed, because it has a nontrivial jump discontinuity, $\hat{g}_4 \notin H^{1/2}(\partial\Omega) \supseteq H^1(\partial\Omega)$ and so $G \notin \mathcal{U}'$.

Unfortunately, violating the energy setting is not simply a mathematical concern. Indeed, we found spuriously concentrated refinements near the discontinuity in $\hat{g}_4$, when using (8.7). Therefore, instead of (8.7), we chose to mollify the physically ideal (but discontinuous) $\hat{g}_4$ to an extent that it obeys the proper energy setting. For our explicit and implicit GMR experiments, specifically, we used the ramp function depicted in Figure 8.12, in place of the function $\hat{g}_4$ defined in (8.7).

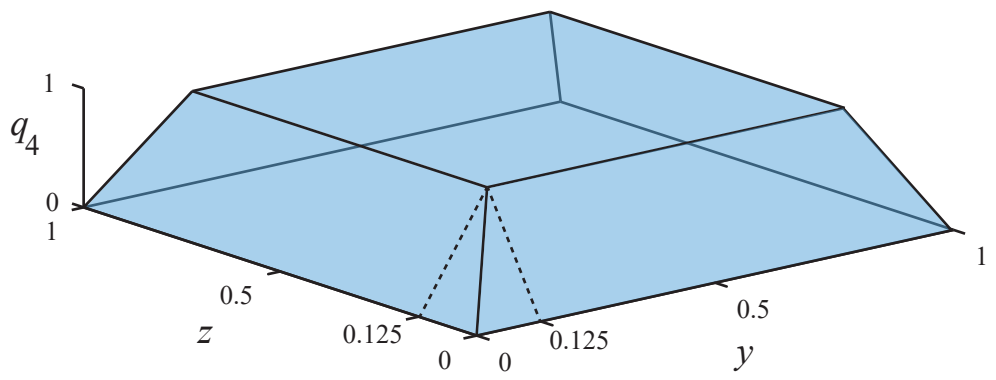

FIGURE 8.12. Illustration of the function $\hat{g}_4 \in H^1(\partial\Omega)$ (on the face $x = 0$) used for explicit and implicit GMR.

8.6.2. *Another alternative functional for ad hoc refinement indicators.* Analogous to the previous experiment, here we considered an alternative to the goal functional (8.7) for the ad hoc DPG* refinement indicator. Namely, we used

$$g_1 = 0 \qquad \boldsymbol{g}_2(\mathbf{x}) = \begin{cases} (h_{K,x}^{-1}, 0, 0)^{\mathsf{T}}, & \mathbf{x} \in \Omega_\Gamma, \\ \mathbf{0}, & \text{otherwise}, \end{cases}, \qquad \hat{g}_3 = 0, \qquad \text{and} \qquad \hat{g}_4 = 0, \tag{8.8}$$

with the same definitions for $\Omega_\Gamma$ and $h_{K,x}$ as in (8.6). Fortunately, with this definition, $\boldsymbol{g}_2 \in L^2(\Omega)$ and so $G \in \mathcal{U}'$ for all meshes and the energy space issues for (8.7) are again avoided.

8.6.3. *Results.* Define the relative error in the QoI to be $\left|\int_\Gamma(\boldsymbol{\sigma}^\star \cdot \boldsymbol{n} - \hat{\sigma}_{h,r})\right| / \left|\int_\Gamma \boldsymbol{\sigma}^\star \cdot \boldsymbol{n}\right|$. An inspection of Figure 8.13 clearly illustrates that the GMR strategies were far more efficient than the conventional SMR strategy for this QoI. In fact, with the conventional strategy, the error in this QoI did not even decrease until the eighth mesh refinement was performed!

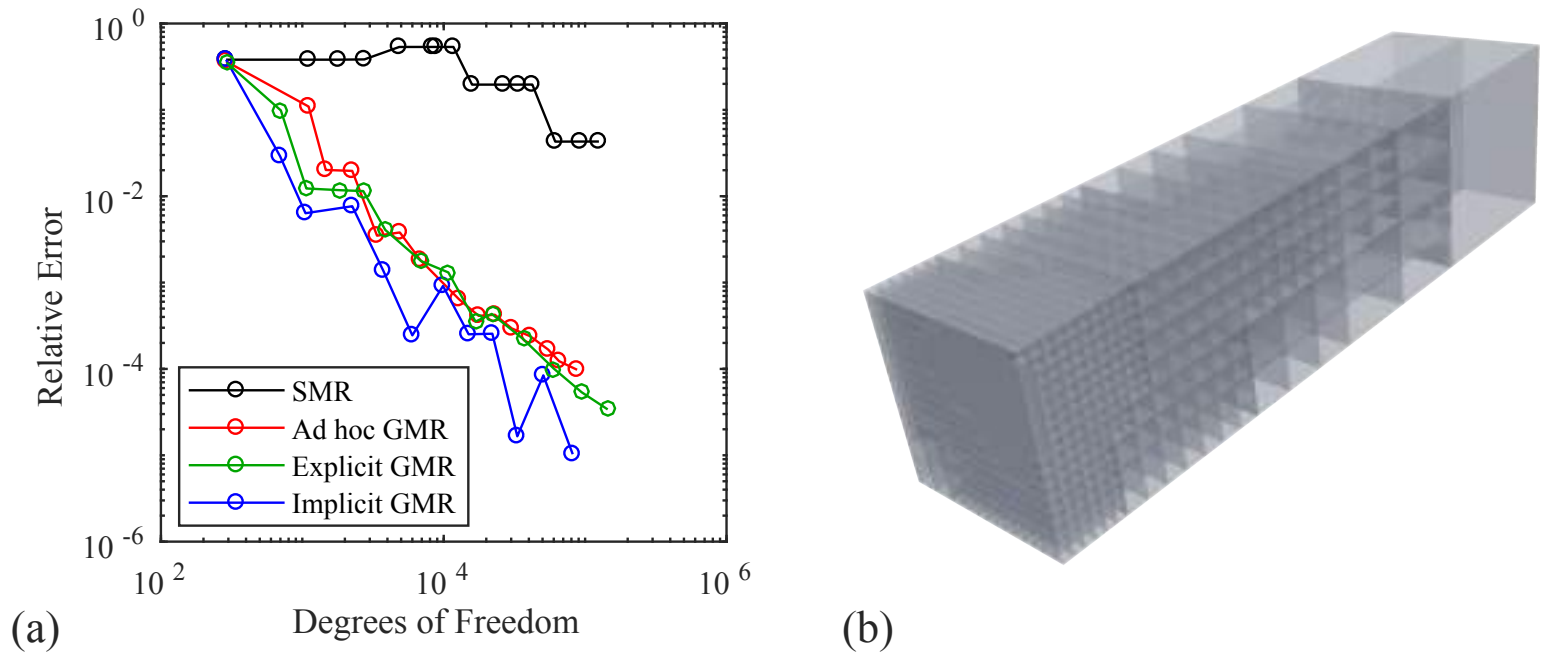

FIGURE 8.13. (a) The error in the QoI for the fourth example: average normal flux $\hat{\sigma}$; (b) the final adaptively refined mesh using the implicit GMR approach.

As for the visual comparison of the error in the flux on the region of interest, Figure 8.14 again demonstrates a significant improvement over conventional DPG SMR.

8.7. **Temperature at a point.** In this final experiment, we considered the goal functional $G_{\widetilde{\mathbf{x}}}(\boldsymbol{u}) = u(\widetilde{\mathbf{x}})$, where $\widetilde{\mathbf{x}} \in \Omega$ is a specified point in the domain. Markedly, this QoI does not fall into the theory of this article because $G_{\widetilde{\mathbf{x}}} \notin \mathcal{U}'$. To overcome this issue, we a mesh-dependent goal functional like (8.6) and (8.8). Define the set of all elements containing the point $\widetilde{\mathbf{x}}$ in their closure as $\mathcal{T}_{\widetilde{\mathbf{x}}} = \{K \in \mathcal{T} \mid \widetilde{\mathbf{x}} \in \overline{K}\}$ and define the subdomain occupied by this set of elements $\Omega_{\widetilde{\mathbf{x}}} = \bigcup_{K \in \mathcal{T}_{\widetilde{\mathbf{x}}}} \overline{K}$. Now, redefine the goal functional as

$$g_1(\mathbf{x}) = \begin{cases} \mathrm{vol}(\Omega_{\widetilde{\mathbf{x}}})^{-1}, & \mathbf{x} \in \Omega_{\widetilde{\mathbf{x}}}, \\ 0, & \text{otherwise}, \end{cases} \qquad \boldsymbol{g}_2 = \mathbf{0}, \qquad \hat{g}_3 = 0, \qquad \text{and} \qquad \hat{g}_4 = 0.$$

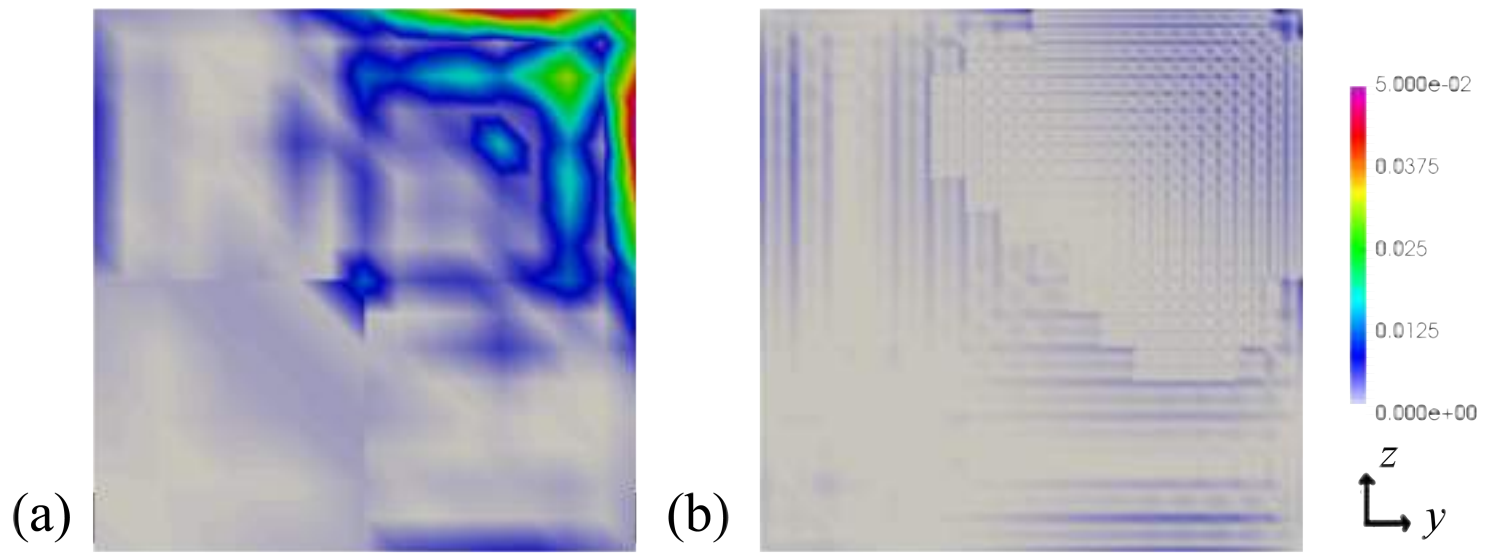

FIGURE 8.14. The final error in the solution $\hat{\sigma}$ at $x = 0$ using: (a) SMR; (b) implicit GMR.

All results presented in this section are from our experiments with the mesh-dependent definition above and $\widetilde{\mathbf{x}} = (0.3, 0.3, 0.3)^{\mathsf{T}}$. The evolution of this functional is visually depicted in Figure 8.15, where the region of interest $\Omega_{\widetilde{\mathbf{x}}}$ is highlighted in red on a selection of meshes.

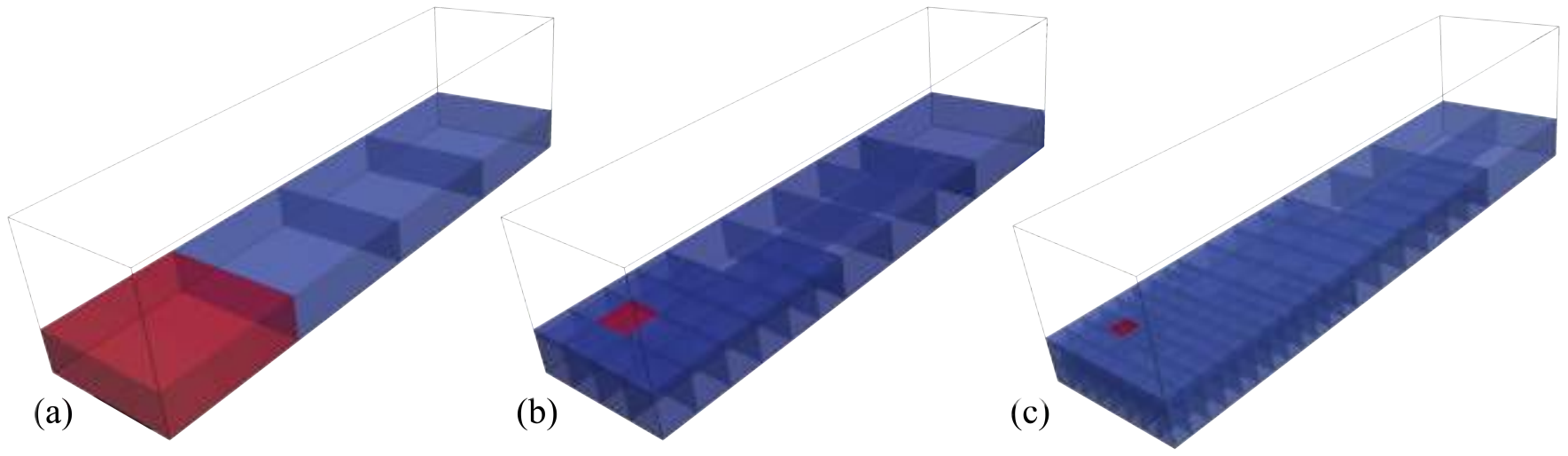

FIGURE 8.15. The region of interest $\Omega_{\widetilde{\mathbf{x}}}$ enclosing the point $\widetilde{\mathbf{x}} = (0.3, 0.3, 0.3)^{\mathsf{T}}$ for: (a) the initial mesh; (b) the mesh after six refinements with the implicit approach; (c) the mesh after twelve refinements. For the sake of illustration, only $y \leq 0.3$ is shown.

Inspect Figure 8.16. Here, the relative error is defined to be $|u^{\star}(\widetilde{\mathbf{x}}) - u_{h,r}(\widetilde{\mathbf{x}})|/|u^{\star}(\widetilde{\mathbf{x}})|$. As in every previous experiment, each GMR approach vastly outperformed conventional SMR. However, in this experiment, the convergence behavior of each GMR approach was quite erratic.

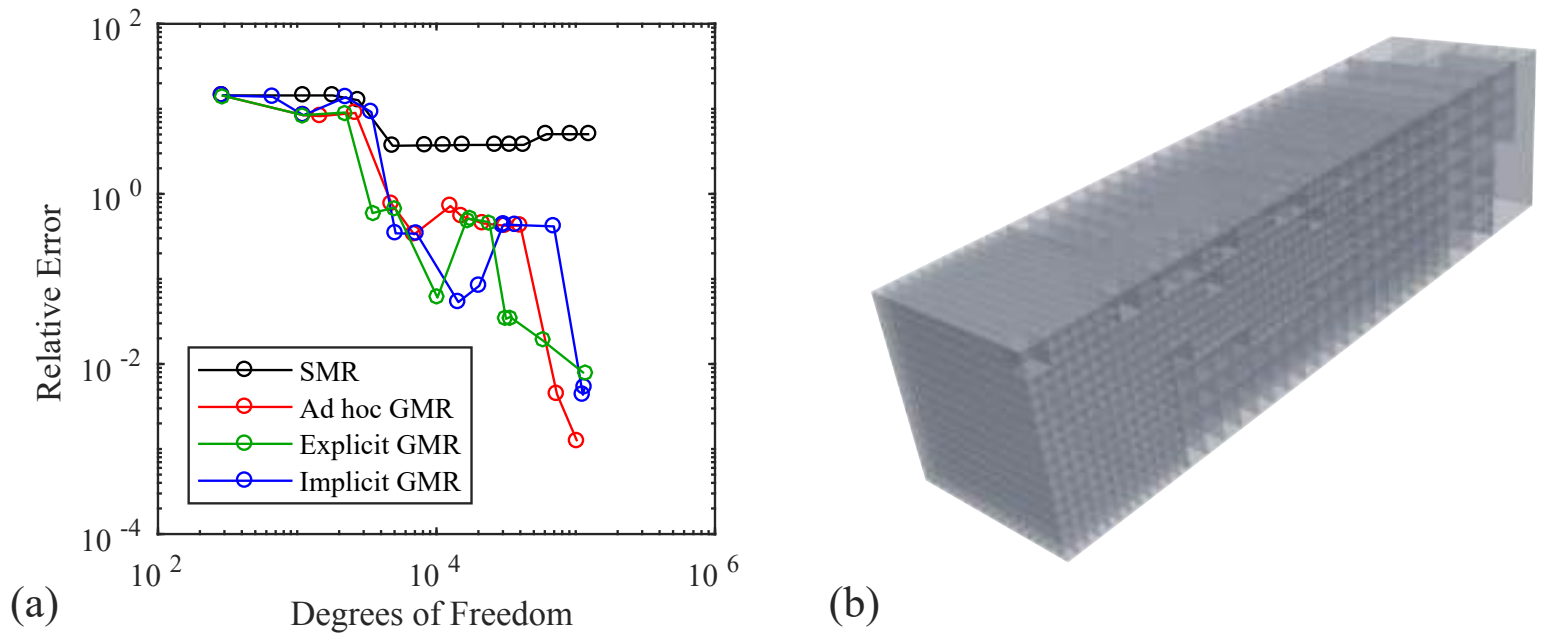

FIGURE 8.16. (a) The error in the QoI for the fifth example: temperature $u$ at a fixed point; (b) the final adaptively refined mesh using the explicit GMR approach.

## APPENDIX A. A PYTHAGOREAN THEOREM FOR BOUNDED PROJECTIONS

The purpose of this appendix is to prove Theorem 6.3. For convenience, we restate it again separately, as Theorem A.1.

**Theorem A.1** (Pythagoras). *Let $\mathcal{W}$ be a Hilbert space and $\mathcal{W}_0 \subseteq \mathcal{W}$ be a nontrivial closed subspace. Let $\mathcal{P} : \mathcal{W} \to \mathcal{W}_0$ be the orthogonal projection onto $\mathcal{W}_0$ and let $\Pi : \mathcal{W} \to \mathcal{W}_0$ be any other bounded projection onto $\mathcal{W}_0 = \Pi(\mathcal{W})$. Then*

$$\|\Pi - \mathcal{P}\|^2 + 1 = \|\Pi\|^2 \,.$$

*Proof.* Throughout this proof, the subscript-$\mathcal{W}$ in norms and inner products is suppressed.

Note that, if $\Pi = \mathcal{P}$, we are done. Therefore, assume that $\Pi \neq \mathcal{P}$ and so $\big\{w \in \mathcal{W} \,\big|\, \|(\mathcal{P} - \Pi)w\| > 0\big\} \neq \varnothing$. Define $\mathcal{W}_1 = \mathcal{W}_0^\perp$. By Lemma 3.1, observe that

$$\sup_{w\in\mathcal{W}} \frac{(\widetilde{w}, \Pi w)^2}{\|w\|^2} = \sup_{w_0\in\mathcal{W}_0} \frac{(\widetilde{w}, \Pi w_0)^2}{\|w_0\|^2} + \sup_{w_1\in\mathcal{W}_1} \frac{(\widetilde{w}, \Pi w_1)^2}{\|w_1\|^2}\,, \quad \forall\, \widetilde{w} \in \mathcal{W}\,. \tag{A.1}$$

Note that $\Pi w_0 = w_0$ and $\Pi w_1 = \Pi(1 - \mathcal{P})w_1 = (\Pi - \mathcal{P})w_1$ for any $w_0 \in \mathcal{W}_0$ and $w_1 \in \mathcal{W}_1$, respectively. Therefore,

$$\begin{aligned}
\|\Pi\|^2 = \sup_{w\in\mathcal{W}} \frac{\|\Pi w\|^2}{\|w\|^2} = \sup_{\widetilde{w},w\in\mathcal{W}} \frac{(\widetilde{w}, \Pi w)^2}{\|\widetilde{w}\|^2\|w\|^2} &= \sup_{\widetilde{w}\in\mathcal{W}} \left( \sup_{w_0\in\mathcal{W}_0} \frac{(\widetilde{w}, \Pi w_0)^2}{\|\widetilde{w}\|^2\|w_0\|^2} + \sup_{w_1\in\mathcal{W}_1} \frac{(\widetilde{w}, \Pi w_1)^2}{\|\widetilde{w}\|^2\|w_1\|^2} \right)\\
&= \sup_{\widetilde{w}\in\mathcal{W}} \left( \sup_{w_0\in\mathcal{W}_0} \frac{(\widetilde{w}, w_0)^2}{\|\widetilde{w}\|^2\|w_0\|^2} + \sup_{w_1\in\mathcal{W}_1} \frac{(\widetilde{w}, (\Pi - \mathcal{P})w_1)^2}{\|\widetilde{w}\|^2\|w_1\|^2} \right),
\end{aligned} \tag{A.2}$$

where the third equality follows from (A.1). Clearly,

$$\|\Pi\|^2 \leq \sup_{\widetilde{w},w\in\mathcal{W}} \frac{(\widetilde{w}, w)^2}{\|\widetilde{w}\|^2\|w\|^2} + \sup_{\widetilde{w},w\in\mathcal{W}} \frac{(\widetilde{w}, (\Pi - \mathcal{P})w)^2}{\|\widetilde{w}\|^2\|w\|^2} = 1 + \|\Pi - \mathcal{P}\|^2\,.$$

Now, define $w_\Pi = \arg\max_{\|w\|=1} \|(\Pi - \mathcal{P})w\| \neq 0$ and then $w_{\Pi,0} = (\Pi - \mathcal{P})w_\Pi \neq 0$ . Consider $\widetilde{w} = w_{\Pi,0}$ in (A.2) and observe that

$$\|\Pi\|^2 \geq \sup_{w_0\in\mathcal{W}_0} \frac{(w_{\Pi,0}, w_0)^2}{\|w_{\Pi,0}\|^2\|w_0\|^2} + \sup_{w_1\in\mathcal{W}_1} \frac{(w_{\Pi,0}, (\Pi - \mathcal{P})w_1)^2}{\|w_{\Pi,0}\|^2\|w_1\|^2}.$$

Note that $w_{\Pi,0} \in \mathcal{W}_0$, $w_\Pi \in \mathcal{W}_1$, $\|w_{\Pi,0}\| = \|\Pi - \mathcal{P}\|$, and $\|w_\Pi\| = 1$. Finally, consider $w_0 = w_{\Pi,0}$ and $w_1 = w_\Pi$ and observe that

$$\|\Pi\|^2 \geq \frac{(w_{\Pi,0}, w_{\Pi,0})^2}{\|w_{\Pi,0}\|^4} + \frac{((\Pi - \mathcal{P})w_\Pi, (\Pi - \mathcal{P})w_\Pi)^2}{\|\Pi - \mathcal{P}\|^2} = 1 + \|\Pi - \mathcal{P}\|^2.$$

□

THE INSTITUTE FOR COMPUTATIONAL ENGINEERING AND SCIENCES (ICES), THE UNIVERSITY OF TEXAS AT AUSTIN, 201 E 24TH ST, AUSTIN, TX 78712, USA

*E-mail address*, corresponding author: brendan@ices.utexas.edu